\documentclass[a4paper,11pt,fleqn]{article}
\usepackage{amsmath,amssymb,amsthm,graphicx,subfigure,float,caption,epstopdf,tabularx,color, bm,amsfonts,epic}
\usepackage[top=1in, bottom=1in, left=1.25in, right=1.25in]{geometry}
\usepackage{appendix}
\usepackage{multirow}
\usepackage{diagbox}
\usepackage{appendix}
\newtheorem{shm}{Scheme}[section]
\allowdisplaybreaks
\newtheorem{thm}{Theorem}[section]
\newtheorem{shn}{Scheme}[section]
\newtheorem{lem}{Lemma}[section]

\newtheorem{rmk}{Remark}[section]
\newtheorem*{prf}{Proof}
\numberwithin{equation}{section}

\usepackage{indentfirst}
\graphicspath{{Fig}}

\begin{document}
\title{Arbitrary high-order linear structure-preserving schemes for the regularized long-wave
equation}
\author{Chaolong Jiang$^{1,2}$, Xu Qian$^1$\footnote{Correspondence author. Email:
qianxu@nudt.edu.cn.},  Songhe Song$^1$, and Jin Cui$^{3,4}$ \\
{\small $^1$ Department of Mathematics, College of Liberal Arts and Science, }\\
{\small National University of Defense Technology, Changsha, 410073, P.R. China}\\
{\small $^2$ School of Statistics and Mathematics, }\\
{\small Yunnan University of Finance and Economics, Kunming 650221, P.R. China}\\
{\small $^3$ Department of Basic Sciences, Nanjing Vocational College of Information Technology,}\\
  {\small Nanjing 210023, China} \\
  {\small $^4$ Jiangsu Key Laboratory for Numerical Simulation of Large Scale Complex Systems,}\\
{\small School of Mathematical Sciences,  Nanjing Normal University,}\\
{\small  Nanjing 210023, P.R. China}\\
}
\date{}
\maketitle

\begin{abstract}
In this paper, a class of arbitrarily high-order linear momentum-preserving and energy-preserving schemes are proposed, respectively, for solving the regularized long-wave
equation. For the momentum-preserving scheme, the key idea is based on the extrapolation/prediction-correction technique and the symplectic Runge-Kutta method in time, together with the standard Fourier pseudo-spectral method in space. We show that the scheme is linear, high-order, unconditionally stable and preserves the discrete momentum of the system. For the energy-preserving scheme, it is mainly based on the energy quadratization approach and the analogous linearized strategy used in the construction of the linear momentum-preserving scheme. The proposed scheme is linear, high-order and can preserve a discrete quadratic energy exactly. Numerical results are addressed to demonstrate the accuracy and efficiency of the proposed scheme.  \\[2ex]
\textbf{AMS subject classification:} 65M06, 65M70\\[2ex]
\textbf{Keywords:} Momentum-preserving, energy-preserving, high-order, linear scheme, regularized long-wave
equation.
\end{abstract}

\section{Introduction}
 In 1966, to describe long-wave behaviour, Peregrine \cite{Peregrine66jfm} first proposed the regularized long-wave (RLW) equation. Later on,  Benjamin, Bona and Mahony \cite{BBM-72-PTRSL} used the RLW equation to model small amplitude long waves on the surface of water in a channel. Nowadays, { the RLW equation has been widely applied} in fluids, plasma, anharmonic lattice and low temperature nonlinear crystals, etc. In this paper, we consider the two-dimensional RLW equation, as follows:
\begin{align}\label{RLW-equation}
\left \{
 \aligned
&\partial_tu({\bf x},t)+\alpha\partial_x u({\bf x},t)+\beta\partial_yu({\bf x},t)+{\color{blue} \alpha} u({\bf x},t)\partial_xu({\bf x},t)+{\color{blue} \beta} u({\bf x},t)\partial_yu({\bf x},t)\\
&~~~~~~~~~~~~~~~~~~~~~~~~~-\mu\partial_{xxt}u({\bf x},t)-\theta\partial_{yyt}u({\bf x},t)=0,\ {\bf x}\in\Omega\subset\mathbb{R}^2,\ t>0,\\%\times \mathbb{R},\ d=1,2
&u({\bf x},0)=u_0({\bf x}),\ {\bf x}\in{\Omega}\subset\mathbb{R}^2,
\endaligned
 \right.
  \end{align}
where $t$ is the time variable, ${\bf x}$ is the spatial variable, $u:=u({\bf x},t)$ is the real-valued wave function, $\alpha,\beta,\mu>0$ and $\theta>0$ are given real constants, $u_0({\bf x})$ is a given initial data, $\Omega=[a,b]\times[c,d]$ is a bounded domain, and the periodic boundary condition is considered on $\partial\Omega$.   %We assume $u(x,y,t)\to 0$ as $x\to \pm\infty$ and $y\to \pm\infty$.

According to \cite{Sun-Qin-2004-jcm}, when $u$ is assumed to be smooth, the RLW equation \eqref{RLW-equation} then can be rewritten as the classical Hamiltonian formulation
\begin{align}\label{RLW-Hamiltonian-system}
u_t=\mathcal S\frac{\delta\mathcal H}{\delta u},
\end{align}
where $\mathcal S = -(1-\mu\partial_{xx}-\theta\partial_{yy})^{-1}(\alpha\partial_x+\beta\partial_y)$ is a Hamiltonian operator and the Hamiltonian functional is
\begin{align}\label{RLW-Hamiltonian-energy}
\mathcal H(t) =\int_{\Omega}\Big(\frac{1}{2}u^2+\frac{1}{6} u^3\Big)d{\bf x},\ t\ge0.
\end{align}
Besides the Hamiltonian energy \eqref{RLW-Hamiltonian-energy}, {the RLW equation \eqref{RLW-equation} also conserves the mass} %\cite{Olver-79-MPPS}
\begin{align}\label{RLW-mass-energy}
\mathcal M(t)=\int_{\Omega}ud{\bf x}\equiv \mathcal M(0),\ t\ge0,
\end{align}
 and momentum
\begin{align}\label{RLW-momentum}
\mathcal I(t) =\int_{\Omega}\Big(\frac{1}{2}u^2+\frac{\mu}{2} u_x^2+\frac{\theta}{2} u_y^2\Big)d{\bf x}\equiv \mathcal I(0),\ t\ge0.
\end{align}

The first attempts to investigate the numerical solution of the RLW equation can date back to 1966, when Peregrine \cite{Peregrine66jfm} constructed a finite difference scheme to
find the development of an undular bore for the equation. %In 1979,  %Alexander and Morris \cite{Alexander-Morris-79-jcp} proposed a dissipative Galerkin method  and some comparative numerical comparisons between the RLW equation and the KdV equation are presented.
Later on, various efficient and accurate numerical schemes were proposed, including the finite element methods \cite{AM-79-jcp,LHQjsc2018,MC-cpc-2012}, the pseudo-spectral method \cite{GCjcp-1988}, the meshless method \cite{DS-CPC-2011}, and so on. However, these schemes do not exactly preserve the momentum \eqref{RLW-momentum} nor the Hamiltonian energy \eqref{RLW-Hamiltonian-energy} of the RLW equation. It has been shown that, in the numerical simulation of  solitary wave interactions for the RLW equation, the schemes that cannot preserve invariant quantities lead to an amplitude error which
grows with time and a quadratic perturbation in the phase \cite{DL-jpamg-2003}. Thus, developing numerical
schemes which can preserve the invariant quantities is necessary for accurately simulating the solitary wave of the RLW equation.

% multisymplectic method \cite{Cai-2009-cpc,Li-Sun-2013-mcm,Sun-Qin-2004-jcm}

 Over the years, many attempts to develop momentum-preserving schemes have been done on the RLW equation. In \cite{CWG-1991-jcp}, Chang, et al. proposed a fully-implicit finite difference scheme which can satisfy the discrete analogue of the momentum. Meanwhile, it is proven rigorously in mathematics that the scheme is second-order accurate in time and space. { More recently, based on the technique
in \cite{Frasca-Caccia-IMA-2020}, Frasca-Caccia et al. \cite{Frasca-Caccia-AMC-2021} proposed a novel class of fully-implicit second-order schemes in time, which can preserve both mass and momentum. For the implicit scheme, a nonlinear equation shall be solved by using a nonlinear iterative method at each time step, and thus it may be time consuming.}  In Ref. \cite{Zhang-amc-2005}, Zhang constructed a linear momentum-preserving scheme, in which a linear system is to be solved at every time step. Thus it is computationally much cheaper than that of the implicit one. Other linear momentum-preserving schemes can be found in \cite{HWG-DCDSB-2019,HWG-NMPDE-2020,JZS-anm-2020,KF-2009-JJIAM,WZC-07-amc,WDG-amc-2019}. Nevertheless, to our best knowledge, all of the existing linear momentum-preserving schemes can achieve at most second-order accurate in time.

 The first aim of this paper is to present a novel paradigm for developing arbitrary high-order linear momentum-preserving (LMP) schemes for the RLW equation. We firstly reformulate the original equation into a skew-adjoint system. The nonlinear terms of the system are then approximated with an extrapolation/prediction-correction technique to obtain a linearized momentum-conserving system. We further apply the symplectic Runge-Kutta (RK) method in time and the Fourier pseudo-spectral method in space to the resulting linearized system and a fully discrete scheme is obtained, which is linear, high-order, and can preserve the discrete momentum exactly.

 On the other hand, Hamiltonian energy is one of most famous geometric characteristics of the Hamiltonian systems. In the past few decades, there has been a series of studies on energy-preserving schemes for Hamiltonian systems, such as, but are not limited to, averaged vector field (AVF) methods \cite{LWQ14,QM08}, Hamiltonian Boundary Value Methods (HBVMs) \cite{BI16,BIT10} and energy-preserving collocation methods \cite{CH11bit,H10}, functionally fitted energy-preserving methods\cite{LWsina16,Miyatake-2014-bit} and energy-preserving continuous stage Runge-Kutta methods \cite{MB16,TS12}. These schemes can be easily extended to propose high-order energy-preserving schemes for the RLW equation, but the resulting schemes are fully implicit in which a large of nonlinear systems need to be solved at every time step and thus it may be very time consuming. Based on the discrete variational derivative method \cite{DO11,Furihata99}, Koide and Furihata \cite{KF-2009-JJIAM}  constructed a linear energy-preserving scheme, however { the scheme is only second-order accurate in time.}

The second aim of this paper is to present a new  strategy for proposing arbitrary high-order linear energy-preserving (LEP) schemes for the RLW equation. Based on the energy quadratization (EQ) approach \cite{GZYW18,YZW17,ZYGW17} and the analogous linearized idea used in the construction of the linear momentum-preserving scheme, we obtain a class of linear energy-preserving schemes. The proposed scheme is high-order, easy to be implemented and can exactly preserve the discrete quadratic energy.

The rest of this paper is organized as follows. In Section \ref{RLW-section-2}, high-order linear structure-preserving schemes for the RLW equation are proposed and their structure-preserving properties are investigated in details. Extensive numerical examples for the one dimensional (1D) and two dimensional (2D) RLW equation are shown to illustrate the accuracy and efficiency of the scheme, respectively. We draw some conclusions in Section \ref{RLW-section-4}.

\section{High-order linear structure-preserving scheme}\label{RLW-section-2}
Denote $t_{n}=n\tau,$ and $t_{ni}=t_n+c_i\tau,\ i=1,2,\cdots,s$, $n = 0,1,2,\cdots,$ where $\tau$ is the time step and $c_1,c_2,\cdots,c_s$ are distinct real numerbers (usually $0\le c_i\le 1$). The approximation of the function $u({\bf x},t)$ at points  $t_n$ and $t_{ni}$ are denoted by $u^n$ and $u_{i}^n$, respectively. Similarly, the approximations of the function $q({\bf x},t)$ at points  $t_n$ and $t_{ni}$ are denoted by $q^n$ and $q_{i}^n$, respectively.

\subsection{High-order linear momentum-preserving scheme}\label{RLW-section-2.1}
In this subsection, we propose a class of temporal high-order linear semi-discrete momentum-preserving schemes based on the combination of the extrapolation/prediction-correction approach with the symplectic RK method. To this end, we first rewrite \eqref{RLW-equation} into an equivalent form
\begin{align}\label{RLW-momentum-equation}
\mathcal Du_t=\mathcal G(u)u,\ \mathcal G(u)=-\big[\alpha\partial_x+\beta\partial_y
+\frac{{ \alpha}}{3}(u\partial_x+\partial_xu)+\frac{{\beta}}{3}(u\partial_y+\partial_yu)\big],
\end{align}
where $\mathcal D =1-\mu\partial_{xx}-\theta\partial_{yy}$ is a self-adjoint operator.%, the operators $u\partial_x+\partial_xu$ and $u\partial_y+\partial_yu$ satisfy $(u\partial_x+\partial_xu)u=u\partial_xu+\partial_x(u^2)$ and $(u\partial_y+\partial_yu)u=u\partial_yu+\partial_y(u^2)$, respectively.

We assume that the numerical solutions of $u$ in \eqref{RLW-momentum-equation} from $t=0$ to $t\le t_n$ are obtained and the Lagrange interpolating
polynomial (denoted by $\mathcal L_{MP}({\bf x},t)$) based on these solutions and the corresponding grid points in time are proposed.
Then we use $\mathcal L_{MP}({\bf x},t)$ to approximate the nonlinear terms of \eqref{RLW-momentum-equation} in $(t_n,t_{n+1}]$ and a linearized system is obtained, as follows:
\begin{align}\label{linear-RLW-momentum-equation}
\mathcal D\tilde u_t=\mathcal G(\mathcal L_{MP}({\bf x},t))\tilde u,
\end{align}
where $\tilde u:=\tilde u({\bf x},t)$ is an approximation of $u({\bf x},t)$ in $(t_n,t_{n+1}]$.

Here, wo should note that since $\mathcal G(u)$ is a skew-adjoint operator for all $u$, $\mathcal G(\mathcal L_{MP}({\bf x},t))$ is also skew-adjoint and the following result holds
\begin{align}
\frac{d}{dt}(\frac{1}{2}\mathcal D \tilde u,\tilde u)=(\mathcal D \tilde u_t,\tilde u)=(\mathcal G(\mathcal L_{MP}({\bf x},t))\tilde u,\tilde u)=0.
\end{align}
{However, we should note that the system \eqref{linear-RLW-momentum-equation} cannot conserve the mass. In fact, with noting the self-adjoint property of the operator $\mathcal D$ and the periodic boundary condition, we can deduce from \eqref{RLW-momentum-equation} that
\begin{align}
\frac{d}{dt}(\tilde u,1)&=(\tilde u_t,1)\nonumber\\
&=(\mathcal D^{-1}\mathcal G(\mathcal L_{MP}({\bf x},t))\tilde u,1)\nonumber\\
&=(\mathcal G(\mathcal L_{MP}({\bf x},t))\tilde u,\mathcal D^{-1}1)\nonumber\\
&=(\mathcal G(\mathcal L_{MP}({\bf x},t))\tilde u,1)\nonumber\\
&=-\frac{\alpha}{3}(\mathcal L_{MP}({\bf x},t)\partial_x\tilde u,1)-\frac{\beta}{3}(\mathcal L_{MP}({\bf x},t)\partial_y\tilde u,1)\nonumber\\
&\not=0,
\end{align}
where $(\cdot,\cdot)$ is the continuous inner production defined by $(u,v)=\int_{\Omega}uvd{\bf x}$ and $\|u\|^2=(u,u)$.}

Finally, we apply an RK method to the linearized system \eqref{linear-RLW-momentum-equation} in time to give a class of semi-discrete linear schemes for \eqref{RLW-momentum-equation}, as follows:
 \begin{shn}\label{RLW-scheme-2.1} Let $b_i,a_{ij}(i,j=1,\cdots,s)$ be real numbers and let $c_i=\sum_{j=1}^sa_{ij}$. For the given $u^n$ and $u_{i}^{n,*}=\mathcal L_{MP}({\bf x},t_n+c_i\tau),i=1,\cdots,s$, an s-stage linear Runge-Kutta method is given by
\begin{align}\label{Linear-momentum-schemes}
\left\{
\begin{aligned}
&k_i^n=\mathcal D^{-1}\mathcal G(u_{i}^{n,*})u_{i}^n,\ u_{i}^n=u^n+\tau\sum_{j=1}^{s}a_{ij}k_j^n,\\
&u^{n+1}=u^n+\tau\sum_{i=1}^{s}b_ik_i^n,
\end{aligned}
\right.
\end{align}
where $i=1,2,\cdots,s$.
\end{shn}

\begin{thm} \label{RLW-equation-th2.1}  If the coefficients of the RK method satisfy
\begin{align}\label{RLW-RK-coeff}
b_ia_{ij}+b_ja_{ji}=b_ib_j, \ \forall \ i,j=1,\cdots,s,
\end{align}
 {\bf Scheme 2.1} exactly preserves the semi-discrete momentum, as follows:
\begin{align}\label{RLW-equation-new-2.5}
 I^{n+1}=I^n, \ I^n=\frac{1}{2}(u^n,\mathcal Du^n),\  n=0,1,2,\cdots.
\end{align}
\end{thm}
\begin{prf} It follows from the second equality of \eqref{Linear-momentum-schemes} that
\begin{align}\label{LMP-RLW-equation-2.4}
 I^{n+1}-{I}^n&=\frac{1}{2}(u^{n+1},\mathcal Du^{n+1})-\frac{1}{2}(u^{n},\mathcal Du^{n})\nonumber\\
&=\frac{1}{2}(u^n+\tau\sum_{i=1}^{s}b_ik_i^n,\mathcal D(u^n+\tau\sum_{j=1}^{s}b_jk_j^n))-\frac{1}{2}(u^{n},\mathcal Du^{n})\nonumber\\
&=\frac{\tau}{2}\sum_{i=1}^{s}b_i(u^n,\mathcal Dk_i^n)+\frac{\tau}{2}\sum_{i=1}^{s}b_i(k_i^n,\mathcal Du^n)+\frac{\tau^2}{2}\sum_{i,j=1}^{s}b_ib_j(k_i^n,\mathcal Dk_j^n).
  \end{align}
  With noting
  \begin{align}\label{LMP-RLW-equation-2.5}
  \frac{\tau}{2}\sum_{i=1}^{s}b_i(u^n,\mathcal Dk_i^n)&=\frac{\tau}{2}\sum_{i=1}^{s}b_i(u_{i}^n-\tau\sum_{j=1}^{s}a_{ij}k_j^n,\mathcal Dk_i^n)\nonumber\\
  &=\frac{\tau}{2}\sum_{i=1}^{s}b_i(u_{i}^n,\mathcal Dk_i^n)-\frac{\tau^2}{2}\sum_{i,j=1}^{s}b_ja_{ji}(k_i^n,\mathcal Dk_j^n).
  \end{align}
  Similarly, we have
  \begin{align}\label{RLW-equation-2.6}
  \frac{\tau}{2}\sum_{i=1}^{s}b_i(k_i^n,\mathcal Du^n)=\frac{\tau}{2}\sum_{i=1}^{s}b_i(k_i^n,\mathcal Du_{i}^n)-\frac{\tau^2}{2}\sum_{i,j=1}^{s}b_ia_{ij}(k_i^n,\mathcal Dk_j^n).
  \end{align}
  {We insert \eqref{LMP-RLW-equation-2.5} and \eqref{RLW-equation-2.6} into \eqref{LMP-RLW-equation-2.4} and then use the self-adjoint property of $\mathcal D$ to obtain}
  \begin{align*}
 I^{n+1}-{I}^n&=\frac{\tau}{2}\sum_{i=1}^{s}b_i\big[(u_{i}^n,\mathcal Dk_i^n)+(k_i^n,\mathcal Du_{i}^n)\big]+\frac{\tau^2}{2}\sum_{i,j=1}^{s}(b_ib_j-b_ia_{ij}-b_ja_{ji})(k_i^n,\mathcal Dk_j^n)\\
  %&=\frac{\tau}{2}\sum_{i=1}^{s}b_i\big[(u_{ni},\mathcal Dk_i)+(k_i,\mathcal Du_{ni})\big]\\
  &=\tau\sum_{i=1}^{s}b_i(u_{i}^n,\mathcal Dk_i^n)+\frac{\tau^2}{2}\sum_{i,j=1}^{s}(b_ib_j-b_ia_{ij}-b_ja_{ji})(k_i^n,\mathcal Dk_j^n).
  %&=0,
  \end{align*}
  {The condition \eqref{RLW-RK-coeff} together with the equality}
\begin{align*}
(u_{i}^n,\mathcal Dk_i^n)=(u_{i}^n,\mathcal G(u_{i}^{n,*})u_{i}^n)=0
\end{align*}
 implies $\ I^{n+1}={I}^n$. This completes the proof. \qed
\end{prf}
\begin{rmk}\label{RLW-rmk-2.1} Assume that the initial condition $u_0({\bf x})$ is sufficiently smooth, then it follows from \eqref{RLW-equation-new-2.5} that the numerical solution of {\bf Scheme 2.1} satisfies
\begin{align*}
\sqrt{\|u^n\|^2+\mu\|\partial_xu^n\|^2+\theta\|\partial_yu^n\|^2}=\sqrt{\|u_0({\bf x})\|^2+\mu\|\partial_xu_0({\bf x})\|^2+\theta\|\partial_yu_0({\bf x})\|^2}\le C,
\end{align*}
which implies that
\begin{align*}
\|u^n\|\le C,\ \|\partial_xu^n\|\le C,\ \|\partial_yu^n\|\le C,
\end{align*}
is uniformly bounded. Thus, {\bf Scheme 2.1} is unconditionally stable.% with respect to the initial value.
\end{rmk}

To improve accuracy as well as stability of  {\bf Scheme} \ref{RLW-scheme-2.1}, the prediction-correction strategy \cite{GZW2020jcp,GZW-20-SISC,SX-CiCP-2018} is employed to {\bf Scheme} \ref{RLW-scheme-2.1} and an $s$-stage linear Runge-Kutta prediction-correction schemes are obtained, as follows:
   \begin{shn}\label{RLW-scheme-2.2} Let $b_i,a_{ij}(i,j=1,\cdots,s)$ be real numbers and let $c_i=\sum_{j=1}^sa_{ij}$. For the given $u^n$, an s-stage linear Runge-Kutta prediction-correction method is given by the following two steps:

   1. Prediction: {set the iterative initial value $k_i^{n,0}=u^n$} and let $M>0$ be a given integer, for $m=0,1,\cdots,M-1$, we calculate $k_i^{n,m+1}$ using
{\begin{align}\label{Linear-momentum-P-schemes}
 k_i^{n,m+1}=\mathcal D^{-1}\mathcal G(u_{i}^{n,m})u_{i}^{n,m},\ u_{i}^{n,m}=u^n+\tau\sum_{j=1}^{s}a_{ij}k_j^{n,m},\ i=1,2,\cdots,s.
\end{align}}
Then, we set $u_{i}^{n,*}=u_{i}^{n,M}$.

2.  Correction: for the predicted value $u_{i}^{n,*},\ i=1,2,\cdots,s$, we obtain $u^{n+1}$ by
\begin{align}\label{Linear-momentum-C-schemes}
\left\{
\begin{aligned}
&k_i^n=\mathcal D^{-1}\mathcal G(u_{i}^{n,*})u_{i}^n,\ u_{i}^n=u^n+\tau\sum_{j=1}^{s}a_{ij}k_j^n,\\
%& k_i^n=-(1-\mu\partial_{xx}-\theta\partial_{yy})^{-1}\big[(\alpha\partial_x+\beta\partial_y)u_{i}^n+\frac{\gamma}{3}(u_{i}^{n,*}\partial_x+\partial_xu_{i}^{n,*})u_{i}^n\\
%&~~~~~~~+\frac{\delta}{3}(u_{i}^{n,*}\partial_y+\partial_yu_{i}^{n,*})u_{i}^n\big],\\
&u^{n+1}=u^n+\tau\sum_{i=1}^{s}b_ik_i^n,
\end{aligned}
\right.
\end{align}
where $i=1,2,\cdots,s$.
\end{shn}
{\begin{rmk} It's worth noting that when the termination criterion for iteration is satisfied in \eqref{Linear-momentum-P-schemes}, we obtain the value $k_i^{n,M},\ i=1,2,\cdots,s$. Then,  $u_{i}^{n,M}$ is updated by
\begin{align*}
u_{i}^{n,M}=u^n+\tau\sum_{j=1}^{s}a_{ij}k_j^{n,M},\ i=1,2,\cdots,s.
\end{align*}
\end{rmk}}
\begin{thm} \label{RLW-equation-th2.2}  If the coefficients of {\color{blue} \eqref{Linear-momentum-C-schemes}} satisfy the condition \eqref{RLW-RK-coeff},
 {\bf Scheme 2.2} exactly conserves the semi-discrete momentum \eqref{RLW-equation-new-2.5}.
\end{thm}

\begin{prf} The aim of the prediction step is to obtain a solution $u_{i}^{n,*},\ i=1,2,\cdots,s$ with the desirable accuracy as well as stability, which is then used to linearize the nonlinear system \eqref{RLW-momentum-equation}. Moreover, the proof of the semi-discrete momentum \eqref{RLW-equation-new-2.5} only depends on the skew-adjoint property of $\mathcal G(u_{i}^{n,*})$. Thus, it is similar to Theorem \ref{RLW-equation-th2.1}. For brevity, we omit it.\qed
\end{prf}

\begin{rmk} By the similar argument to Remark \ref{RLW-rmk-2.1}, we can show that {\bf Scheme 2.2} is unconditionally stable.
\end{rmk}
\begin{rmk}\label{RLW-RMK-2.4} {If we take $c_1,c_2,\cdots,c_s$ as the zeros of the $s$th shifted Legendre polynomial
\begin{align*}
\frac{d^s}{dx^s}\Big(x^s(x-1)^s\Big),
\end{align*}
the RK (or collocation) method based on these nodes has the order $2s$ and satisfies  the condition \eqref{RLW-RK-coeff} (see Refs. \cite{ELW06,Sanzs88,SCbook94} and references therein). The RK coefficients for $s=2$ and $s=3$ (denoted by 4th-Gauss method and 6th-Gauss method, respectively) are  given in Table \ref{Gauss-cllocation-method}} (e.g., see Ref. \cite{ELW06})

\begin{table}[H]
\centering
%\begin{tabular}{c|cc}
%${c}$ & ${A}$  \\
%\hline
%& ${b}^{T}$ \\
%\end{tabular}
%=
\begin{tabular}{c|cc}
$\frac{1}{2}-\frac{\sqrt{3}}{6}$ &$\frac{1}{4}$ & $\frac{1}{4}- \frac{\sqrt{3}}{6}$\\
$\frac{1}{2}+\frac{\sqrt{3}}{6}$ &$\frac{1}{4}+ \frac{\sqrt{3}}{6}$ &$\frac{1}{4}$ \\
\hline
                                 &$\frac{1}{2}$&$\frac{1}{2}$
\end{tabular}\\
\vspace{1mm}
\begin{tabular}{c|ccc}
$\frac{1}{2}-\frac{\sqrt{15}}{10}$ &$\frac{5}{36}$ &  $\frac{2}{9}-\frac{\sqrt{15}}{15}$    &$\frac{5}{36}- \frac{\sqrt{15}}{30}$\\
$\frac{1}{2}$   &$\frac{5}{36}+ \frac{\sqrt{15}}{24}$  &  $\frac{2}{9}$  & $\frac{5}{36}- \frac{\sqrt{15}}{24}$ \\
$\frac{1}{2}+\frac{\sqrt{15}}{10}$ & $\frac{5}{36}+ \frac{\sqrt{15}}{30}$ & $\frac{2}{9}+\frac{\sqrt{15}}{15}$  & $\frac{5}{36}$ \\
\hline
& $\frac{5}{18}$ & $\frac{4}{9}$  & $\frac{5}{18}$
\end{tabular}
\caption{Gauss methods of order 4 and 6.}\label{Gauss-cllocation-method}
\end{table}
\end{rmk}
\begin{rmk}\label{RLW-RMK-2.5} As is well-known, the interpolation polynomial will be highly oscillating if too
many interpolation points are chosen, which makes the extrapolation
not sufficiently accurate. Inspired by Refs. \cite{GZW2020jcp,LS2020jsc}, we only choose $t_{n-1}$, $t_{n-1}+c_i\tau,\ i=1,2,\cdots,s$ and $t_n$ as the interpolation points in this paper. For example, for the 6th-Gauss method, we set the interpolation points $(t_{n-1},u^{n-1}),\ (t_{n-1}+c_1\tau,u_1^{n-1}),\ (t_{n-1}+c_2\tau,u_2^{n-1})$ and $(t_{n-1}+c_3\tau,u_3^{n-1})$, the Lagrange interpolating polynomial at $t_n+c_i\tau,\ i=1,2,3$ is given \cite{BWC2020arixv,LGZJSC-2021}, as follows:
\begin{align}\label{RLW-LI-2.12}
&u_1^{n,*}=(6\sqrt{15}-26)u^{n-1}+(-5\sqrt{15}/3+11)u_1^{n-1}+(16\sqrt{15}/3-24)u_2^{n-1}\nonumber\\
&~~~~~~~~+(-29\sqrt{15}/3+40)u_3^{n-1},\\\label{RLW-LI-2.13}
&u_2^{n,*}=-17u^{n-1}+(5\sqrt{15}/2+\frac{35}{2})u_1^{n-1}-17u_2^{n-1}+(-5\sqrt{15}/2+\frac{35}{2})u_3^{n-1},\\\label{RLW-LI-2.14}
&u_3^{n,*}=(-6\sqrt{15}-26)u^{n-1}+(29\sqrt{15}/3+40)u_1^{n-1}+(-16\sqrt{15}/3-24)u_2^{n-1}\nonumber\\
&~~~~~~~~+(5\sqrt{15}/3+11)u_3^{n-1}.
\end{align}
In particular, when the 6th-Gauss method and the above $u_i^{n,*},\ i=1,2,3$ are used for {\bf scheme 2.1}, the resulting scheme may achieve fourth-order accuracy in time \cite{LS2020jsc} and we denote it as 4th-LMPS. In addition, 4th-LMPS is a linear scheme, thus we obtain the starting values $u_1^1,\ u_2^1$ and $u_3^1$ by using the 6th-Gauss method for \eqref{RLW-momentum-equation}.
\end{rmk}

\begin{rmk} { Here, we introduce two notations, as follows:
 \begin{itemize}
 \item 4th-LMP-PCS: using 4th-Gauss method to {\bf Scheme \ref{RLW-scheme-2.2}};
 \item 6th-LMP-PCS: applying 6th-Gauss method to {\bf Scheme \ref{RLW-scheme-2.2}}.
 \end{itemize}
 In addition, we should note that, to the best of our knowledge, the theoretical result on the choice of iteration step $M$ is missing. From our numerical experience, 4th-LMP-PCS can achieve fourth-order accurate in time by choosing $M= 3$, and 6th-LEP-PCS can arrive at sixth-order accurate in time by choosing $M=5$.}
\end{rmk}

%Meanwhile, the interpolating
%polynomial of $u_N({\bf x},t)$ of order 4 in \cite{BWC2020arixv} is employed for 4th-LMPS and 4th-LEPS.

\subsection{High-order linear energy-preserving scheme}\label{RLW-section-2.2}

In this section, we propose a class of temporal high-order linear semi-discrete energy-preserving schemes for the RLW equation \eqref{RLW-equation} by using the EQ idea \cite{GZYW18,JWG19,JCW19jsc,YZW17,ZYGW17} and the similar linearized strategy mentioned as above. To this end, we start with introducing an auxiliary variable
\begin{align*}
q:=q({\bf x},t)=u^2.
\end{align*}
Then, the Hamiltonian energy \eqref{RLW-Hamiltonian-energy} is rewritten into the following quadratic form
\begin{align}\label{IEQ-RLW-energy}
\mathcal E(t) =\int_{\Omega}\Big(\frac{1}{2}u^2+\frac{1}{6} uq\Big)d{\bf x},\ t\ge0.
\end{align}
Based on the energy variational principle, we obtain the following EQ reformulated system from \eqref{RLW-Hamiltonian-system}
\begin{align}\label{RLW-IEQ-reformulation1}
\left\{
\begin{aligned}
&\partial_t u=\mathcal S\bigg (u+\frac{1}{6}q+\frac{1}{3}u^2\bigg), \\
&\partial_t q=2u\partial_tu,
\end{aligned}
\right.
\end{align}
 {with the consistent initial conditions
\begin{align}
u({\bf x},0)=u_0({\bf x}),\ q({\bf x},0)=(u_0({\bf x}))^2.
\end{align}}
~~~~{We then define the operator $\mathcal L$ as $\mathcal Lw \overset{Def}{=}(\mathcal S u)w\overset{Def}{=}\mathcal S\big(u\cdot w\big)$, and the adjoint operator of
$\mathcal L$ is denoted as $\mathcal L^*$, which satisfies $\mathcal L^*w=-u\cdot \mathcal Sw$. Denoting
\begin{equation*}
\Phi=\begin{bmatrix}
u\\
q
\end{bmatrix},
\end{equation*}
we obtain
\begin{align}
\frac{\delta\mathcal E}{\delta\Phi}=\begin{bmatrix}\vspace{2mm}
\frac{\delta\mathcal E}{\delta u}\\
\frac{\delta\mathcal E}{\delta q}
\end{bmatrix}:=\begin{bmatrix}\vspace{2mm}
u+\frac{1}{6}q\\
\frac{1}{6}u
\end{bmatrix}.
\end{align}
Then we rewrite \eqref{RLW-IEQ-reformulation1} into
\begin{align}\label{RLW-IEQ-reformulation-1}
\left\{
\begin{aligned}
&\partial_t u=\mathcal S\frac{\delta\mathcal E}{\delta u}+2\mathcal L\frac{\delta\mathcal E}{\delta q}, \\
&\partial_t q=-2\mathcal L^*\frac{\delta\mathcal E}{\delta u}+4u\cdot \mathcal L\frac{\delta\mathcal E}{\delta q},
\end{aligned}
\right.
\end{align}
which is equivalent to
\begin{align}\label{RLW-IEQ-reformulation-2}
\left\{
\begin{aligned}
&\partial_t u=\mathcal S\frac{\delta\mathcal E}{\delta u}+2\mathcal L\frac{\delta\mathcal E}{\delta q}, \\
&\partial_t q=-2\mathcal L^*\frac{\delta\mathcal E}{\delta u}-4\mathcal L^{*}\mathcal S^{-1}\mathcal L\frac{\delta\mathcal E}{\delta q}.
\end{aligned}
\right.
\end{align}
Thus,}
%The system \eqref{RLW-IEQ-reformulation1} is equivalent to
%\begin{align}\label{RLW-IEQ-reformulation2}
%\left\{
%\begin{aligned}
%&\partial_t u=\mathcal S\frac{\delta\mathcal E}{\delta u}+2\mathcal L\frac{\delta\mathcal E}{\delta q}, \\
%&\partial_t q=-2\mathcal L^*\frac{\delta\mathcal E}{\delta u}+4\mathcal L^*\mathcal S^{-1} \mathcal L\frac{\delta\mathcal E}{\delta q},
%\end{aligned}
%\right.
%\end{align}
{the system \eqref{RLW-IEQ-reformulation1} can be reformulated into an energy-conserving system, as follows:
\begin{equation}\label{IEQ-skew-system}
\Phi_t=\mathcal M(\Phi) \frac{\delta \mathcal E}{\delta \Phi},
\end{equation}
where $\mathcal M(\Phi)$ is the following skew-adjoint operator
\begin{equation*}
\mathcal M(\Phi) =\begin{bmatrix}
\mathcal S &2\mathcal L\\
-2\mathcal L^* &-4\mathcal L^{*}\mathcal S^{-1} \mathcal L
\end{bmatrix}.
\end{equation*}
}
%with the consistent initial condition
%\begin{align}\label{eq2.6}
%u({\bf x},0)=u_0({\bf x}),  \  q({\bf x},0)=u_0({\bf x}).
%\end{align}
 \begin{thm}\label{LI-E-NLS-thm-2.1} The IEQ reformulation \eqref{RLW-IEQ-reformulation1} preserves the following quadratic energy
  \begin{align}\label{RLW-Q-energy}
  \frac{d}{dt} \mathcal{E}(t)=0,\ \mathcal{E}(t)=\int_{\Omega}\Big(\frac{1}{2}u^2+\frac{1}{6} uq\Big)d{\bf x},\ t\ge0.
  \end{align}
  \end{thm}
  \begin{prf} According to \eqref{IEQ-RLW-energy} and \eqref{RLW-IEQ-reformulation1} together with the skew-adjoint property of $\mathcal S$, we have
  \begin{align*}
  \frac{d}{dt} \mathcal{E}(t)&=\int_{\Omega}\Big(u\partial_tu+\frac{1}{6} q\partial_tu+\frac{1}{6} u\partial_tq\Big)d{\bf x}\\
  &=\int_{\Omega}\Big[(u+\frac{1}{6} q+\frac{1}{3} u^2)\partial_tu\Big]d{\bf x}\\
  &=\int_{\Omega}\big(u+\frac{1}{6} q+\frac{1}{3} u^2\big)\mathcal S\big (u+\frac{1}{6}q+\frac{1}{3}u^2\big)d{\bf x}\\
  &=0.
  \end{align*}
  This completes the proof. \qed
  \end{prf}

Similarly, we assume that the numerical solutions of $u$ in \eqref{RLW-IEQ-reformulation1} from $t=0$ to $t\le t_n$ are obtained and the Lagrange interpolating
polynomial (denoted by $\mathcal L_{EP}({\bf x},t)$) based on these solutions and the corresponding grid points in time are proposed.
Then we apply the Lagrange interpolating
polynomial $\mathcal L_{EP}({\bf x},t)$ to the nonlinear terms of \eqref{RLW-IEQ-reformulation1} and a linearized energy-conserving system in $(t_n,t_{n+1}]$ is obtained, as follows:
\begin{align}\label{linear-RLW-IEQ-equation}
\left\{
\begin{aligned}
&\partial_t \tilde u=\mathcal S\bigg (\tilde u+\frac{1}{6}\tilde q+\frac{1}{3}\mathcal L_{EP}({\bf x},t)\tilde u\bigg), \\
&\partial_t \tilde q=2\mathcal L_{EP}({\bf x},t)\partial_t\tilde u,\ i=1,2,\cdots,s,
\end{aligned}
\right.
\end{align}
where $\tilde u:=\tilde u({\bf x},t)$ and $\tilde q:=\tilde q({\bf x},t)$ {\color{blue} are} approximations of $u({\bf x},t)$ and $q({\bf x},t)$ in $(t_n,t_{n+1}]$, respectively. % Here, $u_{ni}^{*}$ is an (explicit) extrapolation approximation to $u({\bf x},t_n+c_i\tau),\  i=1,2,\cdots,s$ defined as above.

Next, we apply an RK method to the linearized system \eqref{linear-RLW-IEQ-equation} to give a class of semi-discrete linear RK methods for solving \eqref{RLW-IEQ-reformulation1}, as follows:
\begin{shn}\label{RLW-scheme-2.3} Let $b_i,a_{ij}\ (i,j=1,\cdots,s)$ be real numbers and let $c_i=\sum_{j=1}^sa_{ij}$. For the given $(u^n,q^n)$ and $u_{i}^{n,*}=\mathcal L_{EP}({\bf x},t_n+c_i\tau),\ i=1,\cdots,s$, an s-stage linear Runge-Kutta method is given by
\begin{align}\label{E-linear-RLW-IEQ-equation1}
\left\{
\begin{aligned}
&k_i^n=\mathcal S\bigg (u_{i}^n+\frac{1}{6}q_{i}^n+\frac{1}{3}u_{i}^{n,*} u_{i}^n\bigg),\ u_{i}^n=u^n+\tau\sum_{j=1}^{s}a_{ij}k_j^n,\\
 &q_{i}^n=q^n+\tau\sum_{j=1}^{s}a_{ij}l_j^n,\ l_i^n=2u_{i}^{n,*}k_i^n,\ i=1,2,\cdots,s,
\end{aligned}
\right.
\end{align}
and $(u^{n+1},q^{n+1})$ is then updated by
\begin{align}\label{E-linear-RLW-IEQ-equation2}
u^{n+1}=u^n+\tau\sum_{i=1}^{s}b_ik_i^n,\ q^{n+1}=q^n+\tau\sum_{i=1}^{s}b_{i}l_i^n.
\end{align}
\end{shn}
 \begin{thm}\label{RLW-equation-th2.3} If the coefficients of \eqref{E-linear-RLW-IEQ-equation1} and \eqref{E-linear-RLW-IEQ-equation2} satisfy \eqref{RLW-RK-coeff}, {\bf Scheme 2.3} preserves the following semi-discrete energy
  \begin{align}\label{semi-RLW-Q-energy}
   E^{n+1}= E^n,\ {E}^n=\int_{\Omega}\Big(\frac{1}{2}(u^n)^2+\frac{1}{6} u^nq^n\Big)d{\bf x},\ n=0,1,2,\cdots,.
  \end{align}
  \end{thm}
\begin{prf}According to \eqref{E-linear-RLW-IEQ-equation1} and \eqref{E-linear-RLW-IEQ-equation2}, we have
\begin{align*}
{E}^{n+1}- E^n&=\frac{1}{2}(u^{n+1},u^{n+1})+\frac{1}{6} (u^{n+1},q^{n+1})-\frac{1}{2}(u^{n},u^{n})-\frac{1}{6} (u^{n},q^{n})\nonumber\\
&=\frac{\tau}{2}\sum_{i=1}^{s}b_i((u^n,k_i^n)+(k_i^n,u^n))+\frac{\tau^2}{2}\sum_{i,j=1}^{s}b_ib_j(k_i^n,k_j^n)\nonumber\\
&~~~+\frac{\tau}{6}\sum_{i=1}^{s}b_i\big((u^n,l_i^n)+ (k_i^n,q^n)\big)+\frac{\tau^2}{6}\sum_{i,j=1}^{s}b_ib_j(k_i^n,l_j^n)\nonumber\\
&=\tau\sum_{i=1}^{s}b_i\big((u_i^n,k_i^n)+\frac{1}{6}(u_i^n,l_i^n)+\frac{1}{6}(k_i^n,q_i^n)\big)\nonumber\\
&+\frac{\tau^2}{2}\sum_{i,j=1}^{s}(b_ib_j-b_ia_{ij}-b_ja_{ji})(k_i^n,k_j^n)+\frac{\tau^2}{6}\sum_{i,j=1}^{s}(b_ib_j-b_ia_{ij}-b_ja_{ji})(k_i^n,l_j^n).
\end{align*}
   With the condition \eqref{RLW-RK-coeff} and the skew-adjoint property of $\mathcal S$, we have
    \begin{align*}
    (u_i^n,k_i^n)+\frac{1}{6}(u_i^n,l_i^n)+\frac{1}{6}(k_i^n,q_i^n)&=(k_i^n,u_i^n)+\frac{1}{6}(u_i^n,2u_{i}^{n,*}k_i^n)+\frac{1}{6}(k_i^n,q_i^n)\\
    &=(k_i^n,u_i^n)+\frac{1}{6}(k_i^n,2u_{i}^{n,*}u_i^n)+\frac{1}{6}(k_i^n,q_i^n)\\
    &=(k_i^n,u_i^n+\frac{1}{3}u_{i}^{n,*}u_i^n+\frac{1}{6}q_i^n)\\
    &{=(\mathcal S(u_i^n+\frac{1}{3}u_{i}^{n,*}u_i^n+\frac{1}{6}q_i^n),u_i^n+\frac{1}{3}u_{i}^{n,*}u_i^n+\frac{1}{6}q_i^n)}\\
    &=0,
    \end{align*}
which implies $ E^{n+1}={E}^n$.
    This completes the proof. \qed
   \end{prf}

%\section{The construction of the high-order linearly implicit exponential integrator}\label{Sec:LI-EP:3}
%In this section, a class of high-order linearly implicit energy-preserving integrators are proposed for the SAV reformulated system \eqref{NLS-SAV-reformulation}. For simplicity, in this paper, we shall introduce our schemes in two space dimension, i.e., $d=2$ in \eqref{NLS-SAV-reformulation}. Generalizations to $d=1$ or $3$ are straightforward.
%\begin{rmk}When the time step $\tau$ is sufficiently small, the  {\bf Scheme \ref{RLW-scheme-2.3}} admits unique solution. This is because the coefficient matrix of the linear equations \eqref{Linear-momentum-schemes} almost equates to unite matrix,
%
%\end{rmk}

Our next object is to apply  the prediction-correction strategy to {\bf Scheme \ref{RLW-scheme-2.3}} and a class of high-order linear RK prediction-correction schemes are given, as follows:
   \begin{shn}\label{RLW-scheme2.4} Let $b_i,a_{ij}(i,j=1,\cdots,s)$ be real numbers and let $c_i=\sum_{j=1}^sa_{ij}$. For given $u^n$, an s-stage linear Runge-Kutta prediction-correction method is given by

   1. Prediction: set the iterative initial value $k_i^{n,0}=u^n$ and let $M>0$ be a given integer, for $m=0,1,\cdots,M-1$, we calculate $k_i^{n,m+1}$ using
\begin{align*}%\label{linear-RLW-IEQ-equation1}
\left\{
\begin{aligned}
&k_i^{n,m+1}=\mathcal S\bigg (u_{i}^{n,m}+\frac{1}{6}q_{i}^{n,m}+\frac{1}{3}(u_{i}^{n,m})^2\bigg),\ u_{i}^{n,m}=u^n+\tau\sum_{j=1}^{s}a_{ij}k_j^{n,m},\\
 & q_{i}^{n,m}=q^n+\tau\sum_{j=1}^{s}a_{ij}l_j^{n,m},\ l_i^{n,m}=2u_{i}^{n,m}k_i^{n,m},\ i=1,2,\cdots,s.
\end{aligned}
\right.
\end{align*}
Then, we set $u_{i}^{n,*}=u_{i}^{n,M}$.

2.  Correction: for the predicted value $u_{i}^{n,*}$, we obtain the $(u^{n+1},q^{n+1})$ by
\begin{align}\label{linear-RLW-IEQ-equation1}
\left\{
\begin{aligned}
&k_i^n=\mathcal S\bigg (u_{i}^n+\frac{1}{6}q_{i}^n+\frac{1}{3}u_{i}^{n,*} u_{i}^n\bigg),\ u_{i}^n=u^n+\tau\sum_{j=1}^{s}a_{ij}k_j^n,\\
 & q_{i}^n=q^n+\tau\sum_{j=1}^{s}a_{ij}l_j^n,\ l_i^n=2u_{i}^{n,*}k_i,\ i=1,2,\cdots,s,
\end{aligned}
\right.
\end{align}
and $(u^{n+1},q^{n+1})$ is then updated by
\begin{align}\label{linear-RLW-IEQ-equation2}
u^{n+1}=u^n+\tau\sum_{i=1}^{s}b_ik_i^n,\ q^{n+1}=q^n+\tau\sum_{i=1}^{s}b_{i}l_i^n.
\end{align}
\end{shn}
 \begin{thm}\label{RLW-equation-th2.4} If the coefficients of \eqref{linear-RLW-IEQ-equation1} and \eqref{linear-RLW-IEQ-equation2} satisfy \eqref{RLW-RK-coeff}, {\bf Scheme 2.4} satisfies the semi-discrete energy \eqref{semi-RLW-Q-energy}.
 \end{thm}
 \begin{prf}The proof is similar to Theorem \ref{RLW-equation-th2.3}, thus, for brevity, we omit it.\qed

\end{prf}

\begin{rmk} { We should note that the proofs of Theorems \ref{RLW-equation-th2.1}-\ref{RLW-equation-th2.3} follow along similar lines of that given by Sanz-Serna in Ref. \cite{Sanzs88} and Sanz-Serna and Calvo in Ref. \cite{SCbook94}, respectively, that symplectic
Runge-Kutta methods conserve quadratic invariants.}

\end{rmk}

\begin{rmk} We should note that, the quadratic energy \eqref{RLW-Q-energy} is only equivalent to the original Hamiltonian energy \eqref{RLW-Hamiltonian-energy} in the continuous level, which implies that the proposed scheme cannot preserve the folowing original energy
\begin{align}\label{semi-discrete-RLW-Hamiltonian-energy}
 H^n =\int_{\Omega}\Big(\frac{1}{2}(u^n)^2+\frac{1}{6} (u^n)^3\Big)d{\bf x},\ n=0,1,2,\cdots,.
\end{align}

\end{rmk}

Last but not least, to improve readability, we give a summary on the construction of the proposed schemes. There are two key steps:
\begin{itemize}
\item[Step 1:] We linearize the nonlinear system \eqref{RLW-momentum-equation} and \eqref{IEQ-skew-system} by using the extrapolation technique or prediction-correction strategy for the structure matrix $\mathcal G(u)$ and $\mathcal M(\Phi)$, respectively, and a linearized conservative system is obtained, which satisfies a quadratic invariant;

    \item[Step 2:] The linearized system is then discretized by the symplectic RK method and a linear momentum-preserving scheme or energy-preserving is proposed, respectively.

\end{itemize}

\begin{rmk} The schemes obtained by using 4th-Gauss method and 6th-Gauss method (see Remark \ref{RLW-RMK-2.4}) to {\bf Scheme \ref{RLW-scheme2.4}} are denoted by 4th-LEP-PCS and 6th-LEP-PCS, respectively. In our numerical experiments, we choose $M=3$ and $M=5$ for 4th-LEP-PCS and 6th-LEP-PCS, respectively, so that 4th-LEP-PCS and 6th-LEP-PCS can achieve at fourth-order accurate and sixth-order accurate in time, respectively.

\end{rmk}
\begin{rmk} As {stated} in Remark \ref{RLW-RMK-2.5}, we also choose $t_{n-1},$ $t_{n-1}+c_i\tau,\ i=1,2,\cdots,s$ and $t_n$ as the interpolation points. %In particular, for the 6th-Gauss method, the Lagrange interpolating polynomial at $t_n+c_i\tau,\ i=1,2,3$ (i.e., $u_i^{n,*},\ i=1,2,3$) are obtained by using \eqref{RLW-LI-2.12}-\eqref{RLW-LI-2.14}.
We denote the scheme obtained by using the 6th-Gauss method and the Lagrange interpolating polynomial \eqref{RLW-LI-2.12}-\eqref{RLW-LI-2.14} to {\bf scheme 2.3} as 4th-LEPS, which may achieve fourth-order accurate in time, and similarly, we obtain the starting values $u_1^1,\ u_2^1$ and $u_3^1$ by using the 6th-Gauss method for \eqref{RLW-IEQ-reformulation1}.
\end{rmk}

\begin{rmk}
In general, a numerical scheme is called momentum/energy-preserving scheme if it can preserve the
corresponding momentum/energy conservation law exactly after temporal and spatial full-discretizations. Thus, it is very important to choose suitable
spatial discretization. Actually, we shall pay special attention to the following aspects on the spatial discretization:
\begin{itemize}
\item preserves symmetric property of the operator $\mathcal D$;
\item preserves skew-symmetric property of the operators $\mathcal G(u_{i}^{n,*})$ and $\mathcal S$ for the  momentum-preserving scheme and energy-preserving scheme, respectively;
\item the spatial accuracy is comparable to that of the time-discrete discretization.

\end{itemize}
Based on these facts and the periodic boundary condition , the standard Fourier pseudo-spectral
method is a fine choice for the spatial discretization because of its high-order accuracy and the FFT technique. For more details on the Fourier pseudo-spectral method, please refer to Refs. \cite{CQ01,ST06}.
\end{rmk}

\section{Numerical examples}\label{RLW-section-3}
In the rest of this section,  we present extensive numerical performances of the proposed schemes for solving the RLW equation in 1D and 2D.
On one hand, { we show  numerical errors and computational efficiencies of the proposed schemes and in the 1D case, we compare our schemes with the linear-implicit Crank-Nicolson momentum-ponserving  scheme (LCN-MPS) \cite{HWG-NMPDE-2020} and linear mass- and energy-ponserving scheme (LEPS) \cite{KF-2009-JJIAM}, respectively.}
 On the other hand, we apply the proposed schemes to simulate some interesting phenomena, such as
the elastic interaction of two solitary waves, the development of an undular bore and the Maxwellian pulse solitary waves. Here, we should note that LEPS also conserves a modified discrete energy and in our computation, we use the Fourier pseudo-spectral method instead of the finite difference method in space. Moreover, to quantify the numerical errors on the solution, we introduce the error functions
\begin{align*}
&e_{\infty}(t_n)=\|u(\cdot,t_n)-u^n\|_{l^{\infty}},\ e_{2}(t_n)=\|u(\cdot,t_n)-u^n\|_{l^{2}}.
\end{align*}
All simulations are performed on a Win10 machine with Intel Core i7 and 32GB
using MATLAB R2015b.

%\ e_{\mathcal M}(t_n)=|\mathcal M^n-\mathcal M^0|, \\
%&e_{\mathcal I}(t_n)=|\mathcal I^n-\mathcal I^0|,\ e_{\mathcal E}(t_n)=|\mathcal E^n-\mathcal E^0|,\  e_{\mathcal H}(t_n)=|\mathcal H^n-\mathcal H^0|,

%where $\|\cdot\|_{l^2}$ and $\|\cdot\|_{l^{\infty}}$ represent the discrete $L^2$-norm and $L^{\infty}$-norm, respectively.

 %In addition, in our computations, we choose $M= 3$ in the prediction step for 4th-LEPS and 4th-LEP-PCS, while for 6th-LEPS and 6th-LEP-PCS, we set $M=5$.

%For the more details on the construction of $u_N({\bf x},t)$, please refer to \cite{GZW2020jcp,LS2020jsc}.

\subsection{RLW equation in 1D}
We first consider the RLW equation in 1D given by \cite{BBM-72-PTRSL,Olver-79-MPPS,Peregrine66jfm}
\begin{align}\label{RLW-1d-equation}
\left\{
\begin{aligned}
&\partial_tu(x,t)+\alpha\partial_xu(x,t)+\alpha u(x,t)\partial_x u(x,t)-\mu\partial_{xxt}u(x,t)=0,\ x\in\Omega,\\
&u(x,0)=3c\text{sech}^2(k(x-x_0)),\ x\in{\Omega},
\end{aligned}
\right.
\end{align}
which possesses the following solitary wave solution
\begin{align}
u(x,t)=3c\text{sech}^2(k(x-vt-x_0)),\ x\in\Omega,\ t\ge0,
\end{align}
where $k=\frac{1}{2}\sqrt{\frac{ c}{\mu(1+ c)}}$ and $v=\alpha(1+ c)$.

%Here, the initial condition is chose from the exact solution.

\begin{table}[H]
\tabcolsep=9pt
\footnotesize
\renewcommand\arraystretch{1.1}
\centering
\caption{{Numerical errors and convergence rates for the different schemes with various time steps and the Fourier node $2048$ at $T=1$.}}\label{Tab:RLW-equation:1}
\begin{tabular*}{\textwidth}[h]{@{\extracolsep{\fill}} c c c c c c}\hline
{Scheme\ \ } &{$\tau$} &{$e_{2}(t_n=1)$} &{order}& {$e_{\infty}(t_n=1)$}&{order}  \\     %% 第1 行
\hline
\multirow{4}{*}{{LCN-MPS\cite{HWG-NMPDE-2020}}}  &{$\frac{1}{100}$}& {1.684e-03}&{-} &{1.117e-03} & {-}\\[1ex]%&{2.8}\\
  {}  &{$\frac{1}{200}$}& {4.251e-04}&{1.986} &{2.812e-04} & {1.990}\\[1ex] %&{4.3}  \\   %% 第2行
   {}  &{$\frac{1}{400}$}& {1.068e-04}&{1.992} &{7.052e-05} &{1.995} \\[1ex]%\hline
    {}  &{$\frac{1}{800}$}& {2.677e-05}&{1.996} &{1.766e-05} &{1.998} \\
 \multirow{4}{*}{{LEPS\cite{KF-2009-JJIAM}}}  &{$\frac{1}{100}$}& {1.164e-03}&{-} &{5.086e-04} & {-}\\[1ex]%&{2.8}\\
  {}  &{$\frac{1}{200}$}& {2.913e-04}&{1.998} &{1.275e-04} & {1.996}\\[1ex] %&{4.3}  \\   %% 第2行
   {}  &{$\frac{1}{400}$}& {7.287e-05}&{1.999} &{3.190e-05} &{1.998} \\[1ex]%\hline
    {}  &{$\frac{1}{800}$}& {1.822e-05}&{2.000} &{7.980e-06} &{1.999} \\\hline
% \multirow{4}{*}{3rd-LEPS}  &{$\frac{1}{100}$}& {4.072e-05}&{-} &{2.763e-05} & {-}\\[1ex]%&{2.8}\\
%  {}  &{$\frac{1}{200}$}& {5.103e-06}&{2.996} &{3.457e-06} & {2.998}\\[1ex] %&{4.3}  \\   %% 第2行
%   {}  &{$\frac{1}{400}$}& {6.388e-07}&{2.998} &{4.323e-07} &{2.999} \\[1ex]%\hline
%    {}  &{$\frac{1}{800}$}& {7.991e-08}&{2.999} &{5.406e-07} &{3.000} \\\hline
    \multirow{4}{*}{4th-LMPS}  &{$\frac{1}{100}$}& {5.646e-07}&{-} &{3.781e-07} & {-}\\[1ex]%&{2.8}\\
  {}  &{$\frac{1}{200}$}& {3.540e-08}&{3.995} &{2.365e-08} & {3.998}\\[1ex] %&{4.3}  \\   %% 第2行
   {}  &{$\frac{1}{400}$}& {2.216e-09}&{3.998} &{1.479e-09} &{3.999} \\[1ex]%\hline
    {}  &{$\frac{1}{800}$}& {1.386e-10}&{3.999} &{9.244e-11} &{4.000} \\
     \multirow{4}{*}{4th-LEPS}  &{$\frac{1}{100}$}& {5.942e-07}&{-} &{3.794e-07} & {-}\\[1ex]%&{2.8}\\
  {}  &{$\frac{1}{200}$}& {3.724e-08}&{3.996} &{2.372e-08} & {3.999}\\[1ex] %&{4.3}  \\   %% 第2行
   {}  &{$\frac{1}{400}$}& {2.331e-09}&{3.998} &{1.483e-09} &{3.999} \\[1ex]%\hline
    {}  &{$\frac{1}{800}$}& {1.458e-10}&{3.999} &{9.272e-11} &{4.000} \\%\hline
    \multirow{4}{*}{4th-LMP-PCS}  &{$\frac{1}{100}$}& {1.288e-07}&{-} &{7.957e-08} & {-}\\[1ex]%&{2.8}\\
  {}  &{$\frac{1}{200}$}& {8.060e-09}&{3.999} &{4.967e-09} & {4.002}\\[1ex] %&{4.3}  \\   %% 第2行
   {}  &{$\frac{1}{400}$}& {5.040e-10}&{3.999} &{3.102e-10} &{4.001} \\[1ex]%\hline
    {}  &{$\frac{1}{800}$}& {3.151e-11}&{4.000} &{1.938e-11} &{4.001} \\
    \multirow{4}{*}{4th-LEP-PCS}  &{$\frac{1}{100}$}& {1.246e-07}&{-} &{8.185e-08} & {-}\\[1ex]%&{2.8}\\
  {}  &{$\frac{1}{200}$}& {7.775e-09}&{4.002} &{5.115e-09} & {4.000}\\[1ex] %&{4.3}  \\   %% 第2行
   {}  &{$\frac{1}{400}$}& {4.856e-10}&{4.001} &{3.196e-10} &{4.000} \\[1ex]%\hline
    {}  &{$\frac{1}{800}$}& {3.033e-11}&{4.001} &{1.999e-11} &{3.999} \\\hline%\hline
    \multirow{4}{*}{6th-LMP-PCS}  &{$\frac{1}{10}$}& {2.092e-06}&{-} &{1.318e-06} & {-}\\[1ex]%&{2.8}\\
  {}  &{$\frac{1}{20}$}& {3.412e-08}&{5.938} &{2.044e-08} & {6.011}\\[1ex] %&{4.3}  \\   %% 第2行
   {}  &{$\frac{1}{40}$}& {5.454e-10}&{5.967} &{3.185e-10} &{6.004} \\[1ex]%\hline
    {}  &{$\frac{1}{80}$}& {8.622e-12}&{5.983} &{4.952e-12} &{6.006} \\
    \multirow{4}{*}{6th-LEP-PCS}  &{$\frac{1}{10}$}& {2.282e-06}&{-} &{1.484e-06} & {-}\\[1ex]%&{2.8}\\
  {}  &{$\frac{1}{20}$}& {3.564e-08}&{6.000} &{2.357e-08} & {5.976}\\[1ex] %&{4.3}  \\   %% 第2行
   {}  &{$\frac{1}{40}$}& {5.561e-10}&{6.002} &{3.704e-10} &{5.992} \\[1ex]%\hline
    {}  &{$\frac{1}{80}$}& {8.682e-12}&{6.001} &{5.798e-12} &{5.998} \\\hline%\hline

\end{tabular*}
\end{table}
First of all, we verify  the proposed schemes can achieve high-order accuracy in time. We set the computational domain $\Omega=[-100,100]$ and take parameters $\alpha=\mu=1,\ x_0=0$ and $c=3$. Then, we fix the Fourier node $2048$ so that spatial errors play no role here. We calculate the numerical errors between the numerical solution and { the exact solution} and convergence rates with various time steps at $T=1$, which are summarized in Tab. \ref{Tab:RLW-equation:1}. It is clear to see that (i) LCN-MPS and LEPS are second-order accurate in time and their numerical errors are much larger than others; (ii) 4th-LMPS, 4th-LEPS, 4th-LMP-PCS and 4th-LEP-PCS are all fourth-order accurate in time; (iii) 6th-LMP-PCS and 6th-LEP-PCS are sixth-order accurate in time {indeed}; (iv) the error provided by 6th-LMP-PCS is the smallest and LCN-MPS admits the largest numerical error.

\begin{figure}[H]
\centering\begin{minipage}[t]{70mm}
\includegraphics[width=65mm]{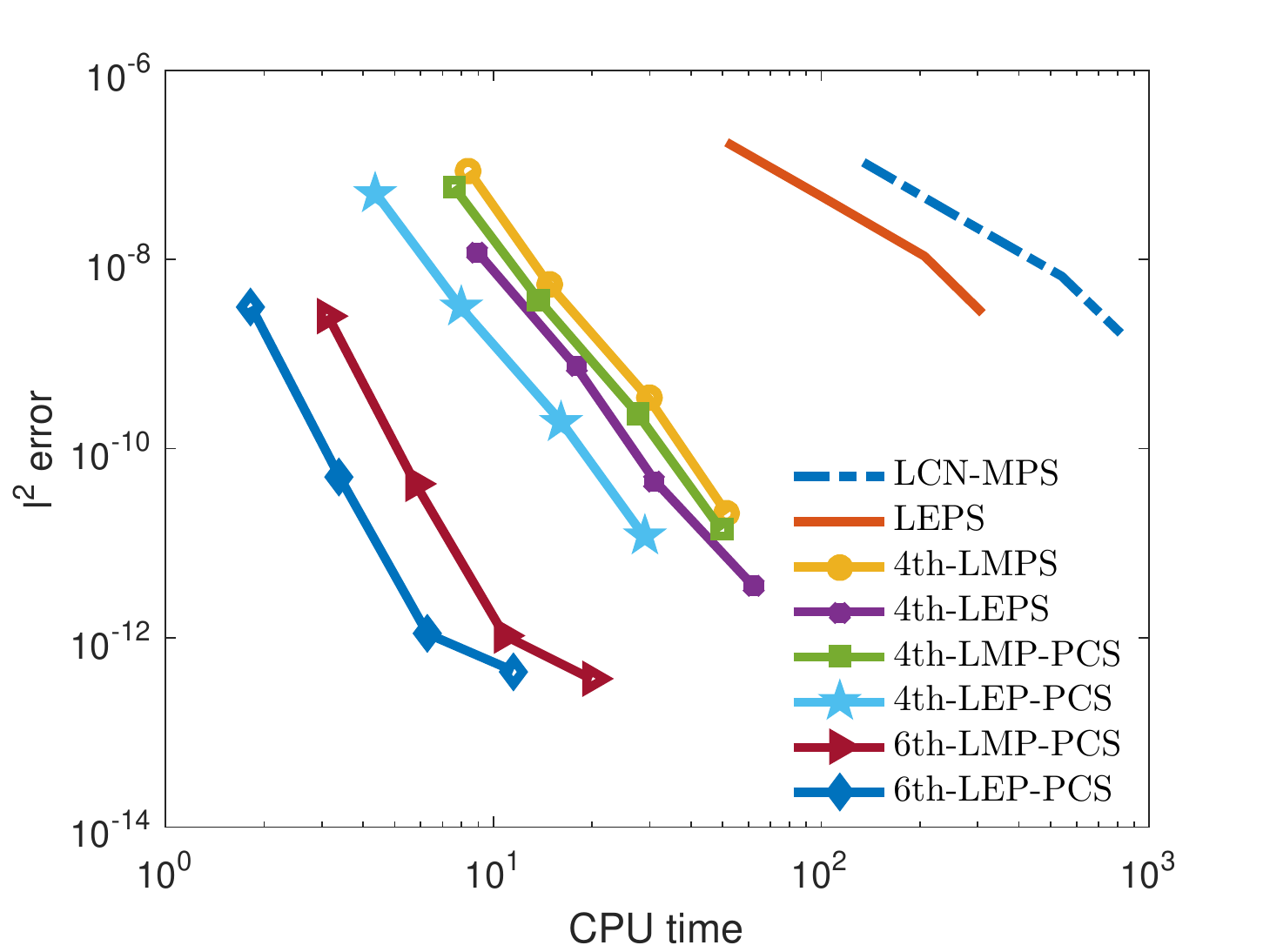}
\end{minipage}
\begin{minipage}[t]{70mm}
\includegraphics[width=65mm]{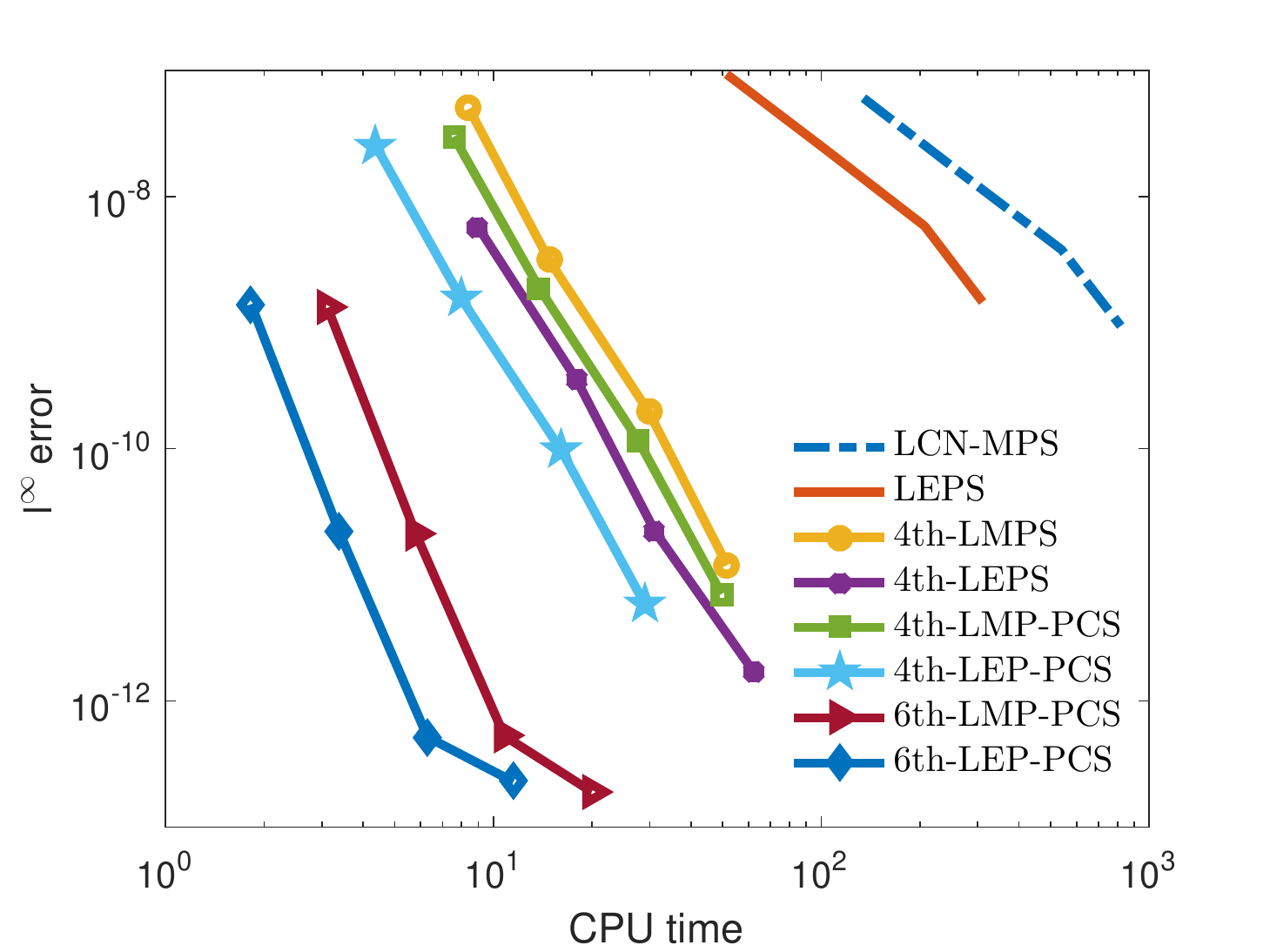}
\end{minipage}
%\centering\begin{minipage}[t]{60mm}
%\includegraphics[width=60mm]{1d_NLS_fig4.eps}
%\end{minipage}
%\begin{minipage}[t]{60mm}
%\includegraphics[width=60mm]{1d_NLS_fig5.eps}
%\end{minipage}
 \caption{{ The numerical error versus the CPU time using the different numerical schemes with the various time steps and the Fourier node 3072 at $T=10$ for the RLW equation in 1D \eqref{RLW-1d-equation}.}}\label{RLW1d-CPU}
\end{figure}

Subsequently, we display the global $l^2$-
 and $l^{\infty}$- errors (i.e., $e_2(t_n)$ and $e_{\infty}(t_n)$) of $u$ versus the CPU time using the eight
numerical schemes with the Fourier node 3072 and various time steps at $T=10$, respectively, in Fig. \ref{RLW1d-CPU}. For a given global error, we observe that {(i) the cost of LCN-MPS is the most expensive;} (ii) the cost of 6th-LEP-PCS is the cheapest; (iii) the costs of 4th-LEPS and 4th-LEP-PCS are much cheaper than the ones provided by 4th-LMPS and 4th-LMP-PCS, respectively.

\begin{figure}[H]
\centering\begin{minipage}[t]{65mm}
\includegraphics[width=65mm]{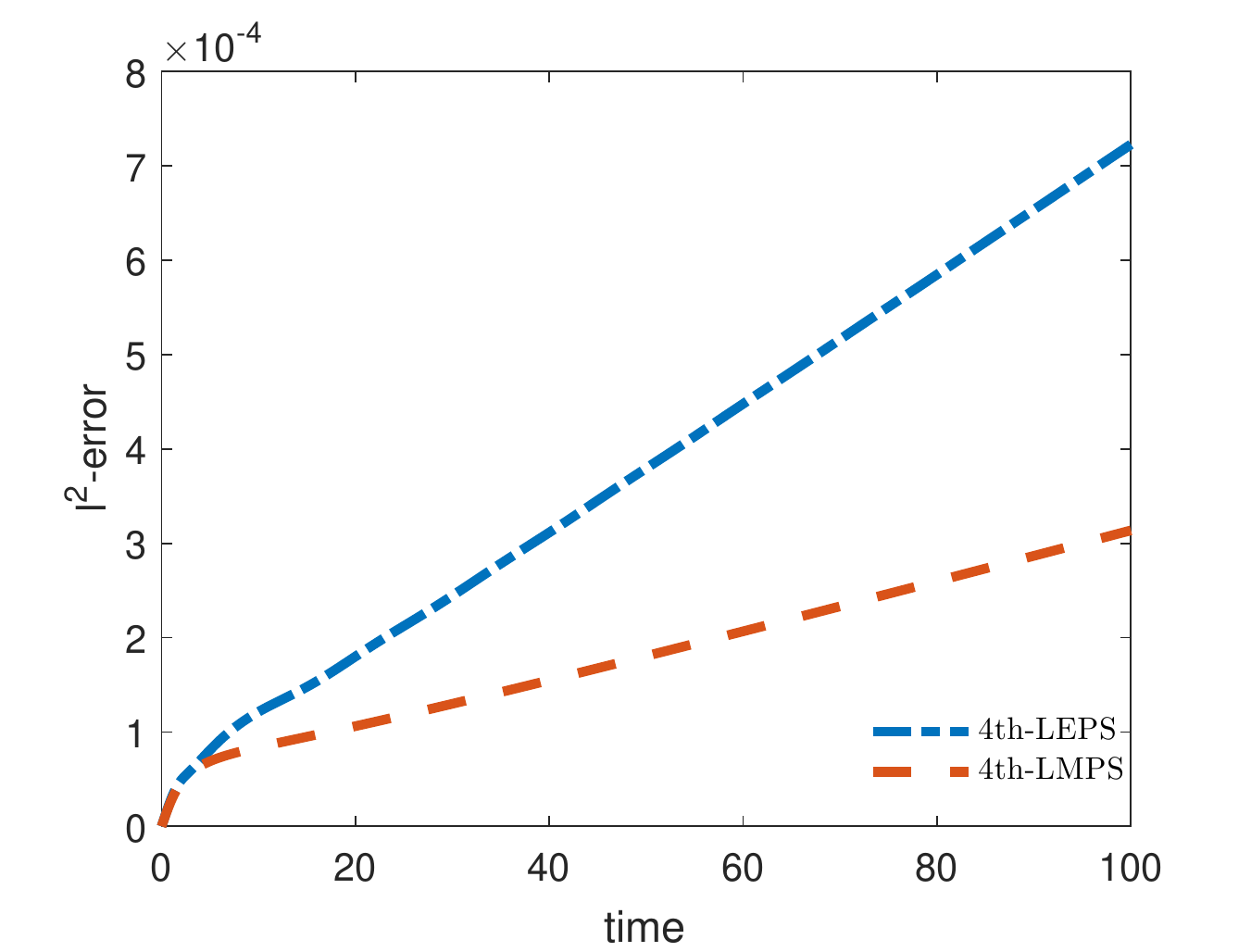}
\end{minipage}
\begin{minipage}[t]{65mm}
\includegraphics[width=65mm]{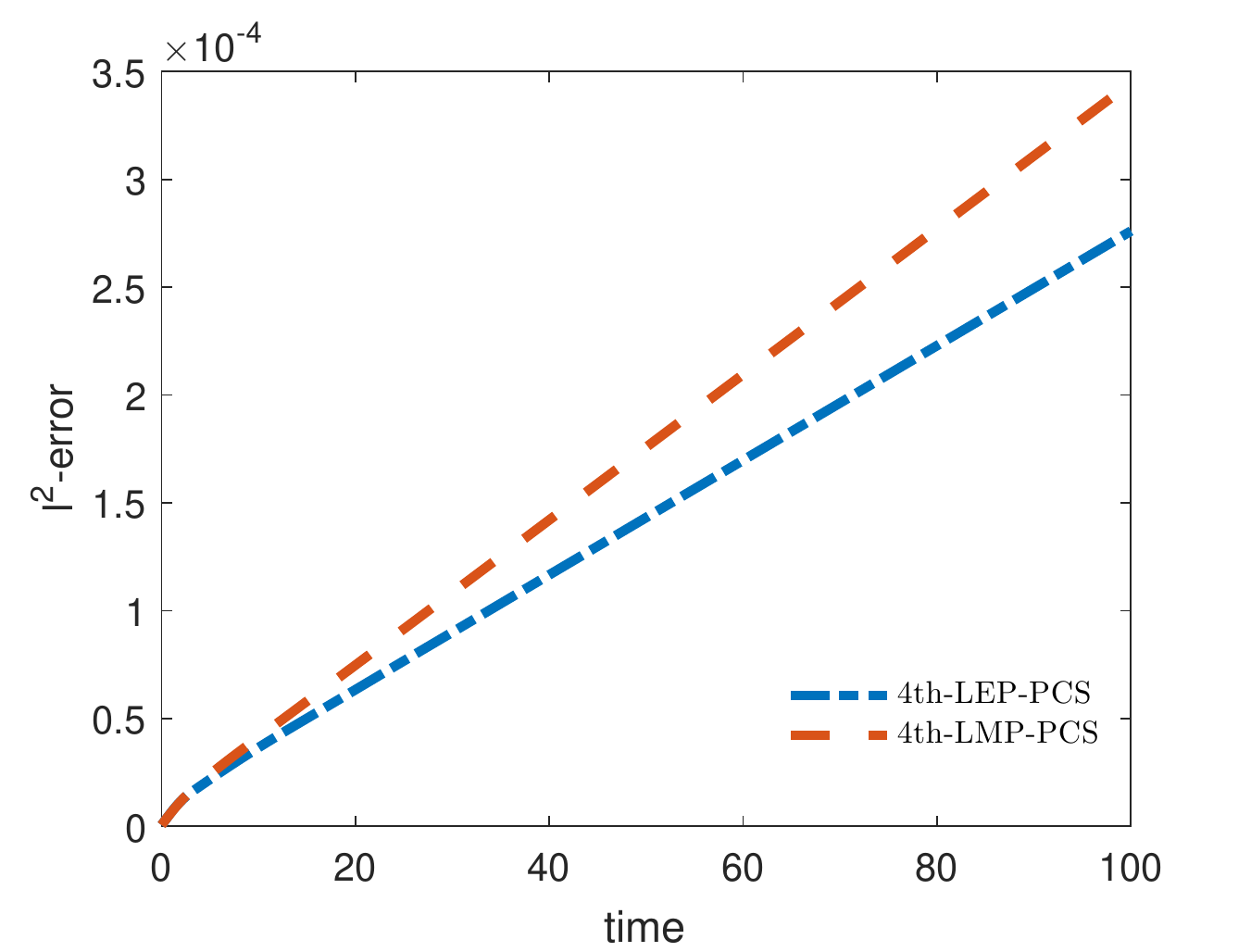}
\end{minipage}
\begin{minipage}[t]{65mm}
\includegraphics[width=65mm]{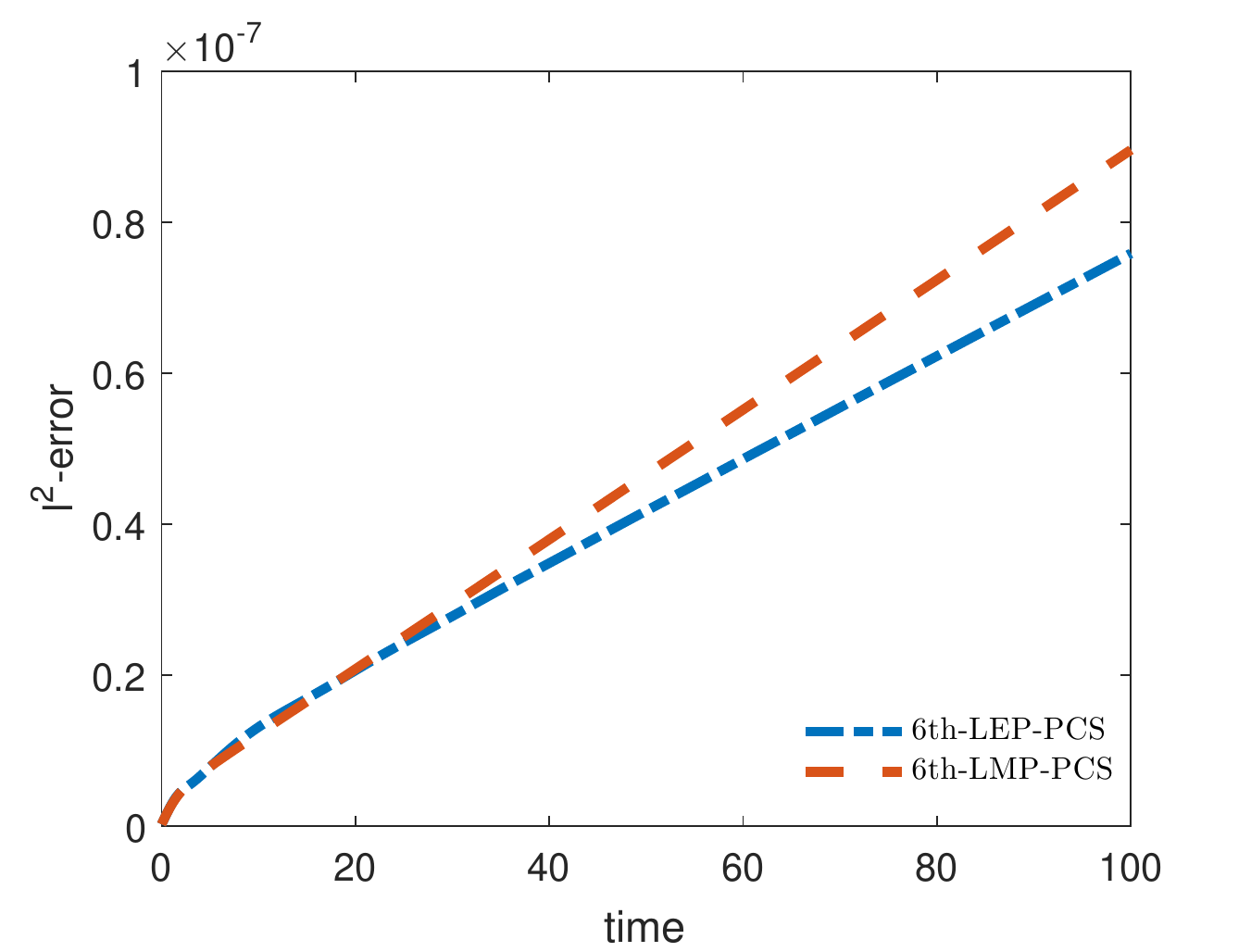}
\end{minipage}
%\centering\begin{minipage}[t]{60mm}
%\includegraphics[width=60mm]{1d_NLS_fig4.eps}
%\end{minipage}
%\begin{minipage}[t]{60mm}
%\includegraphics[width=60mm]{1d_NLS_fig5.eps}
%\end{minipage}
 \caption{{The numerical error using the different numerical schemes with the time step $\tau=0.1$ and the Fourier node 2048 at $T=100$ for the RLW equation in 1D \eqref{RLW-1d-equation}.}}\label{RLW1d-error}
\end{figure}

{Furthermore, we also examine how errors accumulate over a long time by taking parameters $\alpha=\mu=1,\ x_0=0$ and $c=1$ with the time step $\tau=0.1$ and the Fourier node 2048, respectively. The results are summarized in Fig. \ref{RLW1d-error}, which implies that the solution error of the proposed schemes grows linearly with time. Here, we set the computational domain $\Omega=[-250,250]$ to prevent the significant error introduced by the boundary.}

\begin{figure}[H]
\centering\begin{minipage}[t]{70mm}
\includegraphics[width=65mm]{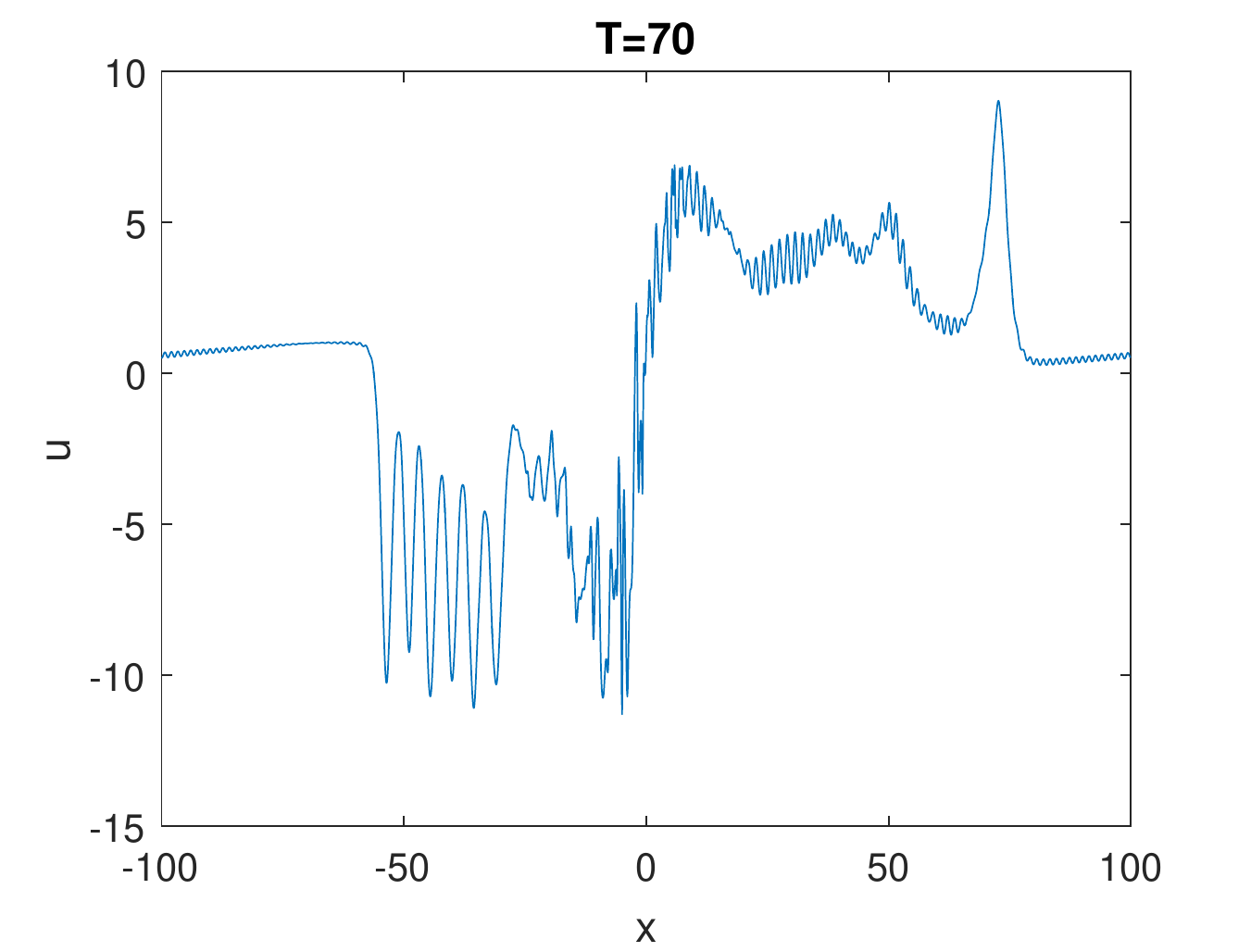}
\caption*{(a) 4th-LEPS}
\end{minipage}
\begin{minipage}[t]{70mm}
\includegraphics[width=65mm]{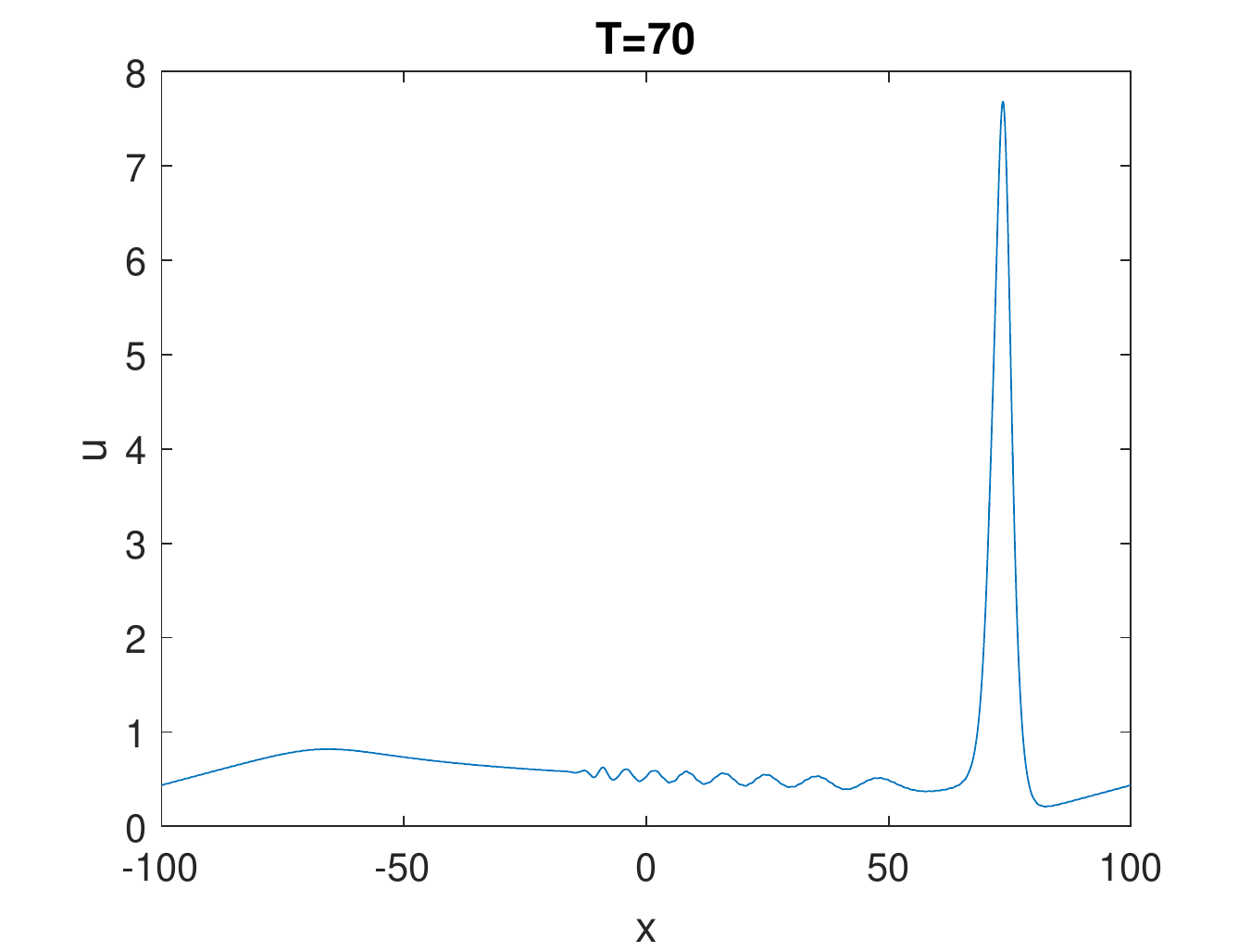}
\caption*{(b) 4th-LMPS}
\end{minipage}
\centering\begin{minipage}[t]{70mm}
\includegraphics[width=65mm]{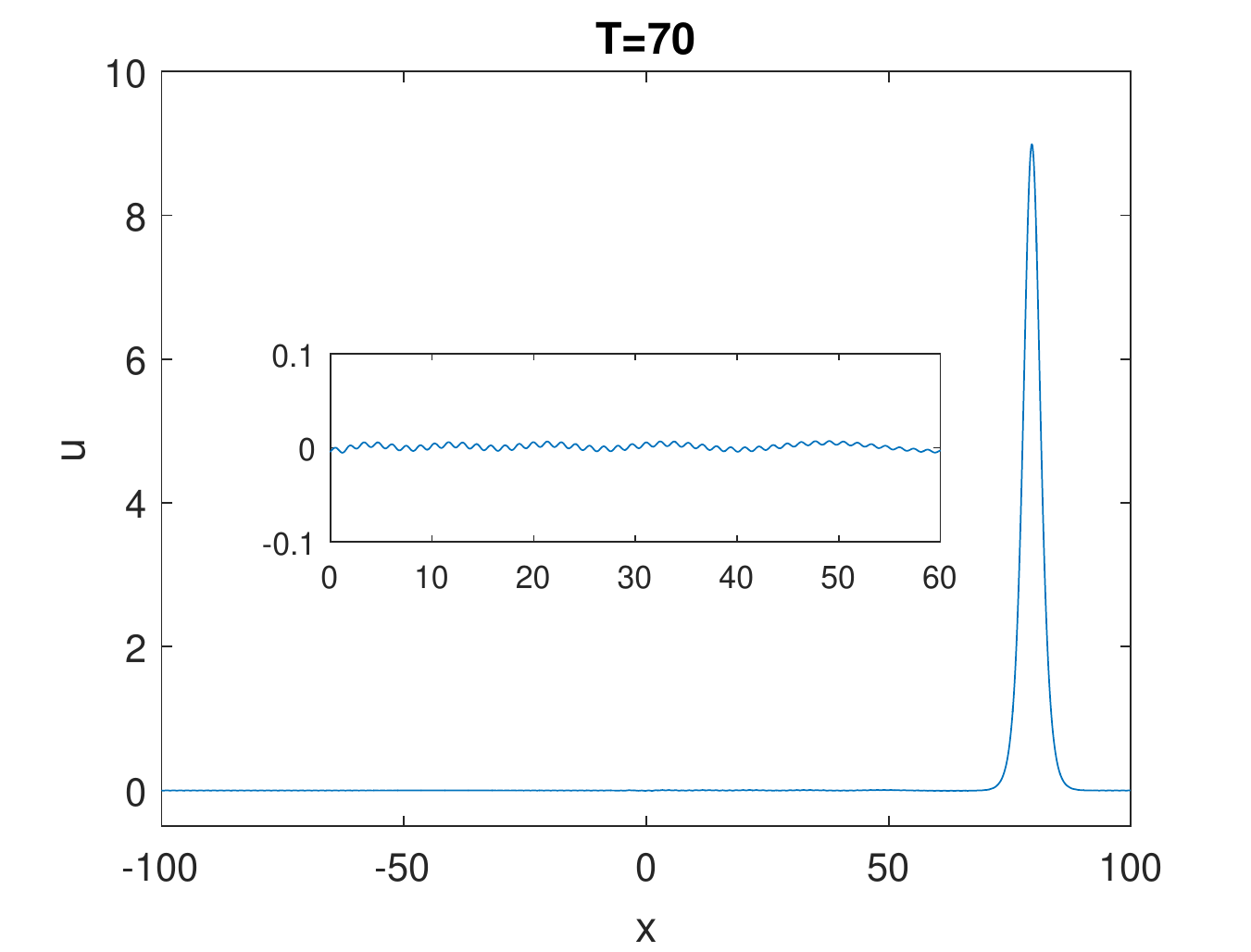}
\caption*{(c) 4th-LEP-PCS}
\end{minipage}
\begin{minipage}[t]{70mm}
\includegraphics[width=65mm]{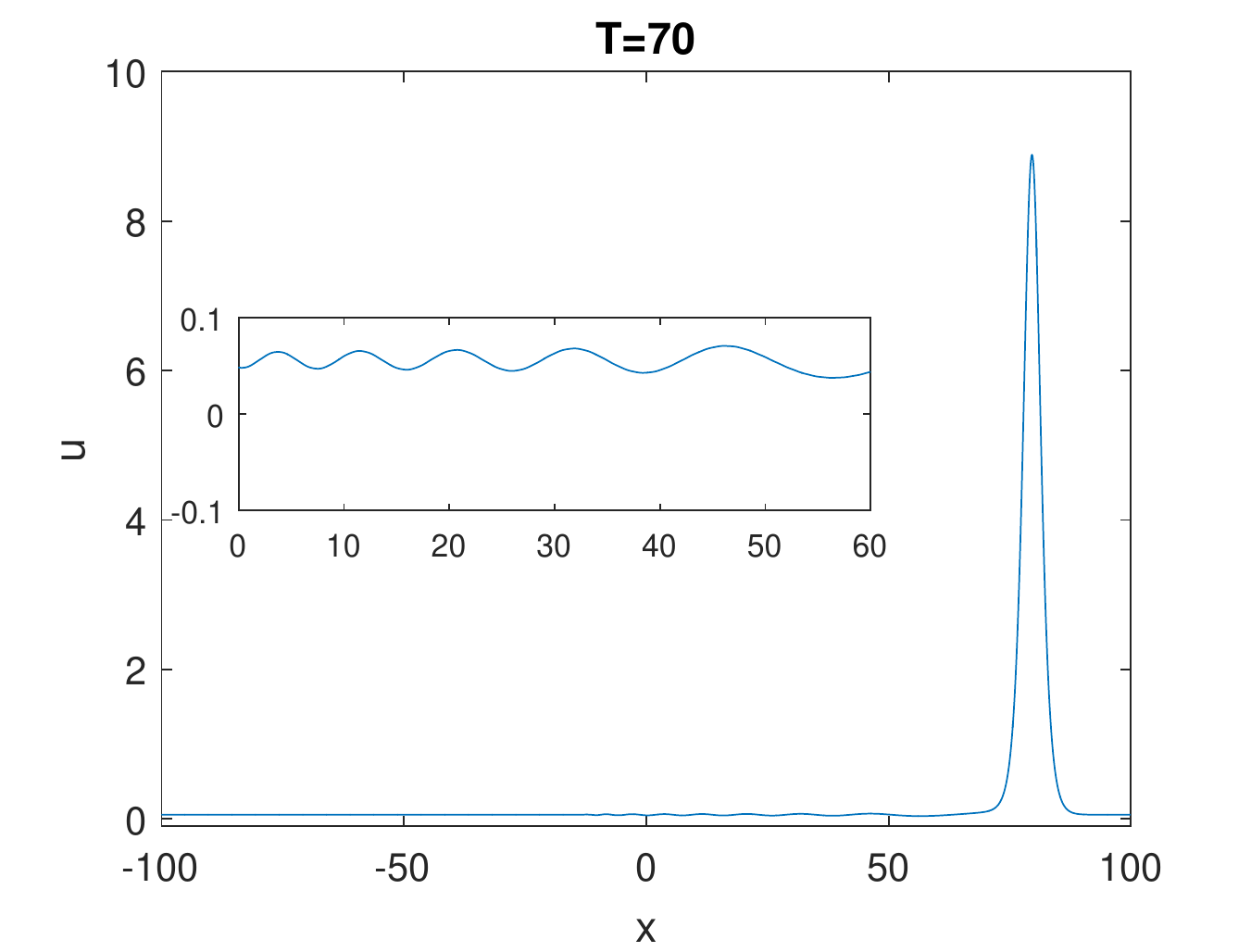}
\caption*{(d) 4th-LMP-PCS}
\end{minipage}
\centering\begin{minipage}[t]{70mm}
\includegraphics[width=65mm]{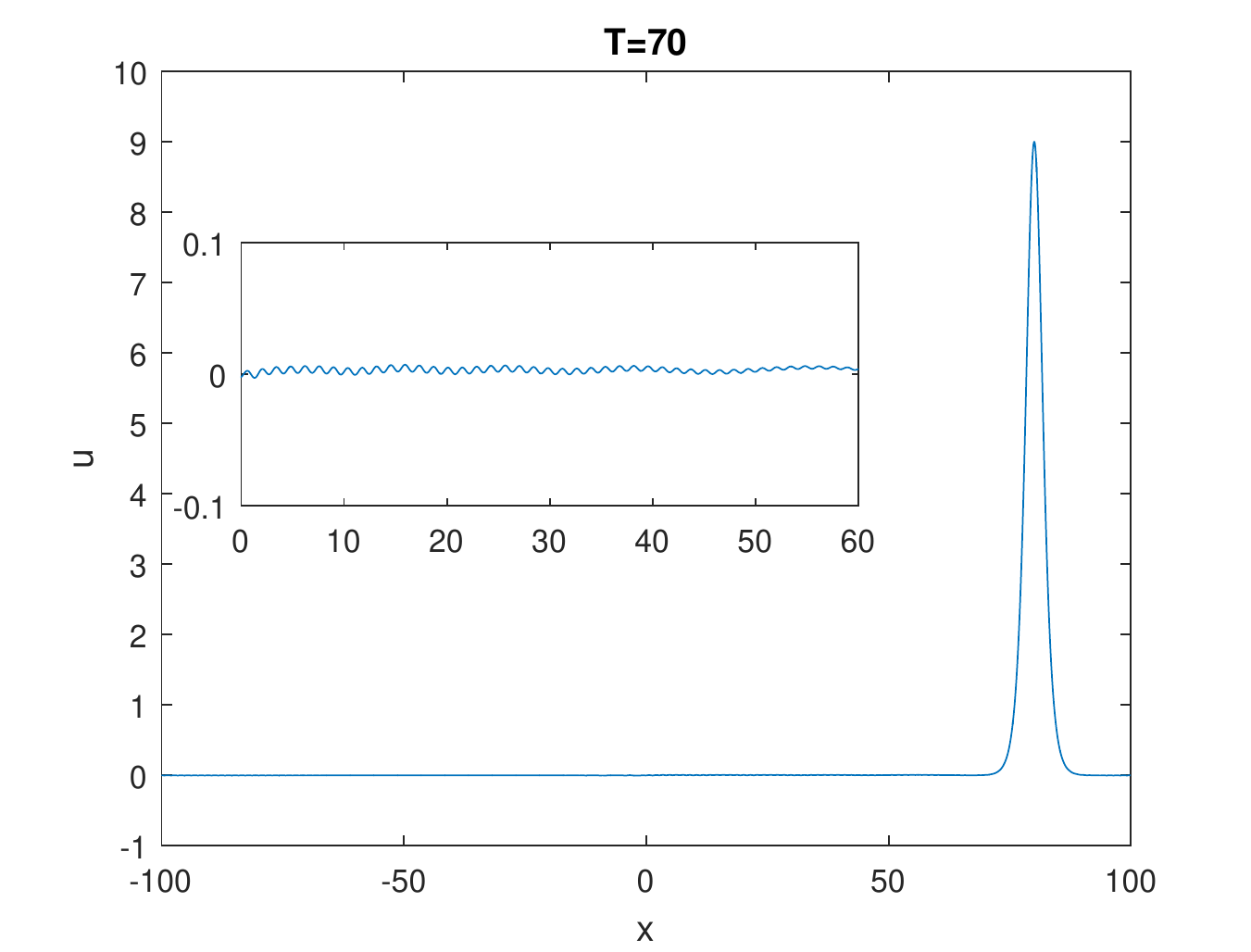}
\caption*{(e) 6th-LEP-PCS}
\end{minipage}
\begin{minipage}[t]{70mm}
\includegraphics[width=65mm]{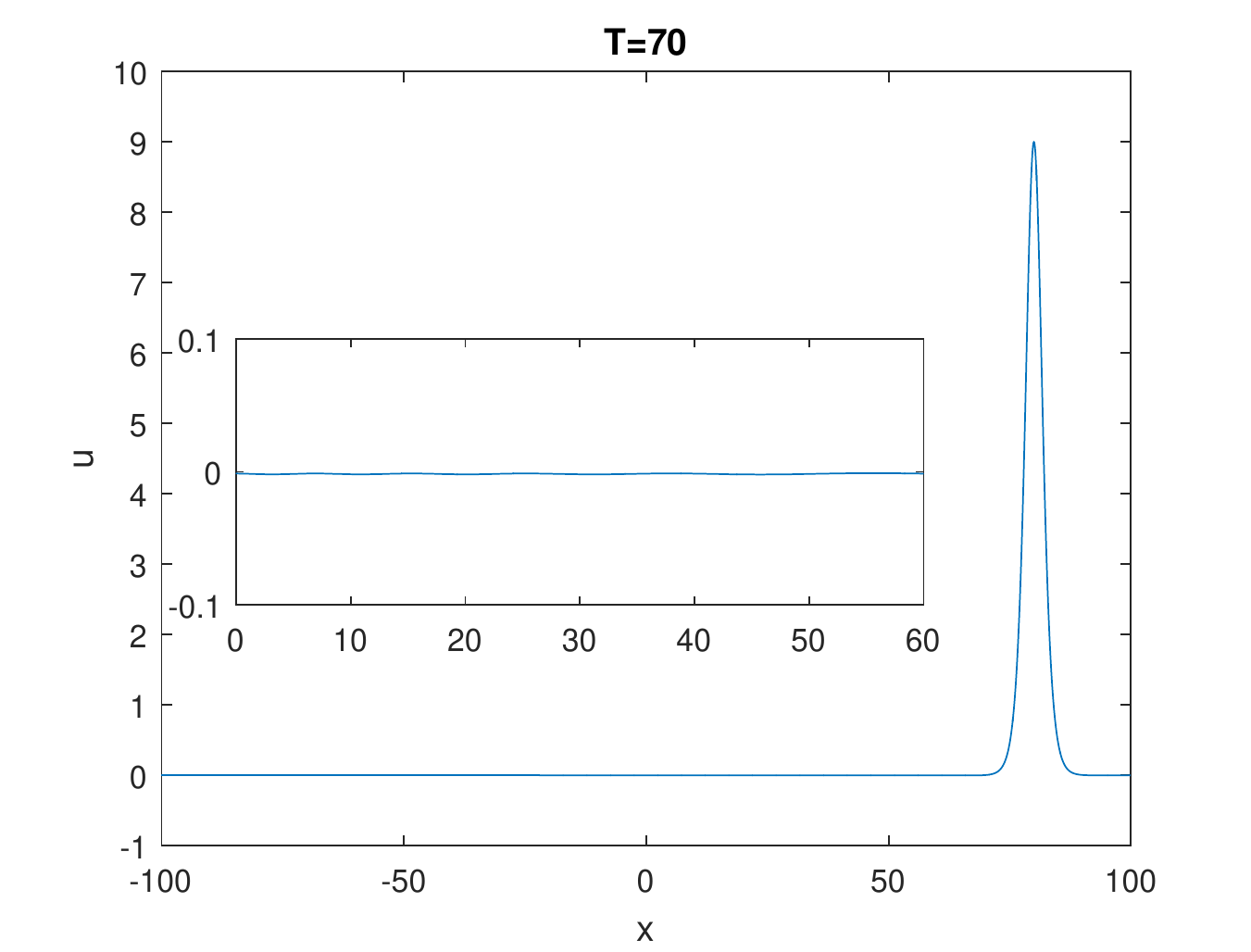}
\caption*{(f) 6th-LMP-PCS}
\end{minipage}
\centering\begin{minipage}[t]{70mm}
\includegraphics[width=65mm]{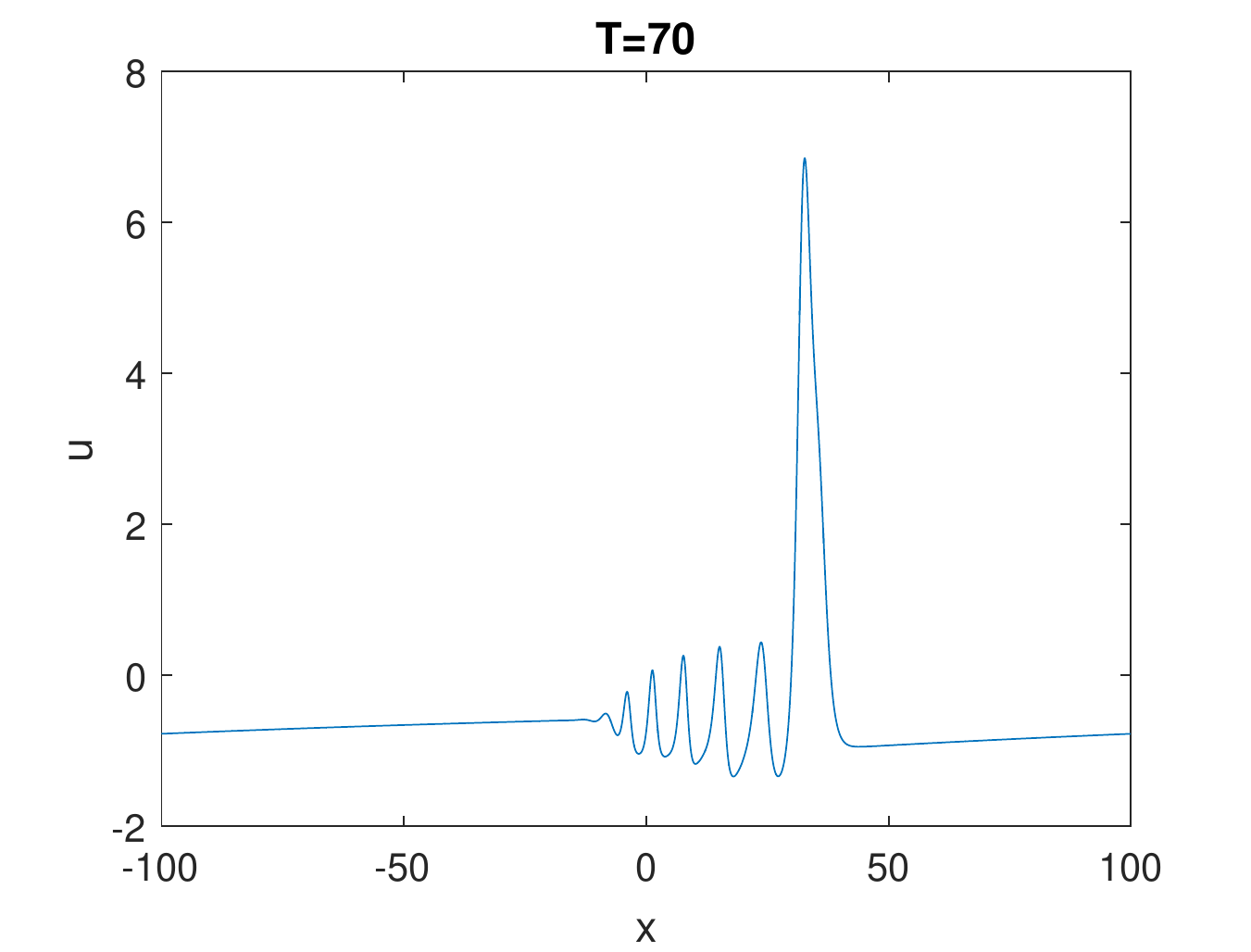}
\caption*{{(g) LCN-MPS}}
\end{minipage}
\centering\begin{minipage}[t]{70mm}
\includegraphics[width=65mm]{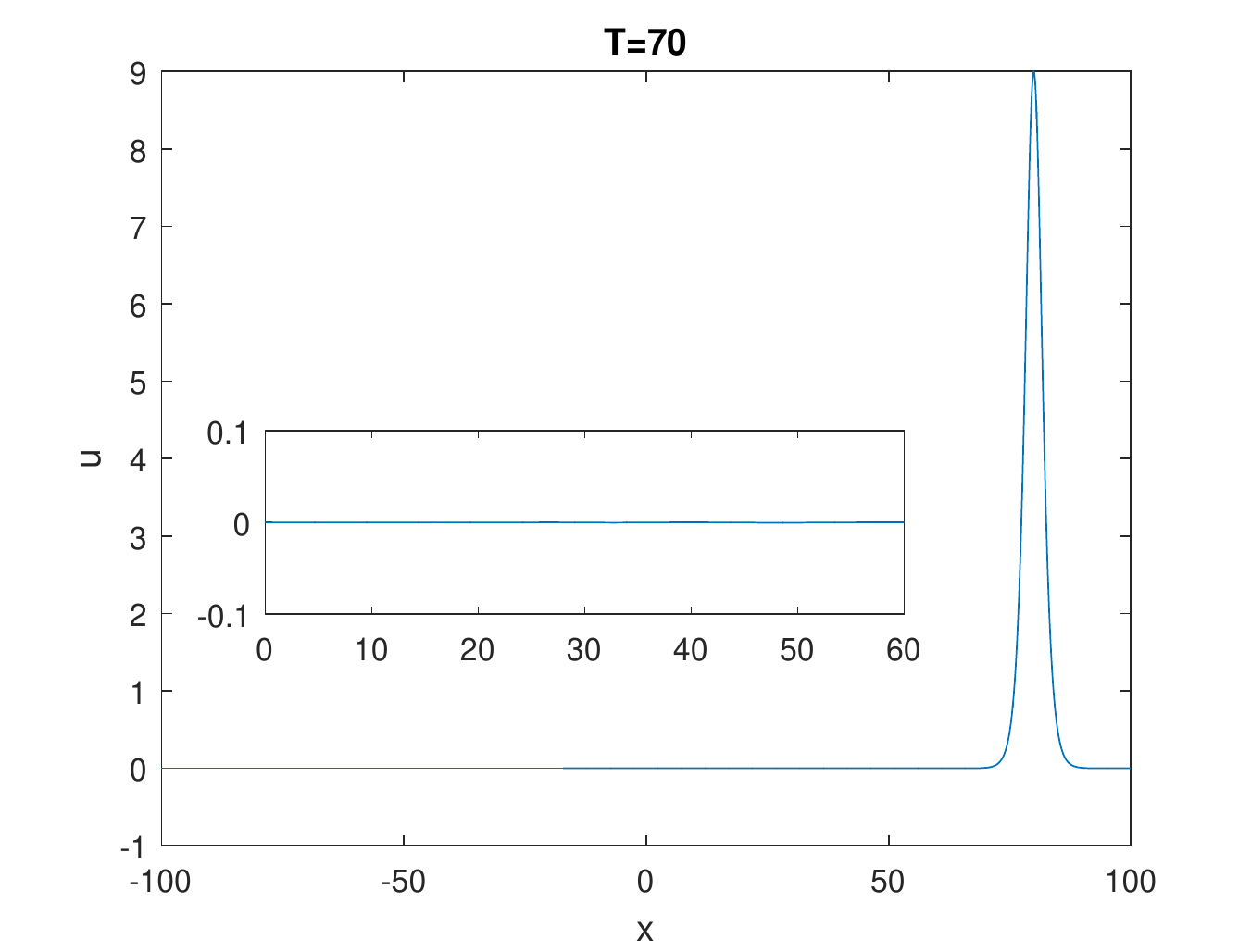}
\caption*{{(h) The reference solution obtained by LEPS with the time step $\tau=0.001$}}
\end{minipage}
\caption{ The profile of $u$ provided by the proposed schemes with the time step $\tau=0.35$ and the Fourier node 2048 at $T=70$ for the RLW equation in 1D \eqref{RLW-1d-equation}.}\label{RLW1d-soliton}
\end{figure}

Next, we consider to test the { robustness} of the proposed schemes in simulating evolution of the soliton as a large time step $\tau=0.35$ is chosen. The results are summarized in Fig. \ref{RLW1d-soliton}. In particular, we also display the reference solution obtained by LEPS with the small time step $\tau=0.001$ in Fig. \ref{RLW1d-soliton} (h).
We find that (i) LCN-MPS, 4th-LEPS, 4th-LMPS and 4th-LEP-PCS { fail to} capture the waveform of the soliton;  (ii) although  4th-LMP-PCS and 6th-LEP-PCS can capture the amplitude and waveform, but the small disturbance emergences; (iii)  6th-LMP-PCS can simulate the soliton well. Based on these observations, we find the following facts: (i) the sixth-order schemes are more robust than the fourth-order ones at the large time step; (ii) the scheme proposed by the prediction-correction strategy shows more robust than the one provided by the extrapolation technique; (iii)
the linear momentum-preserving scheme is more robust than the linear energy-preserving at large time step. It is remarked that the result provided by LEPS is omitted since the computing result is below-up.

%\begin{rmk} We should note that when we set $M=4$, 4th-LMP-PCS can accurately simulate the soltion, however, the result provided by 4th-LEP-PCS  {\color{blue} still  shows a small disturbance}.

%\end{rmk}

Finally, we consider the interaction of two solitary waves of the RLW equation in 1D \eqref{RLW-1d-equation} by choosing the initial condition
\begin{align}
u(x,0)=3c_1\text{sech}^2(k_1(x-x_1))+3c_2\text{sech}^2(k_2(x-x_2)), \ x\in\Omega,
\end{align}
where $k_1=\frac{1}{2}\sqrt{\frac{c_1}{\mu(1+c_1)}}$ and $k_2=\frac{1}{2}\sqrt{\frac{c_2}{\mu(1+c_2)}}$.

\begin{figure}[H]
\centering
\begin{minipage}[t]{75mm}
\includegraphics[width=70mm]{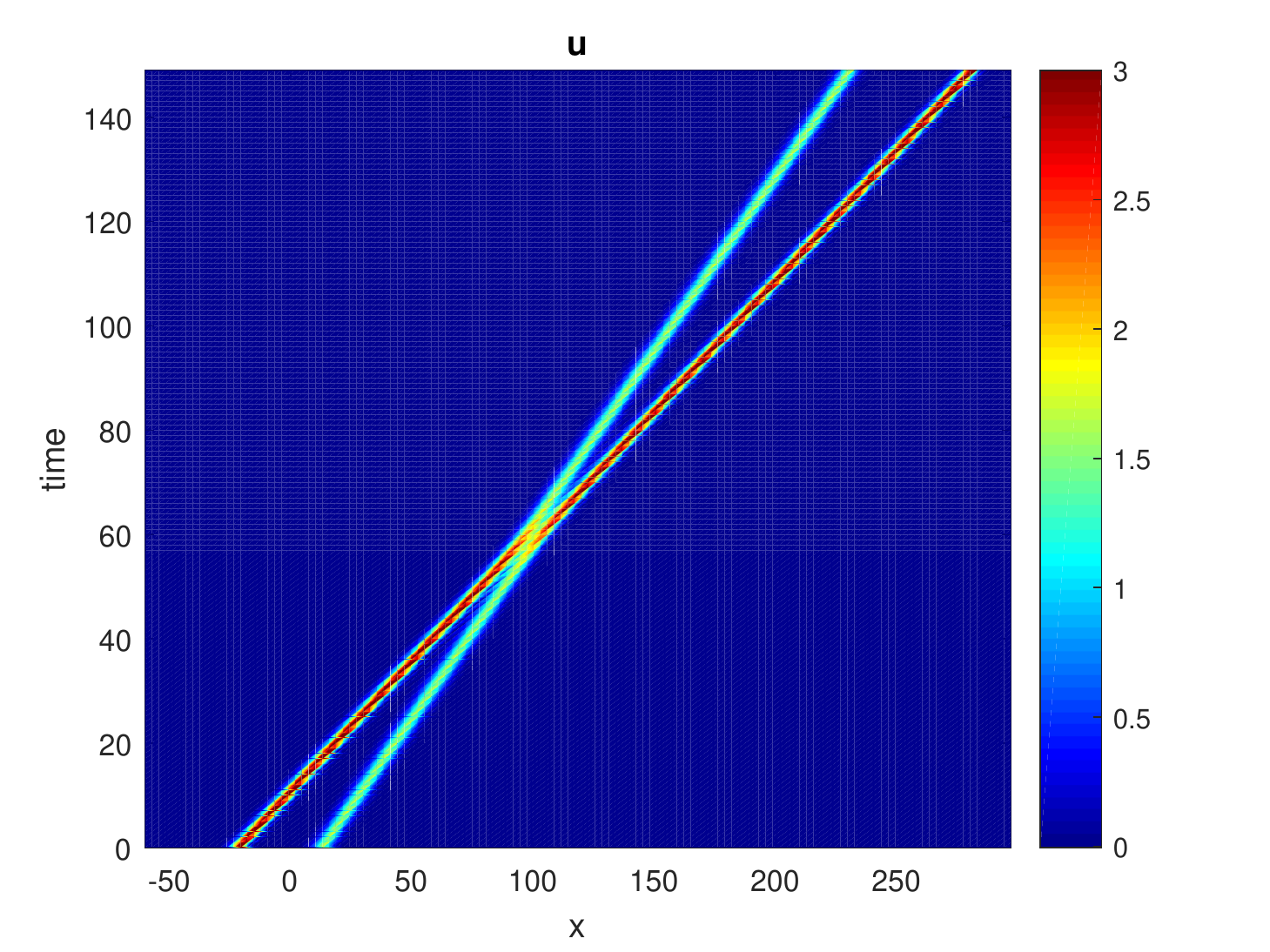}
\end{minipage}
\caption{ The profile of $u$ provided by 6th-LEP-PCS over time interval $[0,150]$ with the time step $\tau=0.01$ and the Fourier node 1024, respectively, for the RLW equation in 1D \eqref{RLW-1d-equation}.}\label{RLW1d-fig-1}
\end{figure}

\begin{figure}[H]
\centering\begin{minipage}[t]{60mm}
\includegraphics[width=60mm]{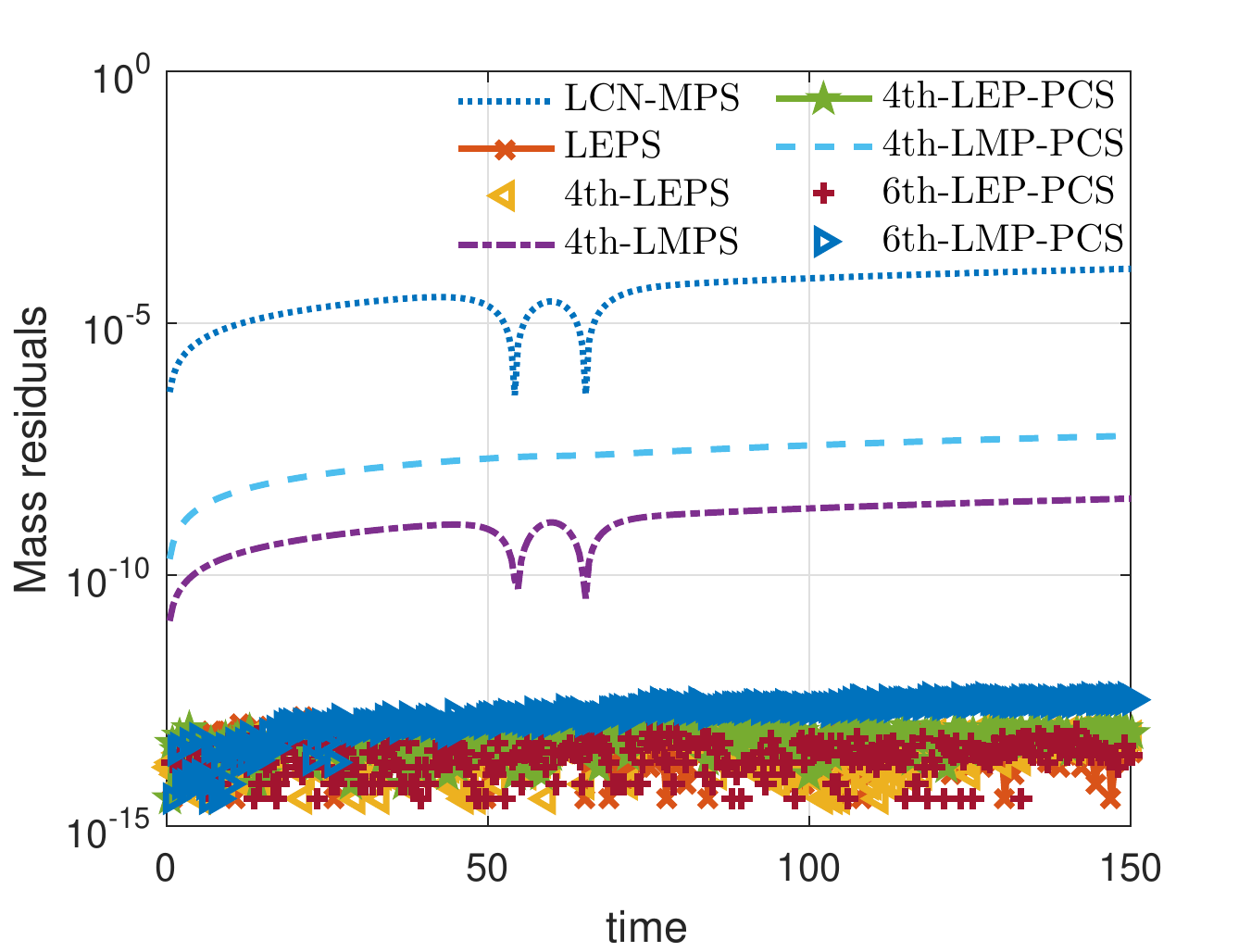}
%\caption*{(a) Mass}
\end{minipage}
\begin{minipage}[t]{60mm}
\includegraphics[width=60mm]{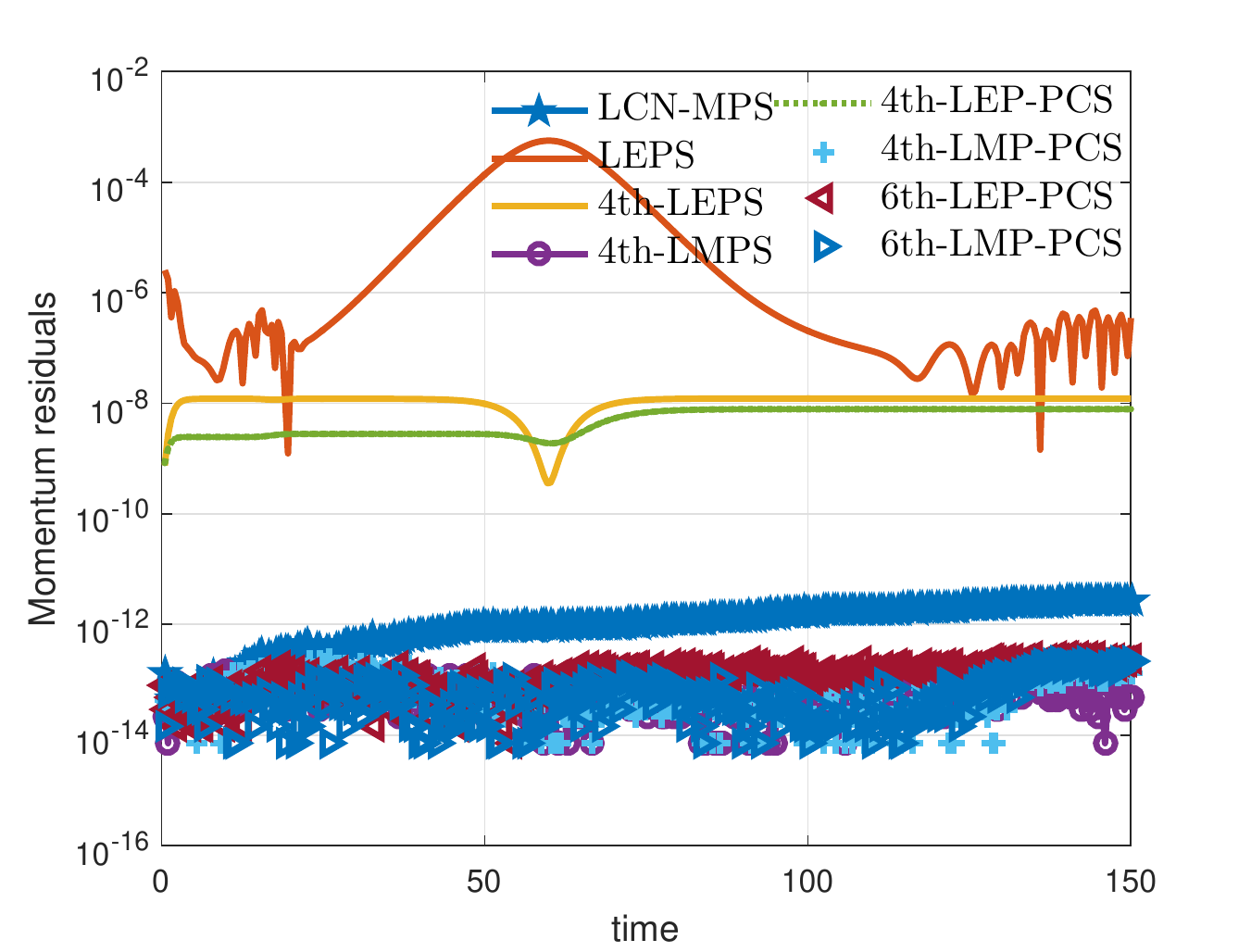}
%\caption*{(b) Momentum}
\end{minipage}
\centering\begin{minipage}[t]{60mm}
\includegraphics[width=60mm]{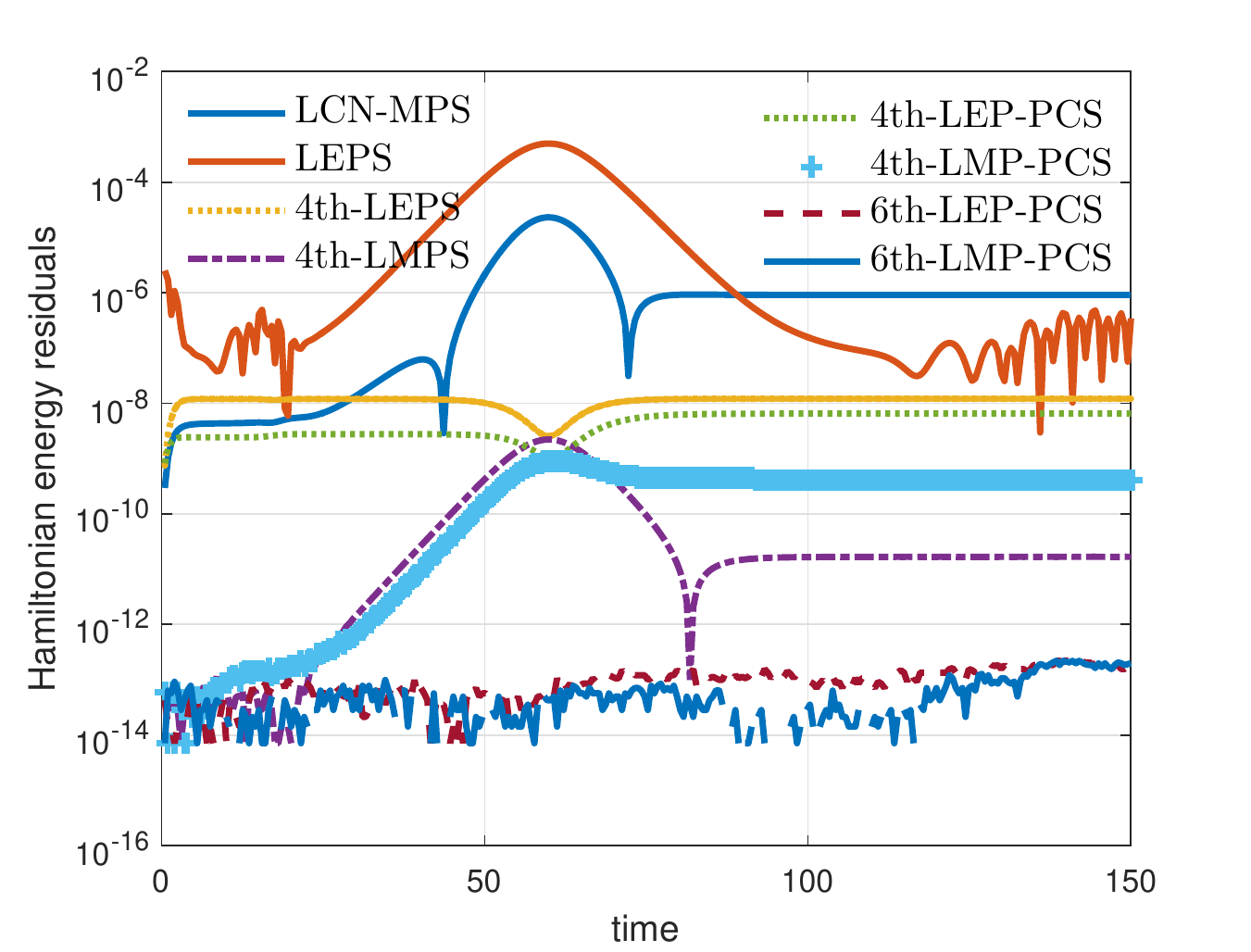}
%\caption*{(c) Hamiltonian energy}
\end{minipage}
\begin{minipage}[t]{60mm}
\includegraphics[width=60mm]{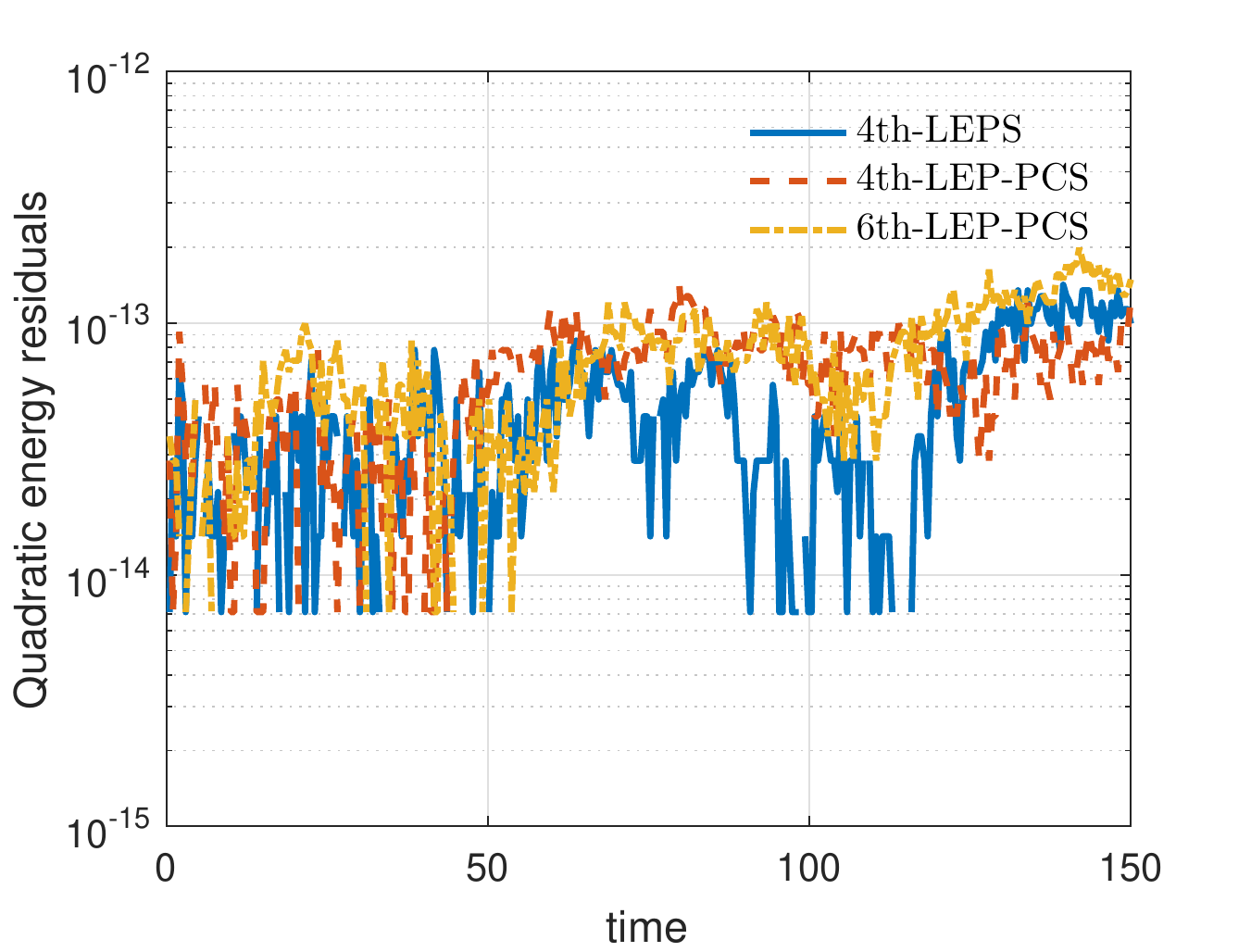}
%\caption*{(d) Quadratic energy}
\end{minipage}
\caption{{The residuals on discrete conservation laws using the different numerical schemes with the time step $\tau=0.01$ and the Fourier node 1024, respectively, for the one-dimensional RLW equation \eqref{RLW-1d-equation}.}}\label{RLW1d-fig-2}
\end{figure}

We take the computational domain $\Omega=[-60,300]$ with a periodic boundary condition and take parameters $\alpha=1,\ \mu=1,\ x_1=-20,\ x_2=15, c_1=1,$ and $c_2=0.5$, respectively.
In the implementation, we choose the time step $\tau=0.01$ and the Fourier node 1024, respectively. The evolution of the solitary waves is summarized Fig. \ref{RLW1d-fig-1}. We observe that an interaction occurs when the larger one is catching up with the smaller one and eventually passes it. We point out that the profile of $u$
calculated using other schemes are similar to Fig. \ref{RLW1d-fig-1}, thus for brevity, we omit them here. In addition, we further visualize the mass, momentum, Hamiltonian energy and quadratic energy residuals (i.e., $|M^n-M^0|,|I^n-I^0|,|H^n-H^0|$ and $|E^n-E^0|$, hereafter), which are summarized in Fig. \ref{RLW1d-fig-2}. We
can draw the following observations: (i) 4th-LEPS, 4th-LEP-PCS and 6th-LEP-PCS preserve the discrete mass and quadratic energy exactly, and they cannot preserve the discrete momentum exactly but the residuals reduce as the accuracy of the scheme is improved; (ii) LCN-MPS, 4th-LMPS, 4th-LMP-PCS and 6th-LMP-PCS preserve discrete momentum up to the machine precision and conserve the discrete mass approximately; (iii) LEPS preserves the discrete mass and all of the proposed schemes cannot preserve the discrete Hamiltonian energy exactly.

{Moreover, inspired by \cite{Frasca-Caccia-AMC-2021}, we treat the numerical solution of 6th-LMP-PCS on a finer grid: $\tau=0.001$ and Fourier node 2048 at $T=30$ as the reference exact solution and some comparisons with the solutions of the different schemes on a much coarser grid: $\tau=0.1$ and Fourier node 1024 at $T=30$ are summarized in Tab. \ref{Tab:RLW-equation:2}. It is clear to observe: (i) the momentum errors of LCN-MPS, 4th-LMPS, 4th-LMP-PCS and 6th-LMP-PCS are preserved to machine accuracy and LEPS, 4th-LEPS, 4th-LEP-PCS and 6th-LEP-PCS can preserve the discrete mass exactly; (ii) all of the schemes cannot preserve the discrete Hamiltonian energy exactly and the energy error provided by 6th-LMP-PCS is the smallest; (iii) the solution error of 6th-LEP-PCS is the smallest and the one provided by LEPS is the largest.}

\begin{table}[H]
\tabcolsep=9pt
\footnotesize
\renewcommand\arraystretch{1.1}
\centering
\caption{{The interaction of two solitary waves of the RLW equation: the time step $\tau=0.1$ and the Fourier node 1024 at $T=30$.}}\label{Tab:RLW-equation:2}
\begin{tabular*}{\textwidth}[h]{@{\extracolsep{\fill}} c c c c c c}\hline
{Scheme\ \ } &{$|M^n-M^0|$} &{$|H^n-H^0|$} &{$|I^n-I^0|$}& {$e_2(t_n=30)$}  \\[2ex]     %% 第1 行
\hline
 {LCN-MPS}  &{2.564e-02}& {5.256e-05}&{3.411e-13} &{1.004e-02} \\[2ex]%&{2.8}\\
  {LEPS}  &{3.553e-14}& {7.590e-04}&{6.679e-04} &{1.741e-02}\\[2ex] %&{4.3}  \\   %% 第2行
   {4th-LMPS}  &{7.313e-05}& {1.630e-08}&{2.132e-14} &{1.318e-04} \\[2ex]%\hline
    {4th-LEPS}  &{1.776e-14}& {1.316e-04}&{1.325e-04} &{2.462e-04} \\[2ex]
  {4th-LMP-PCS}  &{1.067e-04}& {5.684e-09}&{2.132e-14} &{1.092e-04} \\[2ex] %&{4.3}  \\   %% 第2行
   {4th-LEP-PCS}  &{2.487e-14}& {2.068e-05}&{2.092e-05} &{8.816e-05} \\[2ex]
   {6th-LMP-PCS}  &{6.968e-08}& {1.378e-12}&{5.684e-14} &{2.947e-08}\\[2ex]
   {6th-LEP-PCS}  &{3.553e-15}& {1.060e-08}&{1.065e-08} &{2.727e-08} \\[2ex]\hline
\end{tabular*}
\end{table}

%{\bf Example 3:} Finally,
%\begin{align}
%u(x,0)=\exp(-(x-7)^2), \ x\in\mathbb{R},
%\end{align}
%
%
%
%\begin{figure}[H]
%\centering
%\begin{minipage}[t]{70mm}
%\includegraphics[width=80mm]{1d_RWL_fig2.eps}
%\end{minipage}
%\caption{ The profile of $u$ provided by 4th-LEPS over time interval $[0,150]$ with the time step $\tau=0.1$ and the Fourier node 1024, respectively, for the RLW equation in 1D \eqref{RLW-1d-equation}.}\label{RLW1d-fig-1}
%\end{figure}
%
%
%\begin{figure}[H]
%\centering\begin{minipage}[t]{60mm}
%\includegraphics[width=60mm]{RLW1d_example3_M.eps}
%\end{minipage}
%\caption{(a) Mass}
%\begin{minipage}[t]{60mm}
%\includegraphics[width=60mm]{RLW1d_example3_Mom.eps}
%\end{minipage}
%\centering\begin{minipage}[t]{60mm}
%\includegraphics[width=60mm]{RLW1d_example3_H.eps}
%\end{minipage}
%\begin{minipage}[t]{60mm}
%\includegraphics[width=60mm]{RLW1d_example3_E.eps}
%\end{minipage}
%\caption{ The residuals on discrete conservation laws using the different numerical schemes with the time step $\tau=0.1$ and the Fourier node 1024, respectively, for the RLW equation in 1D \eqref{RLW-1d-equation}.}\label{RLW1d-fig-2}
%\end{figure}

\subsection{RLW equation in 2D}
In this subsection,  we firstly investigate temporal accuracies and computational efficiencies of the proposed schemes for two-dimensional RLW equation \eqref{RLW-equation}, respectively.
We consider the computational domain $\Omega=[0,2\pi]^2$ and fix the parameters $\alpha=1, \beta=1, \theta=1$ and $\mu=1$.
Although the two-dimensional RLW equation \eqref{RLW-equation} admits soliton solution \cite{Huang-AMM-2002}, we shall point out that even a large enough computational domain cannot prevent the periodic
boundary condition from introducing a significant (aliasing) error relative to the whole
space problem. Thus, we choose the initial profiles as
\begin{align}
u(x,y,0)=(1+\sin(x))(1+\sin(y)),
\end{align}
and the exact solution at $T=10$ is generated by 6th-LEP-PCS with the time step $\tau=0.001$ and the Fourier node $128\times 128$, respectively.

The $l^2$ and $l^{\infty}$ -errors and convergence rates with various time steps and the Fourier node $128\times 128$ at $T=10$ are summarized in Tab. \ref{Tab:2dRLW-equation:1}. As shown in the Table, we observe that (i) 4th-LMPS, 4th-LEPS, 4th-LMP-PCS and 4th-LEP-PCS can achieve fourth-order accurate in time, and 6th-LMP-PCS and 6th-LEP-PCS are sixth-order accurate in time; (ii) the error provided by 6th-LEP-PCS is the smallest and the error provided by 4th-LEPS is the largest; (iii) 4th-LEP-PCS admits much smaller numerical error than 4th-LMP-PCS.

{Moreover, some comparisons with the solutions of the different schemes on a much coarser grid: the time step $\tau=0.2$ and the Fourier node $64\times 64$ at $T=50$ are given in Tab. \ref{Tab:2dRLW-equation:2}, where the exact solution at $T=50$ is obtained by 6th-LEP-PCS on a finer grid: $\tau=0.001$ and the Fourier node $128\times 128$. It is clear to see that (i) 4th-LEPS, 4th-LEP-PCS and 6th-LEP-PCS can preserve the discrete mass exactly and the mass errors provided by the momentum-preserving schemes are reduced as the accuracy of the scheme is improved; (ii) 4th-LMPS, 4th-LMP-PCS and 6th-LMP-PCS conserve the discrete momentum (up to rounding error); (iii) the solution error of 6th-LEP-PCS is the smallest and the one provided by 4th-LMPS is the largest.}

\begin{figure}[H]
\centering\begin{minipage}[t]{60mm}
\includegraphics[width=60mm]{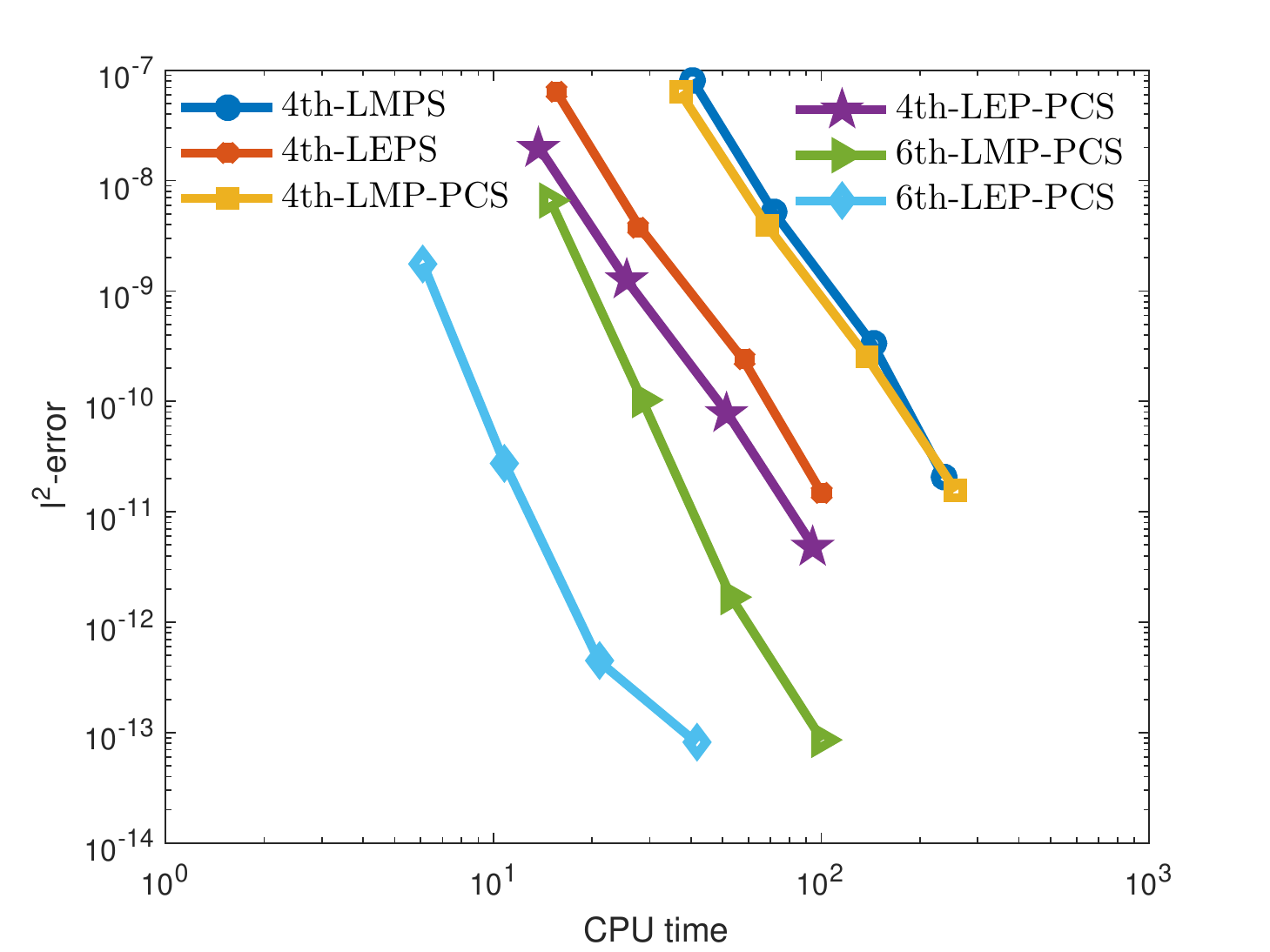}
\end{minipage}
\begin{minipage}[t]{60mm}
\includegraphics[width=60mm]{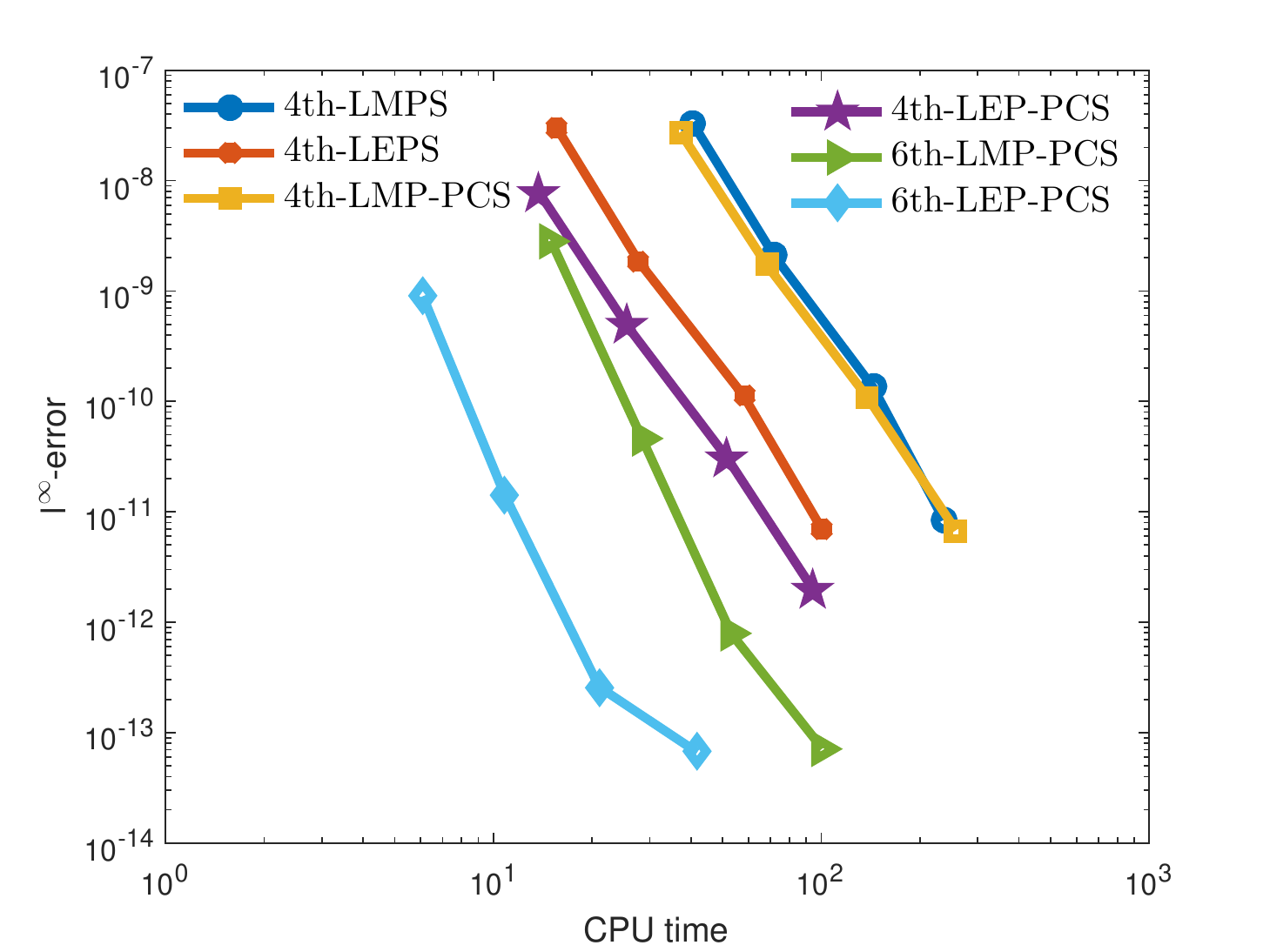}
\end{minipage}
%\centering\begin{minipage}[t]{60mm}
%\includegraphics[width=60mm]{1d_NLS_fig4.eps}
%\end{minipage}
%\begin{minipage}[t]{60mm}
%\includegraphics[width=60mm]{1d_NLS_fig5.eps}
%\end{minipage}
\caption{ The numerical error versus the CPU time using the different numerical schemes with various time steps and the Fourier node $128\times 128$ at $T=10$ for the two-dimensional RLW equation  \eqref{RLW-equation}.}\label{RLW2d-CPU}
\end{figure}

In Fig. \ref{RLW2d-CPU}, we display the global $l^2$ and $l^{\infty}$ -errors (i.e., $e_2(t_n)$ and $e_{\infty}(t_n)$) of $u$ versus the CPU time using the proposed schemes with the Fourier node $128\times 128$ and various time steps at $T=10$, respectively. For a given global error, we can observe that (i) the cost of 6th-LEP-PC is the cheapest and the cost of 4th-LMPS is most expensive; (ii) the costs of 4th-LEP-PCS and 4th-LEPS are much cheaper than the ones provided by 4th-LMP-PCS and 4th-LMPS, respectively; (iii) the cost provided by the linear energy-preserving scheme is much cheaper than the linear momentum-preserving one of the same order.

\begin{table}[H]
\tabcolsep=9pt
\footnotesize
\renewcommand\arraystretch{1.1}
\centering
\caption{{Numerical errors and convergence rates for the different schemes with various time steps and the Fourier node $128\times 128$ at $T=10$.}}\label{Tab:2dRLW-equation:1}
\begin{tabular*}{\textwidth}[h]{@{\extracolsep{\fill}} c c c c c c}\hline
{Scheme\ \ } &{$\tau$}&{$e_{2}(t_n=10)$} &{order}& {$e_{\infty}(t_n=10)$}&{order}  \\     %% 第1 行
\hline
% \multirow{4}{*}{3rd-LMPS}  &{$\frac{1}{100}$}& {4.517e-05}&{-} &{2.362e-05} & {-}\\[1ex]%&{2.8}\\
%  {}  &{$\frac{1}{200}$}& {5.652e-06}&{2.998} &{2.968e-06} & {2.992}\\[1ex] %&{4.3}  \\   %% 第2行
%   {}  &{$\frac{1}{400}$}& {7.068e-07}&{2.999} &{3.720e-07} &{2.996} \\[1ex]%\hline
%    {}  &{$\frac{1}{800}$}& {8.838e-08}&{3.000} &{4.657e-08} &{2.998} \\
% \multirow{4}{*}{3rd-LEPS}  &{$\frac{1}{100}$}& {4.072e-05}&{-} &{2.763e-05} & {-}\\[1ex]%&{2.8}\\
%  {}  &{$\frac{1}{200}$}& {5.103e-06}&{2.996} &{3.457e-06} & {2.998}\\[1ex] %&{4.3}  \\   %% 第2行
%   {}  &{$\frac{1}{400}$}& {6.388e-07}&{2.998} &{4.323e-07} &{2.999} \\[1ex]%\hline
%    {}  &{$\frac{1}{800}$}& {7.991e-08}&{2.999} &{5.406e-07} &{3.000} \\\hline
    \multirow{4}{*}{4th-LMPS}  &{$\frac{1}{10}$}& {5.190e-04}&{-} &{2.228e-04} & {-}\\[1ex]%&{2.8}\\
  {}  &{$\frac{1}{20}$}& {4.308e-05}&{3.590} &{1.808e-05} & {3.623}\\[1ex] %&{4.3}  \\   %% 第2行
   {}  &{$\frac{1}{40}$}& {3.009e-06}&{3.839} &{1.246e-06} &{3.859} \\[1ex]%\hline
    {}  &{$\frac{1}{80}$}& {1.977e-07}&{3.928} &{8.133e-08} &{3.937} \\
     \multirow{4}{*}{4th-LEPS}  &{$\frac{1}{10}$}& {1.022e-03}&{-} &{3.564e-04} & {-}\\[1ex]%&{2.8}\\
  {}  &{$\frac{1}{20}$}& {6.140e-05}&{4.057} &{2.215e-05} & {4.008}\\[1ex] %&{4.3}  \\   %% 第2行
   {}  &{$\frac{1}{40}$}& {3.796e-06}&{ 4.016} &{1.394e-06} &{3.990} \\[1ex]%\hline
    {}  &{$\frac{1}{80}$}& {2.364e-07}&{4.005} &{8.759e-08} &{3.992} \\%\hline
    \multirow{4}{*}{4th-LMP-PCS}  &{$\frac{1}{10}$}& {5.648e-04}&{-} &{2.442e-04} & {-}\\[1ex]%&{2.8}\\
  {}  &{$\frac{1}{20}$}& {3.761e-05}&{3.909} &{1.637e-05} & {3.898}\\[1ex] %&{4.3}  \\   %% 第2行
   {}  &{$\frac{1}{40}$}& {2.418e-06}&{3.959} &{1.056e-06} &{3.955} \\[1ex]%\hline
    {}  &{$\frac{1}{80}$}& {1.531e-07}&{3.981} &{6.697e-08} &{ 3.979} \\
    \multirow{4}{*}{4th-LEP-PCS}  &{$\frac{1}{10}$}& {2.614e-04}&{-} &{1.128e-04} & {-}\\[1ex]%&{2.8}\\
  {}  &{$\frac{1}{20}$}& {1.531e-05}&{4.094} &{6.778e-06} & {4.057}\\[1ex] %&{4.3}  \\   %% 第2行
   {}  &{$\frac{1}{40}$}& {9.248e-07}&{4.049} &{4.144e-07} &{4.032} \\[1ex]%\hline
    {}  &{$\frac{1}{80}$}& {5.681e-08}&{4.025} &{ 2.560e-08} &{4.017} \\\hline%\hline
    \multirow{4}{*}{6th-LMP-PCS}  &{$\frac{1}{10}$}& {3.828e-07}&{-} &{1.694e-07} & {-}\\[1ex]%&{2.8}\\
  {}  &{$\frac{1}{20}$}& {6.461e-09}&{5.889} &{ 2.872e-09} & {5.882}\\[1ex] %&{4.3}  \\   %% 第2行
   {}  &{$\frac{1}{40}$}& {1.044e-10}&{5.952} &{4.646e-11} &{5.950} \\[1ex]%\hline
    {}  &{$\frac{1}{80}$}& {1.664e-12}&{5.973} &{7.847e-13} &{5.889} \\
    \multirow{4}{*}{6th-LEP-PCS}  &{$\frac{1}{10}$}& {9.059e-08}&{-} &{3.474e-08} & {-}\\[1ex]%&{2.8}\\
  {}  &{$\frac{1}{20}$}& {1.501e-09}&{5.915} &{5.829e-10} & {5.897}\\[1ex] %&{4.3}  \\   %% 第2行
   {}  &{$\frac{1}{40}$}& {2.447e-11}&{5.939} &{9.413e-12} &{5.953} \\[1ex]%\hline
    {}  &{$\frac{1}{80}$}& {3.969e-13}&{5.946} &{1.819e-13} &{5.697} \\\hline%\hline

\end{tabular*}
\end{table}

{\color{red}\begin{table}[H]
\tabcolsep=9pt
\footnotesize
\renewcommand\arraystretch{1.1}
\centering
\caption{{The numerical solution of the RLW equation in 2D: $\tau=0.2$ and Fourier node $64\times 64$ at $T=50$.}}\label{Tab:2dRLW-equation:2}
\begin{tabular*}{\textwidth}[h]{@{\extracolsep{\fill}} c c c c c c}\hline
{Scheme\ \ } &{$|M^n-M^0|$} &{$|H^n-H^0|$} &{$|I^n-I^0|$}& {$e_2(t_n=50)$}  \\[2ex]     %% 第1 行
\hline
 %{LCN-MPS}  &{2.552e-05}& {1.365e-08}&{4.832e-13} &{1.286e-04} \\[2ex]%&{2.8}\\
%  {LEPS}  &{3.552e-14}& { 5.843e-07}&{6.986e-07} &{1.728e-04} \\[2ex] %&{4.3}  \\   %% 第2行
   {4th-LMPS}  &{1.423e-01}& {1.189e-01}&{8.527e-14} &{2.341e-01}  \\[2ex]%\hline
    {4th-LEPS}  &{1.634e-13}& {1.096e-01}&{1.584e-01} &{1.279e-01}  \\[2ex]
  {4th-LMP-PCS}  &{2.672e-02}& {2.310e-02}&{2.842e-14} &{8.039e-02} \\[2ex] %&{4.3}  \\   %% 第2行
   {4th-LEP-PCS}  &{7.105e-15}& {4.111e-02}&{6.874e-02} &{1.426e-02} \\[2ex]
   {6th-LMP-PCS}  &{1.637e-04}& {1.411e-04}&{5.684e-14} &{3.567e-04} \\[2ex]
   {6th-LEP-PCS}  &{2.132e-14}& {2.446e-04}&{4.133e-04} &{4.572e-05} \\[2ex]\hline
\end{tabular*}
\end{table}}

\begin{figure}[H]
\centering\begin{minipage}[t]{60mm}
\includegraphics[width=60mm]{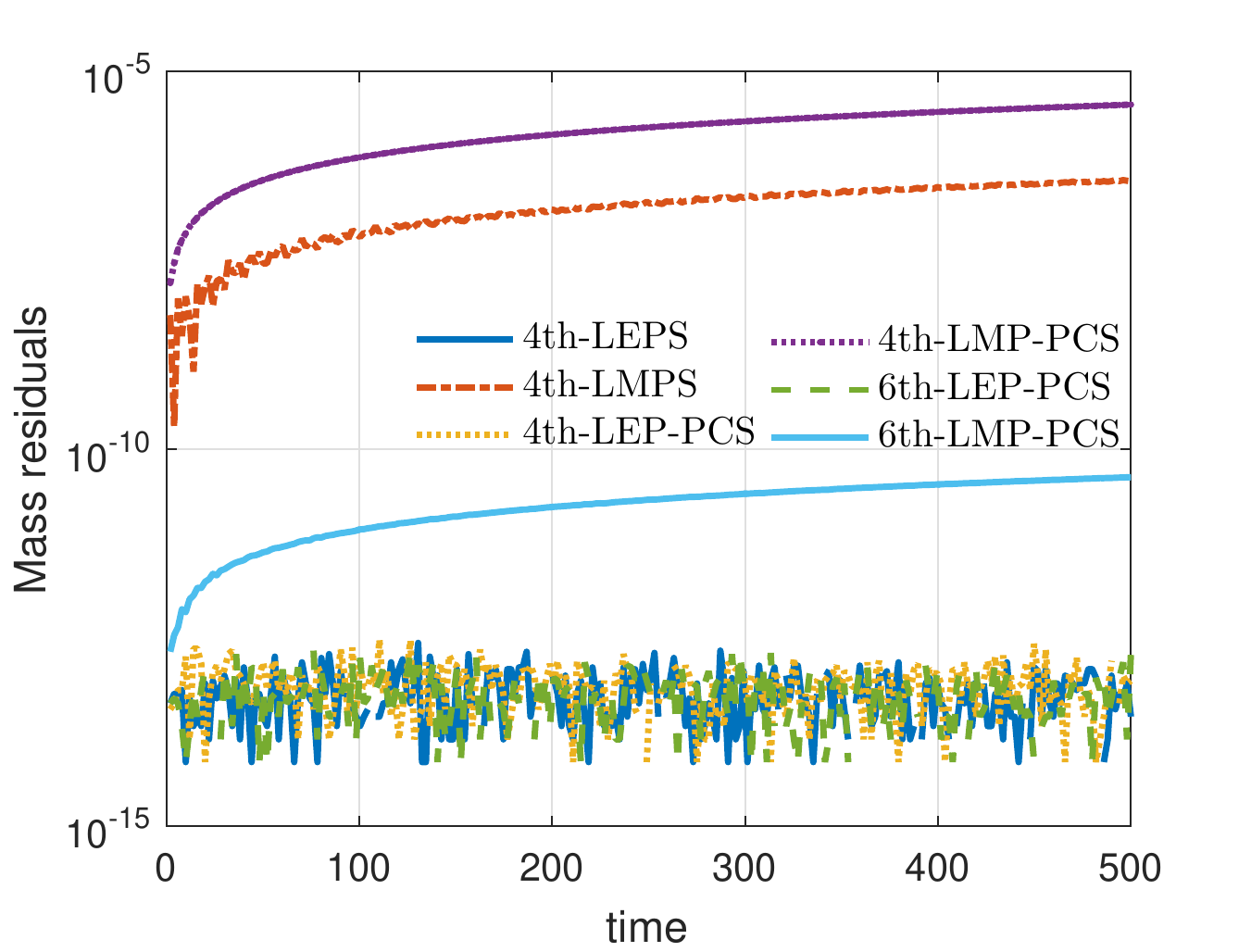}
%\caption*{(a) Mass}
\end{minipage}
\begin{minipage}[t]{60mm}
\includegraphics[width=60mm]{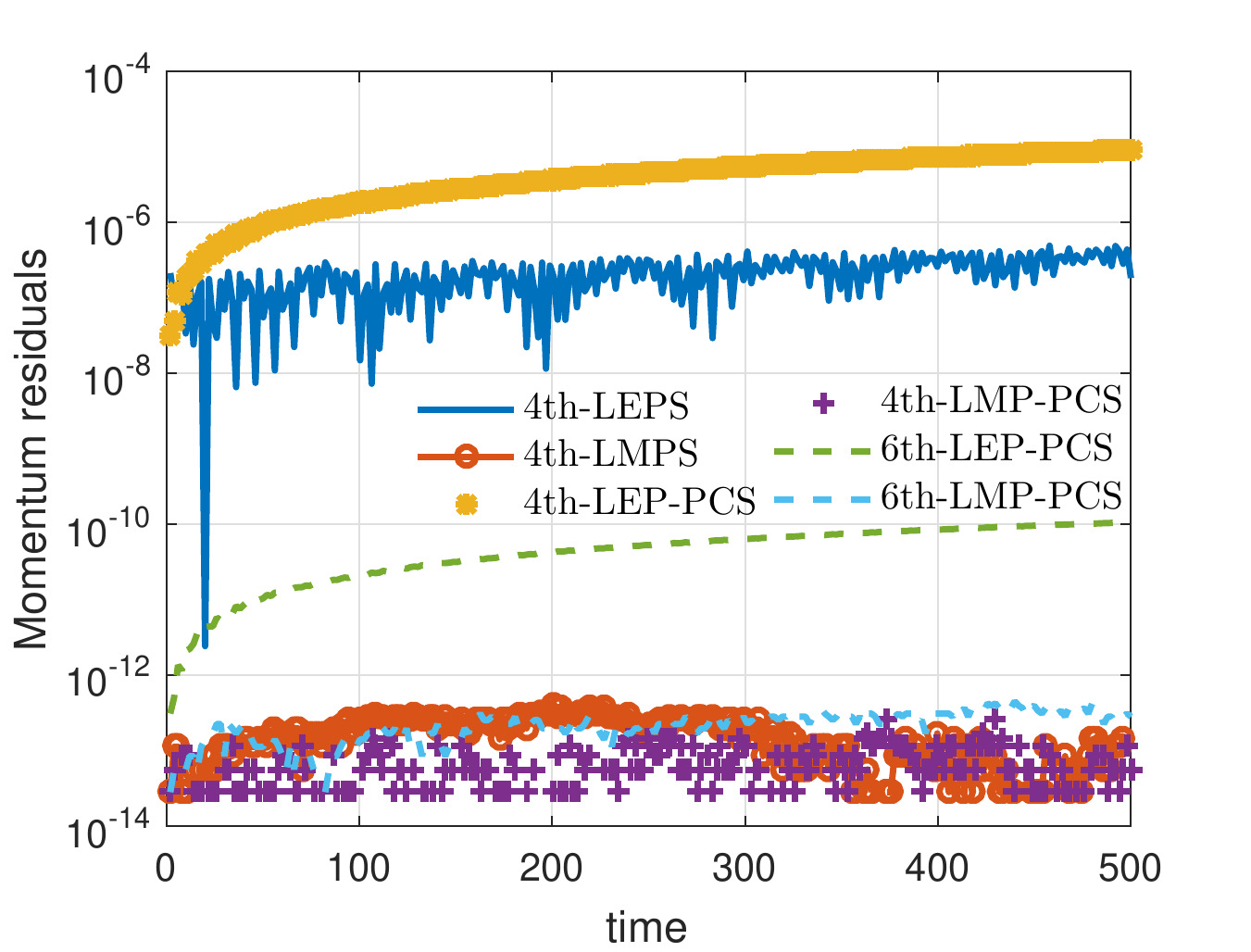}
%\caption*{(b) Momentum}
\end{minipage}
\centering\begin{minipage}[t]{60mm}
\includegraphics[width=60mm]{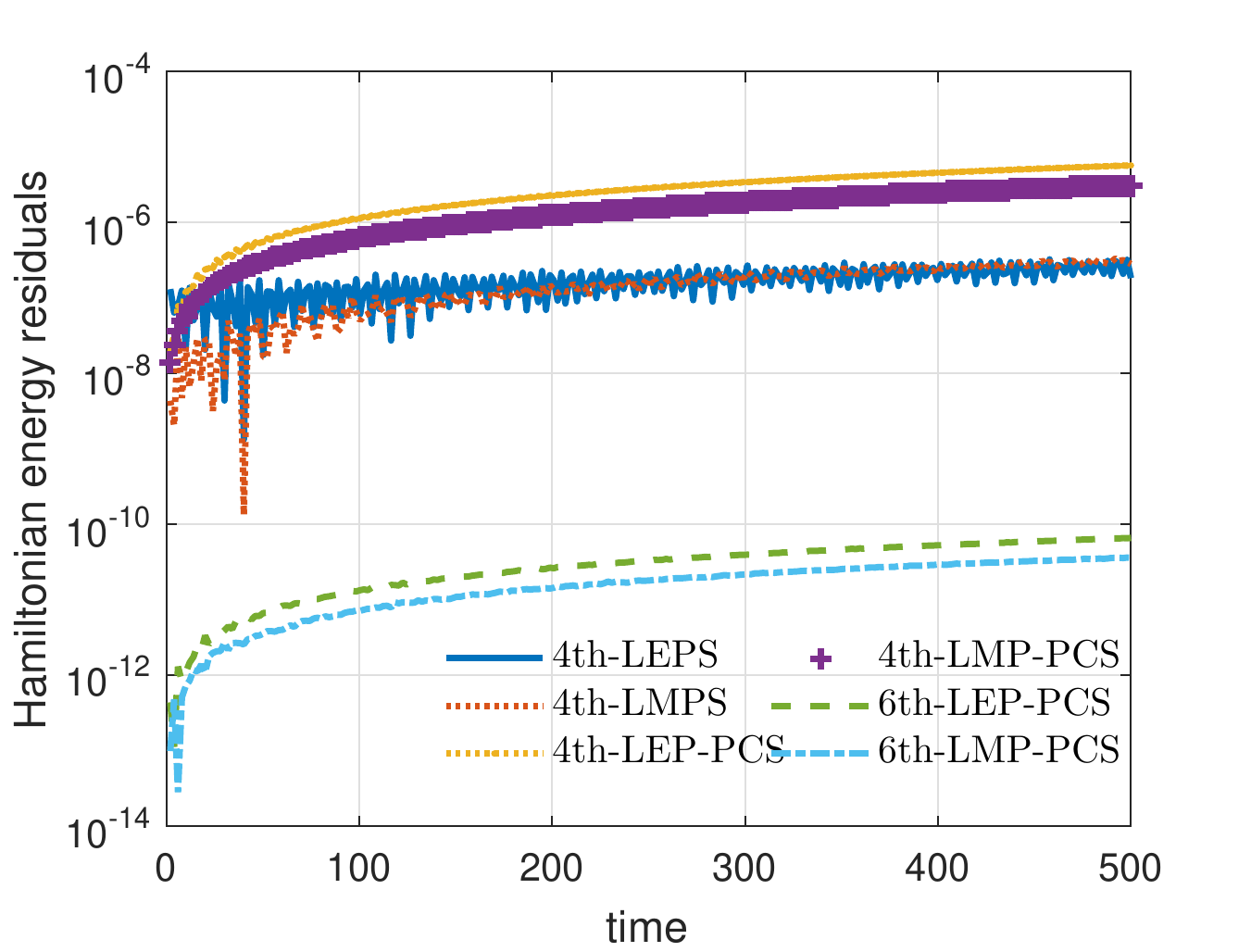}
%\caption*{(c) Hamiltonian energy}
\end{minipage}
\begin{minipage}[t]{60mm}
\includegraphics[width=60mm]{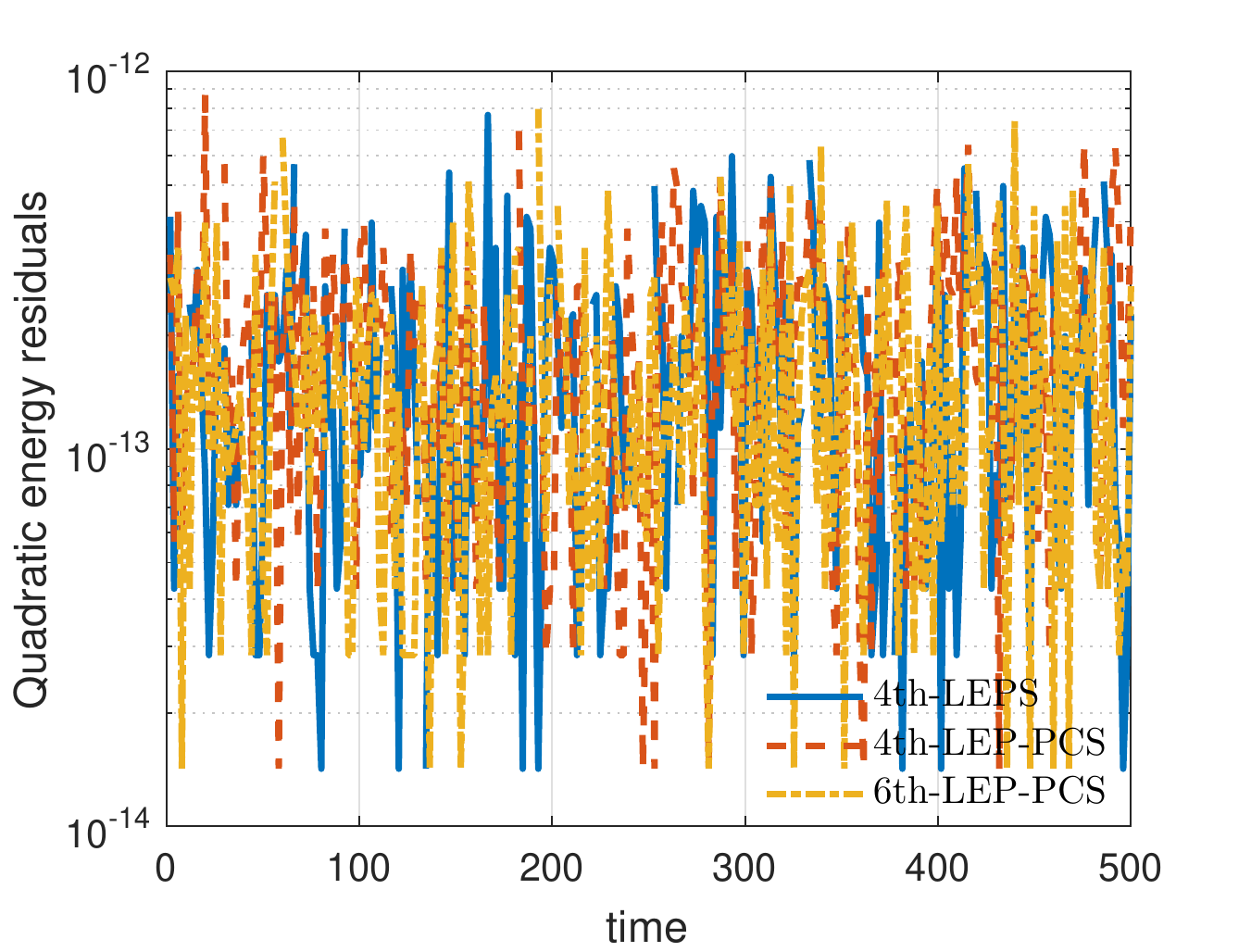}
%\caption*{(d) Quadratic energy}
\end{minipage}
\caption{ {The residuals on discrete conservation laws using the different numerical schemes with the time step $\tau=0.01$ and the Fourier node $64\times 64$, respectively, for the two-dimensional RLW equation \eqref{RLW-equation}.}}\label{RLW2d-example1-fig1}
\end{figure}

Meanwhile, the residuals on discrete conservation laws from $t=0$ and $t=500$ using the different numerical schemes are summarized in Fig. \ref{RLW2d-example1-fig1}.
 We observe that 4th-LMPS, 4th-LMP-PCS and 6th-LMP-PCS can preserve the discrete momentum exactly, and 4th-LEPS, 4th-LEP-PCS and 6th-LEP-PCS preserve the discrete mass and quadratic energy (i.e., the modified energy \eqref{RLW-Q-energy}) exactly. In addition, it is clear to see that all schemes can preserve the Hamiltonian energy well.

\begin{figure}[H]
\centering\begin{minipage}[t]{60mm}
\includegraphics[width=60mm]{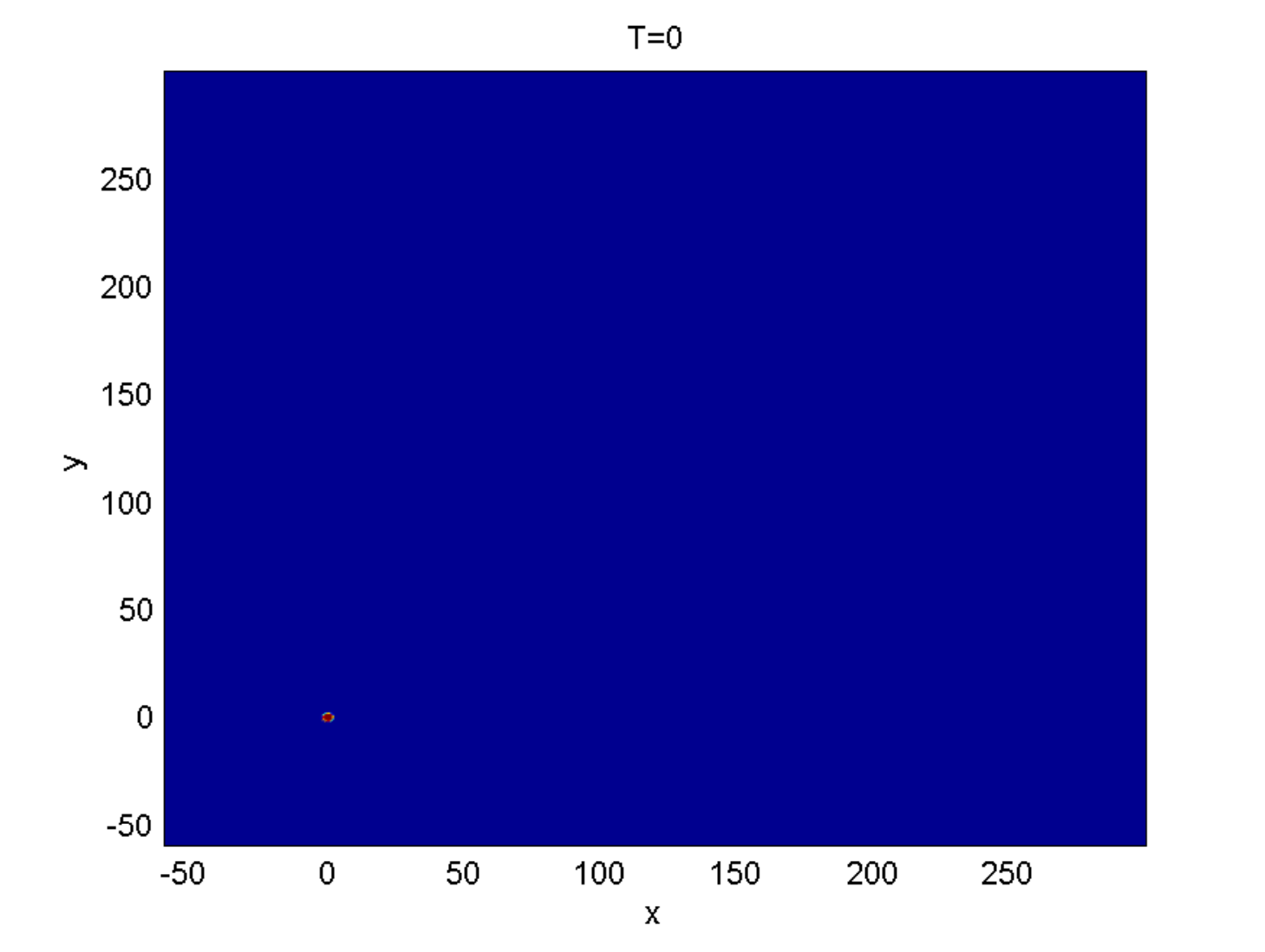}
\end{minipage}
\begin{minipage}[t]{60mm}
\includegraphics[width=60mm]{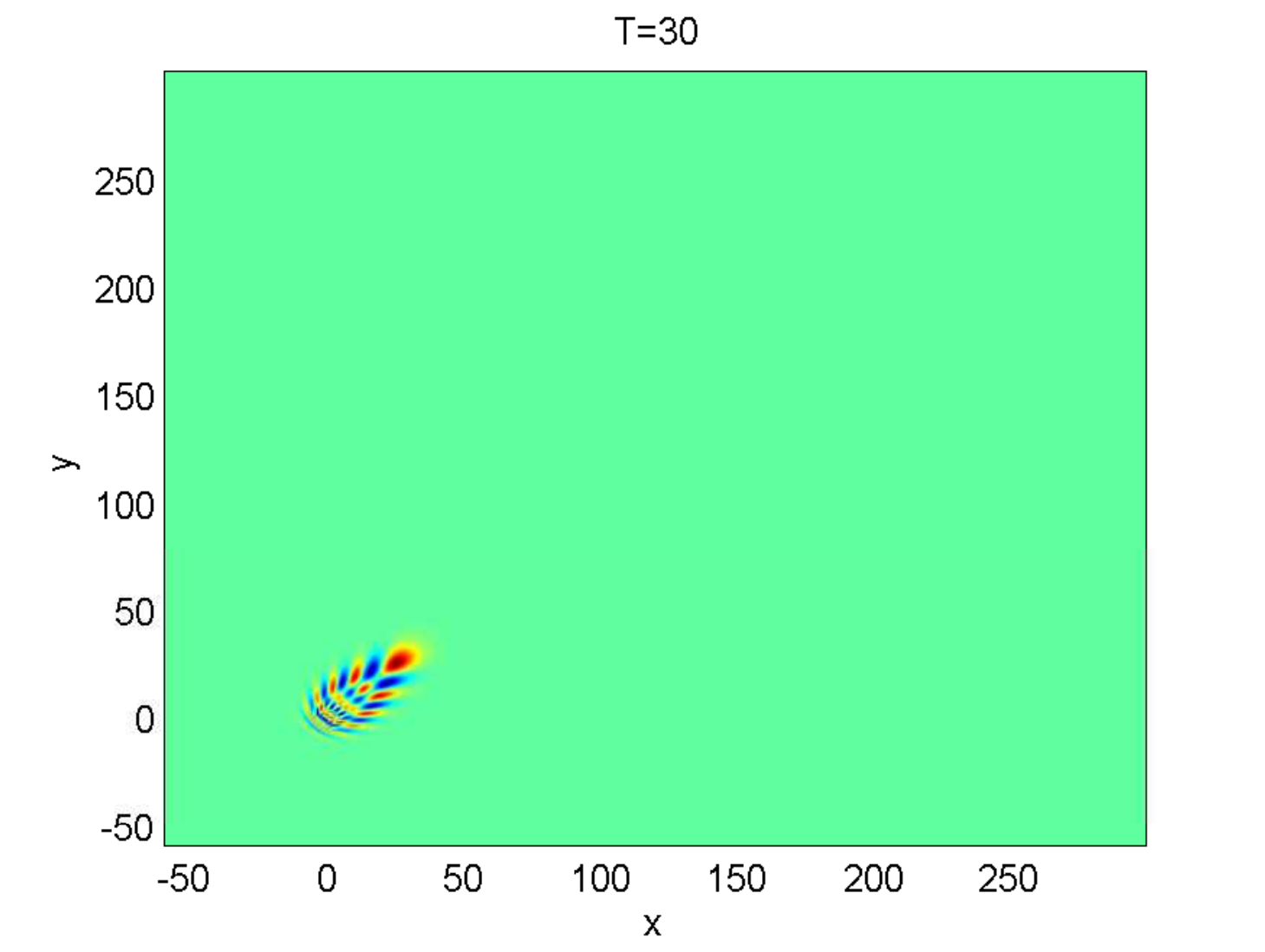}
\end{minipage}
\centering\begin{minipage}[t]{60mm}
\includegraphics[width=60mm]{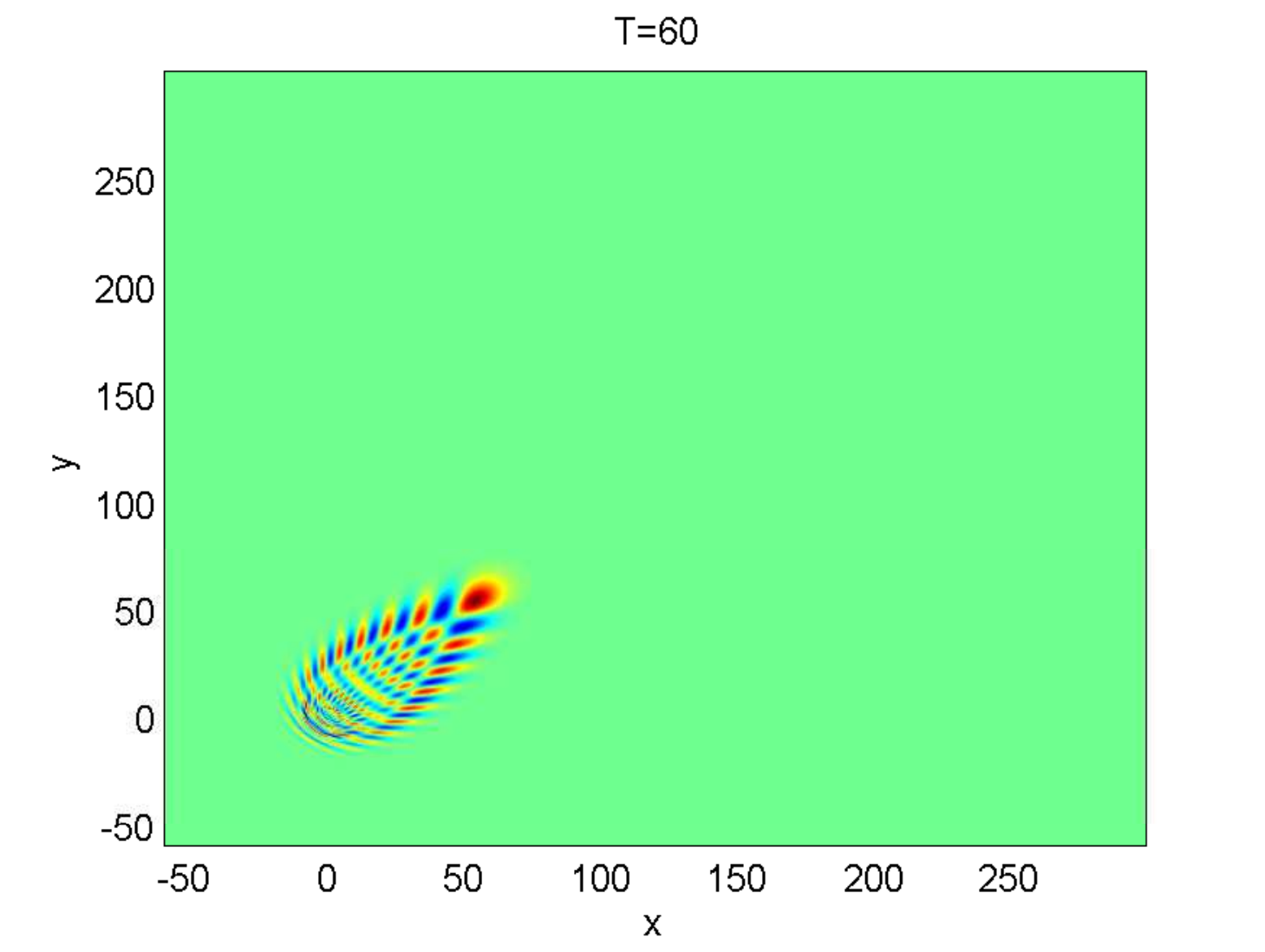}
\end{minipage}
\begin{minipage}[t]{60mm}
\includegraphics[width=60mm]{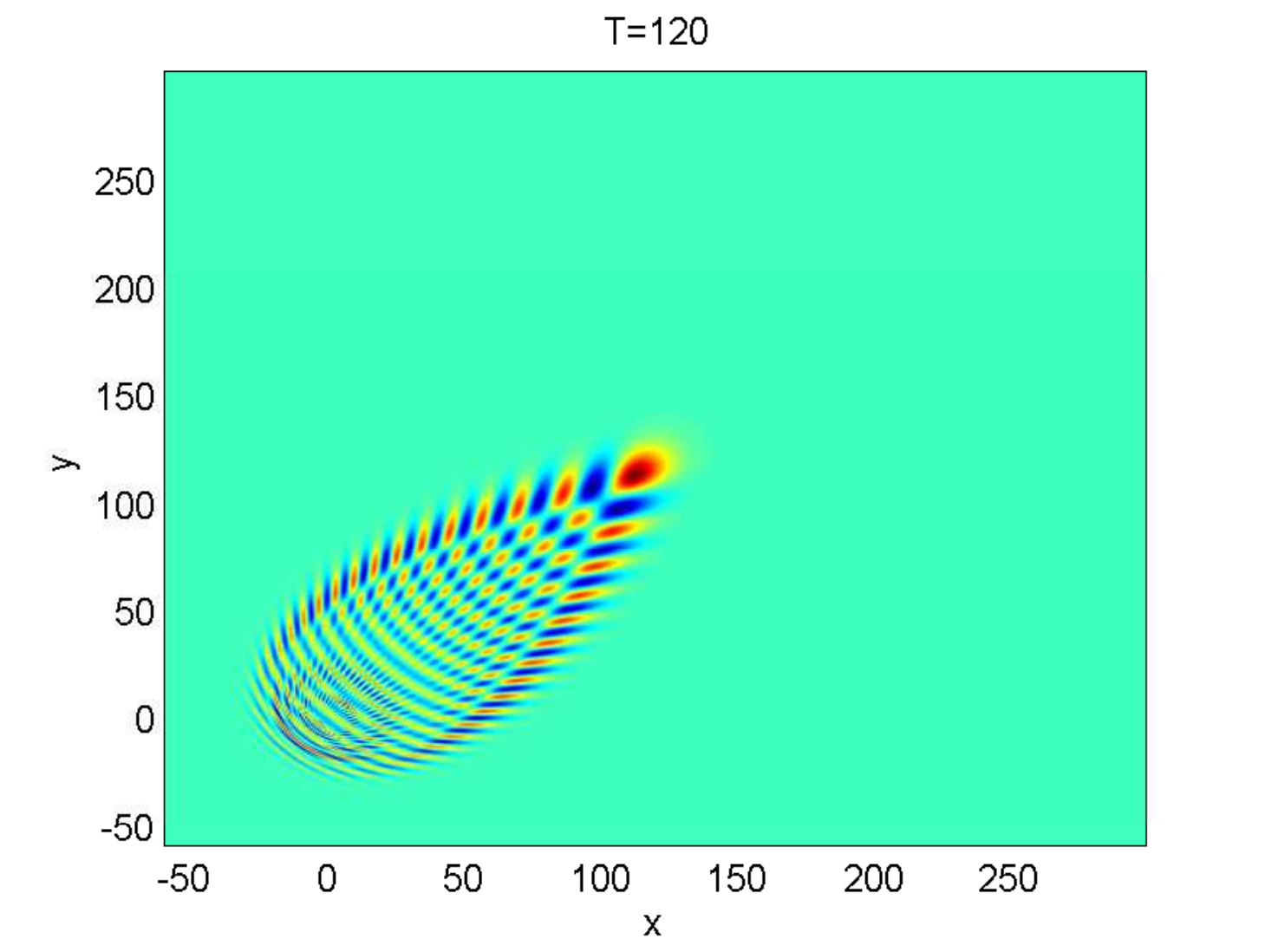}
\end{minipage}
\centering\begin{minipage}[t]{60mm}
\includegraphics[width=60mm]{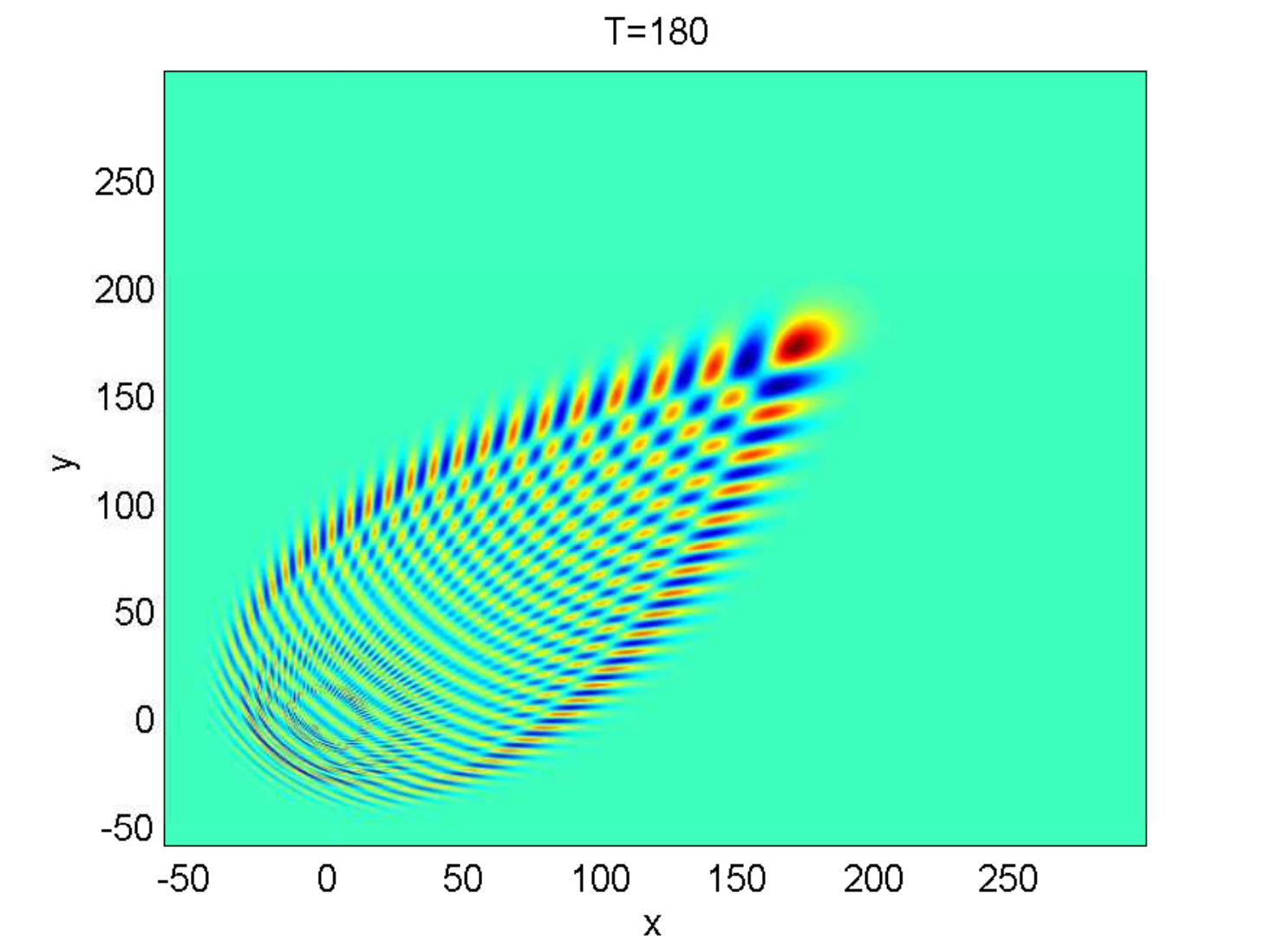}
\end{minipage}
\begin{minipage}[t]{60mm}
\includegraphics[width=60mm]{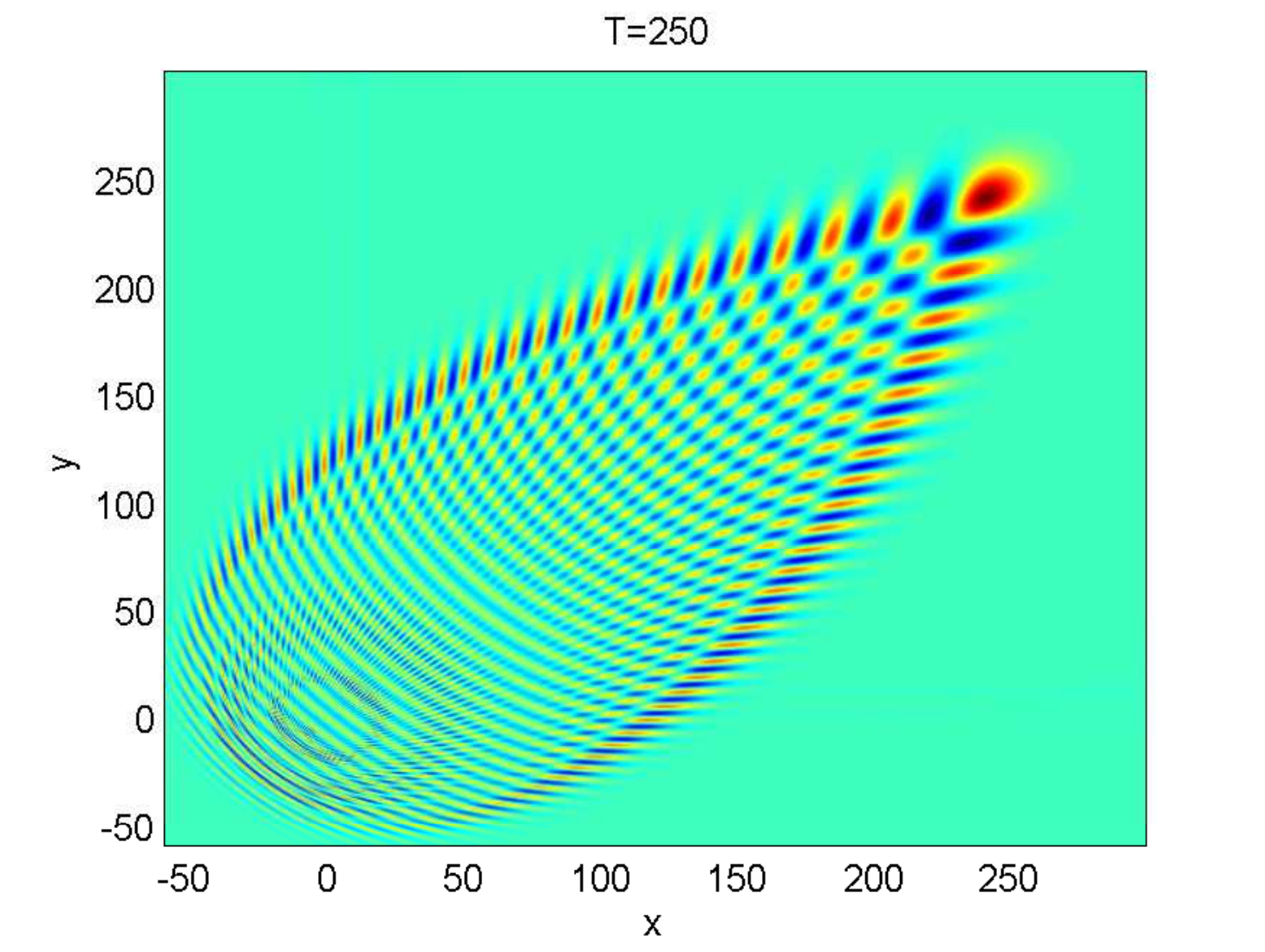}
\end{minipage}
\caption{ The profile of $u$ provided by 6th-LMP-PCS at times $t=0,30,60,120,180$ and 250 with the time step $\tau=0.1$ and the Fourier node $512\times 512$, respectively, for the two-dimensional RLW equation \eqref{RLW-equation}.}\label{RLW2d-example2-fig1}
\end{figure}

Next, we consider the numerical evolution of an undular bore for the two-dimensional RLW equation with the initial condition \cite{LHQjsc2018}
\begin{align}
u(x,y,0)=0.05\Big(1-\tanh\big((x-x_0)^2+(y-y_0)^2-d^2\big)\Big)
\end{align}
with $x_0=0, y_0=0$ and $d=2$.  We set computational domain $\Omega=[-60,300]^2$ with a periodic boundary condition, take parameters $\alpha=1, \beta=1, \theta=1$ and $\mu=1$, and use the time step $\tau=0.1$ and the Fourier node $512\times 512$, respectively. The profile of $u$ at different times are shown in Fig. \ref{RLW2d-example2-fig1}. It is clear to see that the 2D undular bore expands and propagates in a certain direction. This fact
is consistent with the results obtained in \cite{LHQjsc2018}. We should note that the profile of $u$ at different times calculated using other
schemes are similar to Fig. \ref{RLW2d-example2-fig1}, thus for brevity, we omit them here. Meanwhile, the residuals on discrete conservation laws from $t=0$ to $t=250$ are summarized in Fig. \ref{RLW2d-example2-fig2}, which agrees with the theoretical analysis.

\begin{figure}[H]
\centering\begin{minipage}[t]{60mm}
\includegraphics[width=60mm]{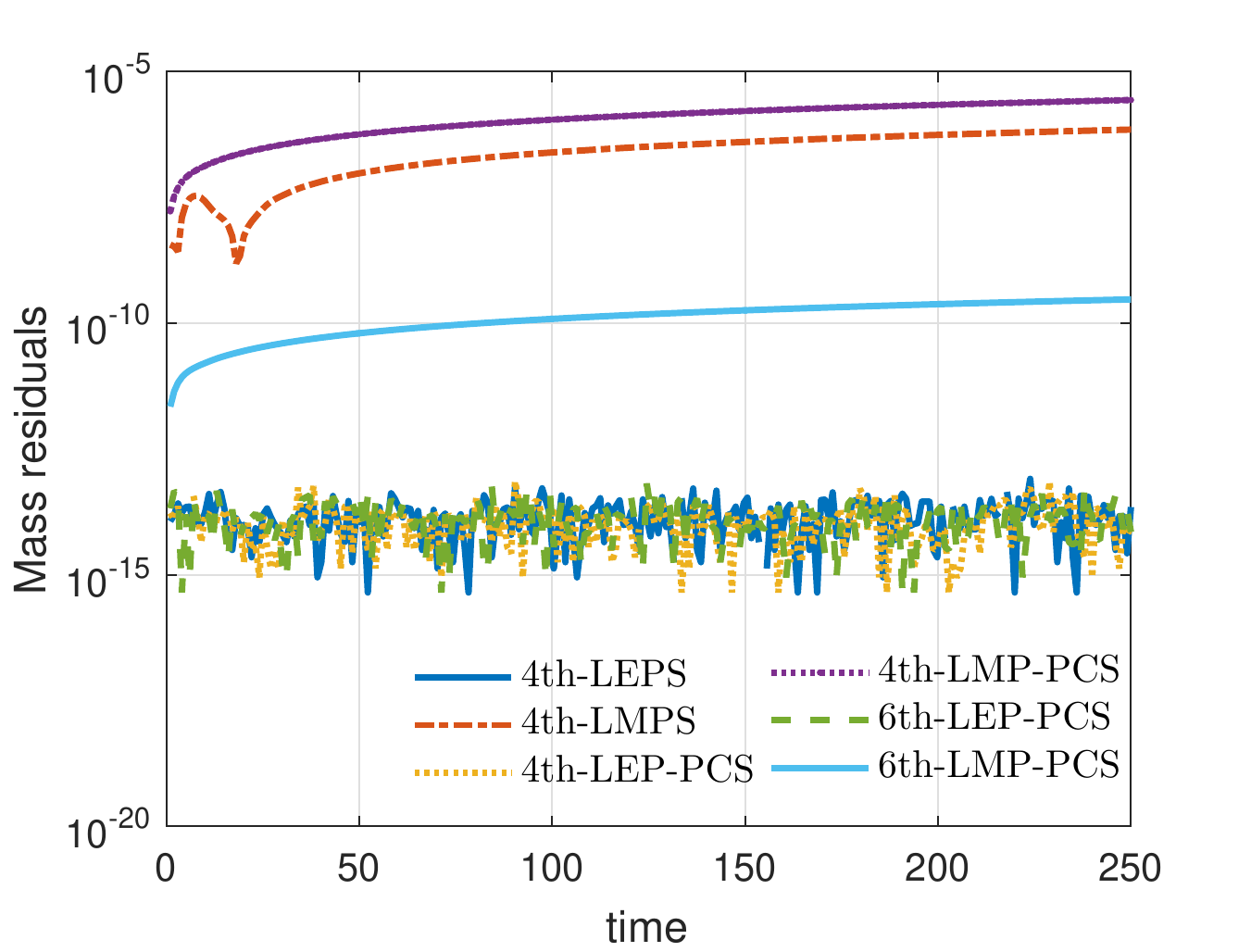}
%\caption*{(a) Mass}
\end{minipage}
\begin{minipage}[t]{60mm}
\includegraphics[width=60mm]{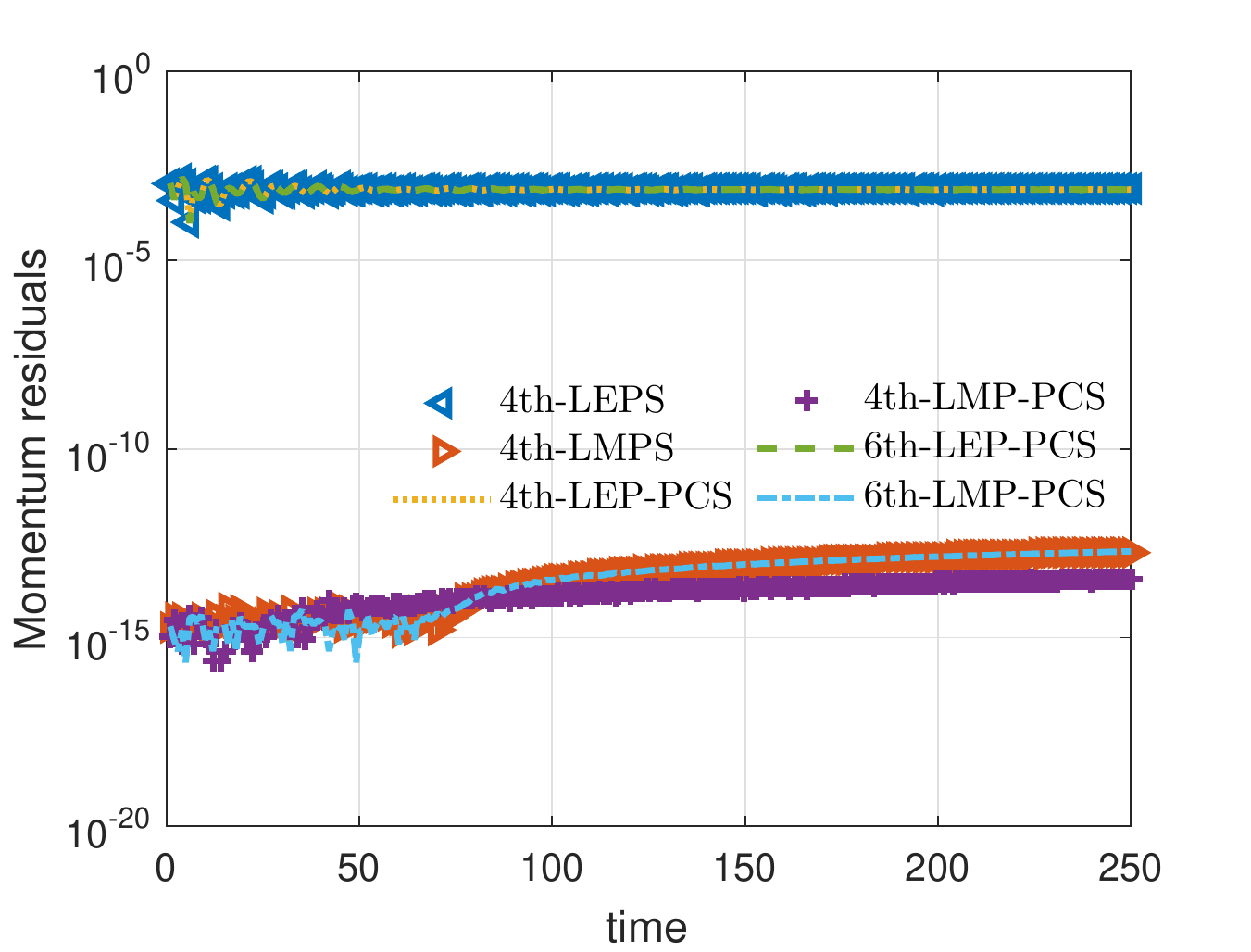}
%\caption*{(b) Momentum}
\end{minipage}
\end{figure}
\begin{figure}[H]
\centering\begin{minipage}[t]{60mm}
\includegraphics[width=60mm]{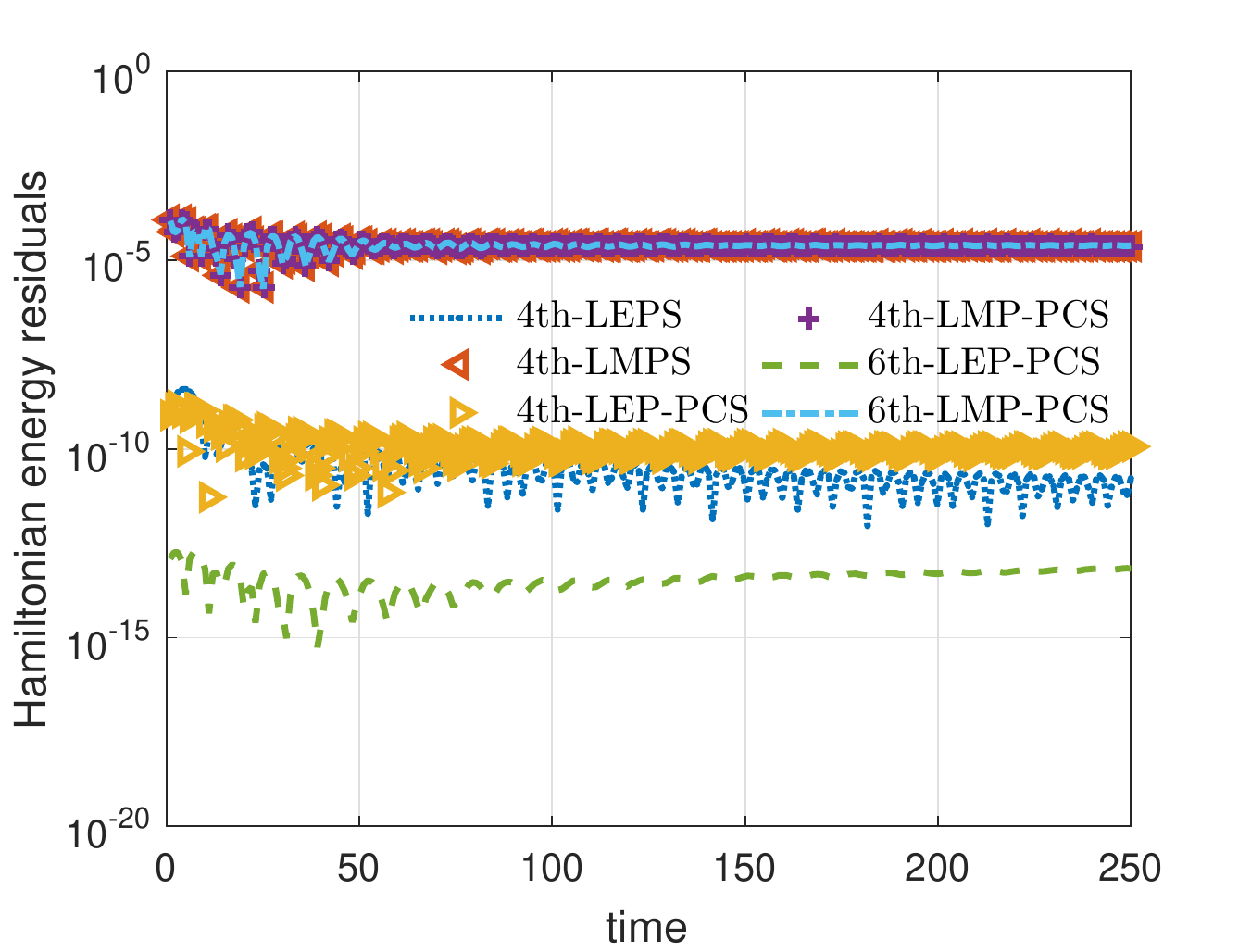}
%\caption*{(c) Hamiltonian energy}
\end{minipage}
\begin{minipage}[t]{60mm}
\includegraphics[width=60mm]{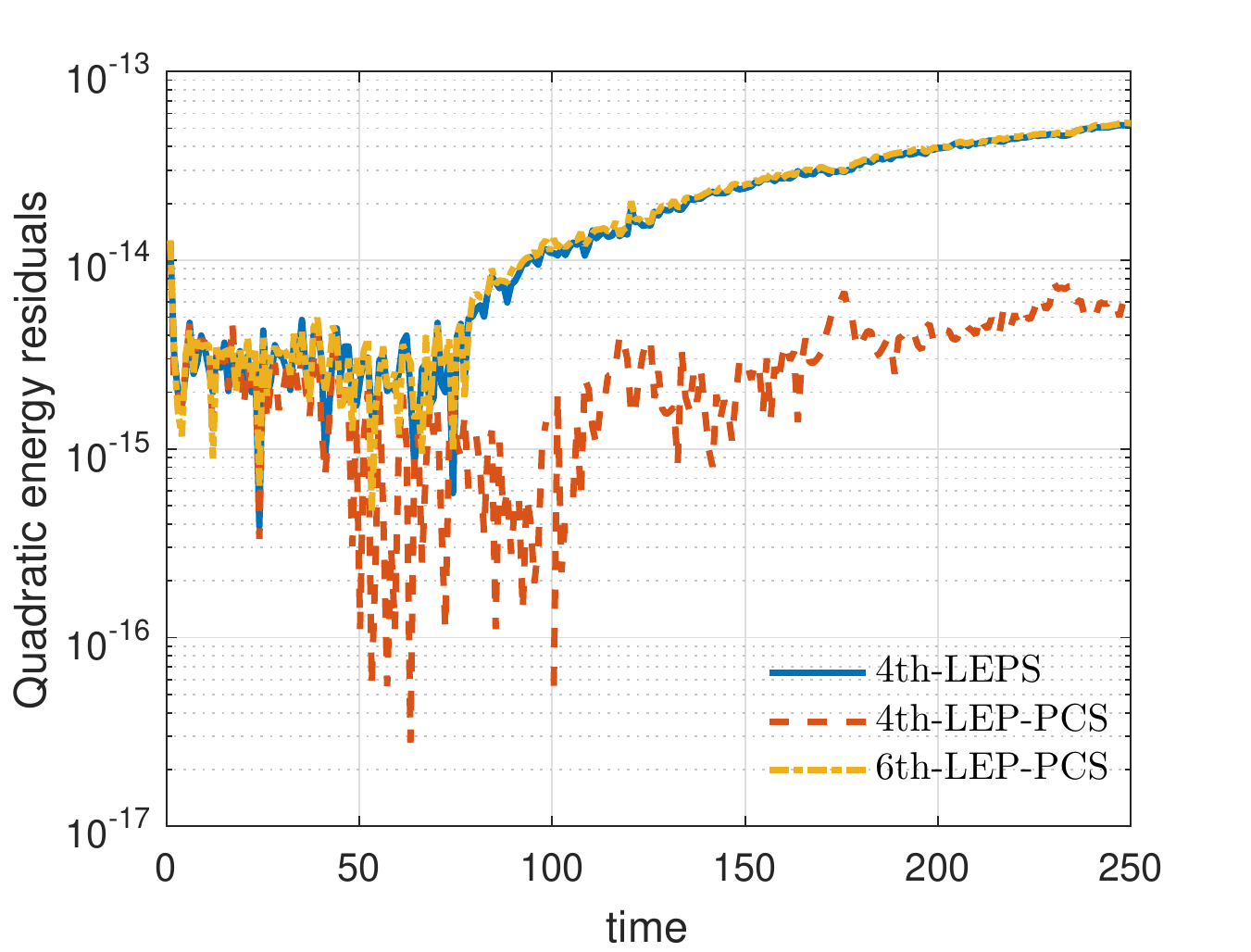}
%\caption*{(d) Quadratic energy}
\end{minipage}
\caption{ The residuals on discrete conservation laws using the different numerical schemes over time interval $[0,250]$ with the time step $\tau=0.1$ and the Fourier node $512\times 512$, respectively, for the two-dimensional RLW equation \eqref{RLW-equation}.}\label{RLW2d-example2-fig2}
\end{figure}

Finally, we consider the Maxwellian initial condition \cite{LHQjsc2018}, as follows:
\begin{align}
u(x,y,0)=e^{-((x-x_0)^2+(y-y_0)^2)},
\end{align}
with $x_0=40$ and $y_0=40$.

We take the computational domain $\Omega=[-100,100]^2$ with a periodic boundary condition and choose parameters $\alpha=1, \beta=1, \theta=1$ and $\mu=1$, and use the time step $\tau=0.1$ and the Fourier node $512\times 512$, respectively. The profiles of $u$  provided by 6th-LEP-PCS at various times are displayed in Fig. \ref{RLW2d-example3-fig1}. We can observe that the train of solitary waves expands and propagates in a certain direction.  The results are in agreement with
those given in \cite{LHQjsc2018}. Here, we  point out that the profile of $u$ at various times calculated using other
schemes are similar to Fig. \ref{RLW2d-example3-fig1}, thus for brevity, we omit them. In Fig. \ref{RLW2d-example3-fig2}, we investigate the residuals on discrete conservation laws form $t=0$ to $t=200$ of the proposed schemes, which behaves similarly as that given in Fig. \ref{RLW2d-example2-fig2}.

\begin{figure}[H]
\centering\begin{minipage}[t]{60mm}
\includegraphics[width=60mm]{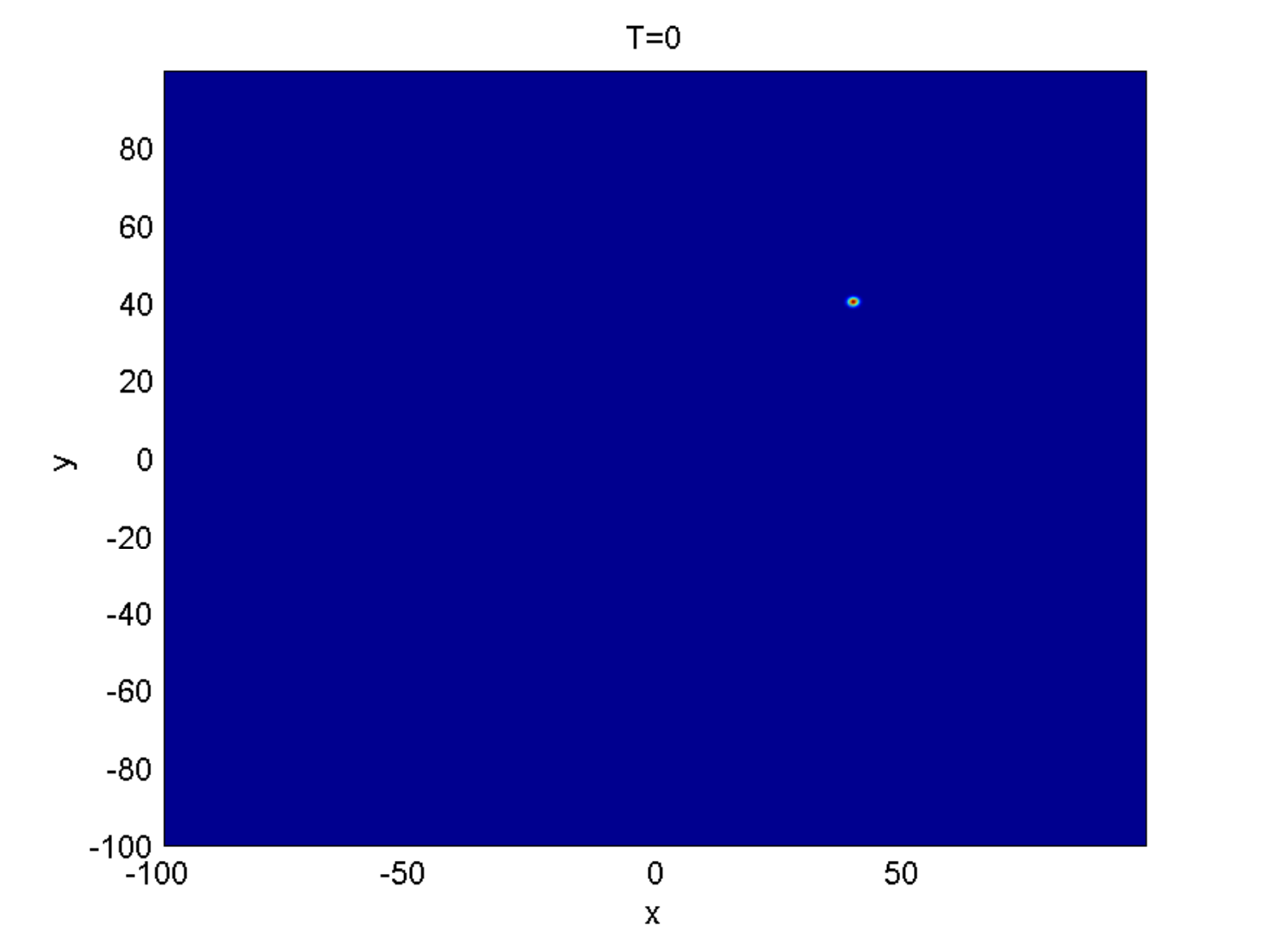}
\end{minipage}
\begin{minipage}[t]{60mm}
\includegraphics[width=60mm]{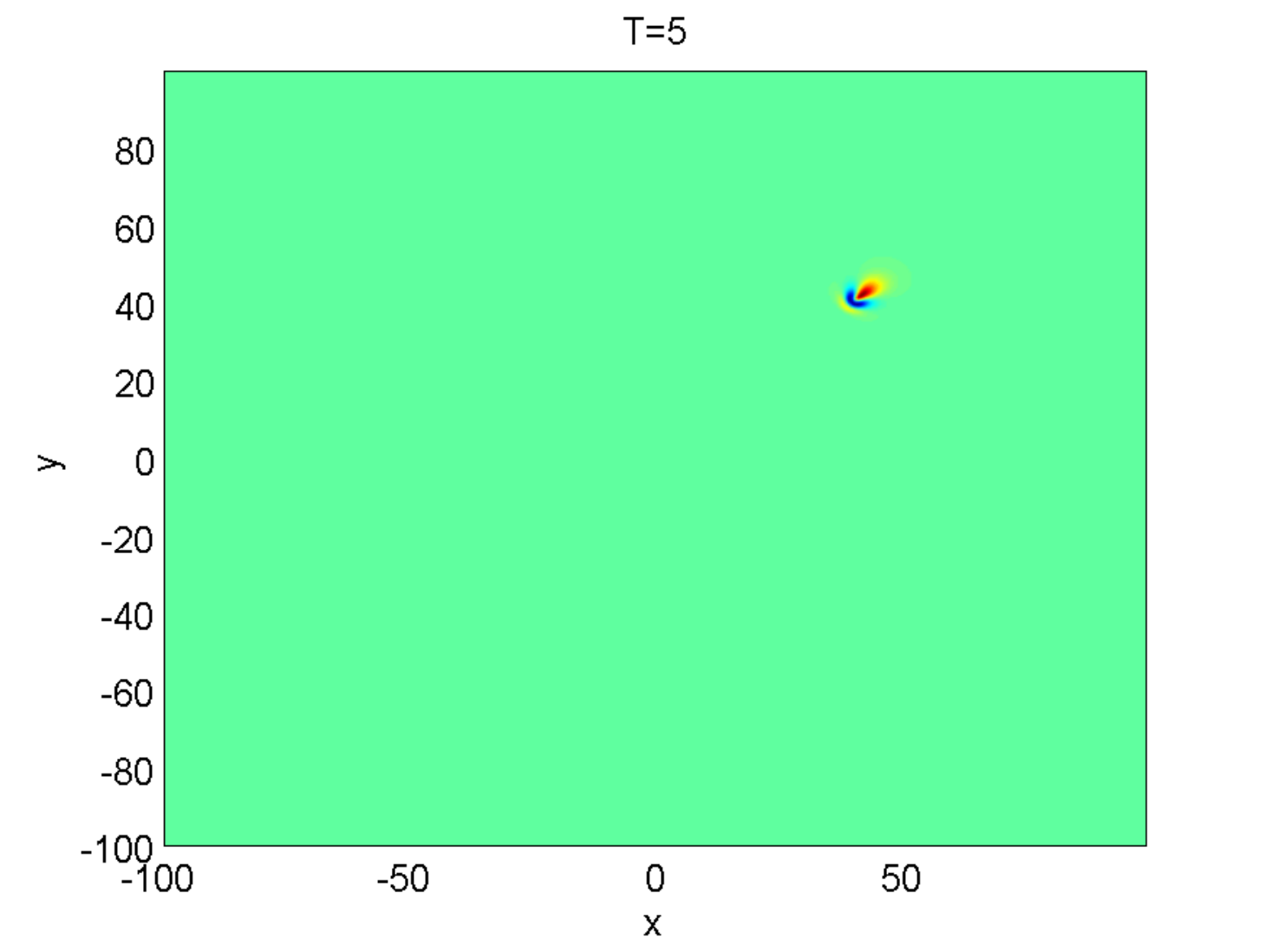}
\end{minipage}
\centering\begin{minipage}[t]{60mm}
\includegraphics[width=60mm]{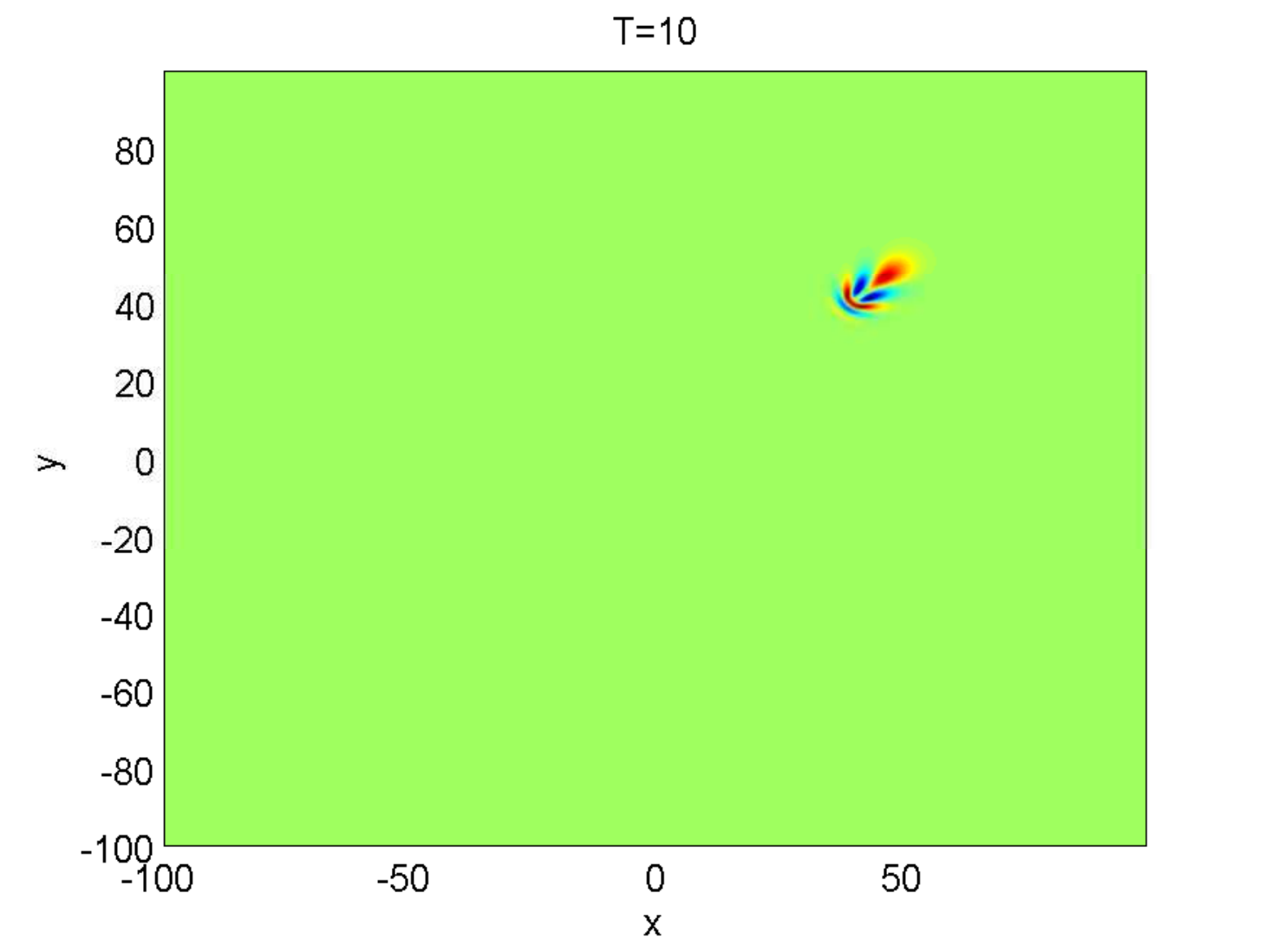}
\end{minipage}
\begin{minipage}[t]{60mm}
\includegraphics[width=60mm]{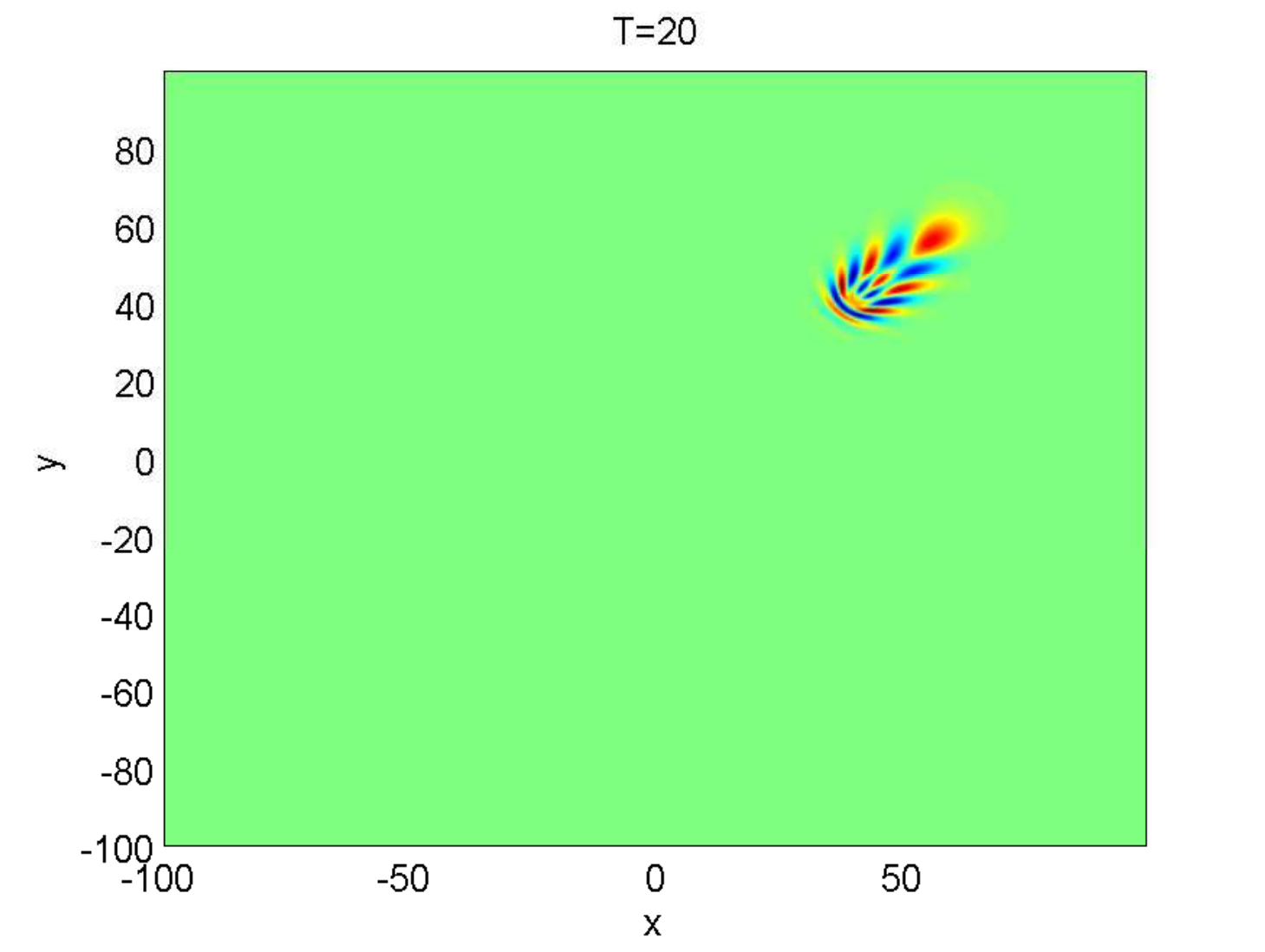}
\end{minipage}
\end{figure}
\begin{figure}[H]
\centering\begin{minipage}[t]{60mm}
\includegraphics[width=60mm]{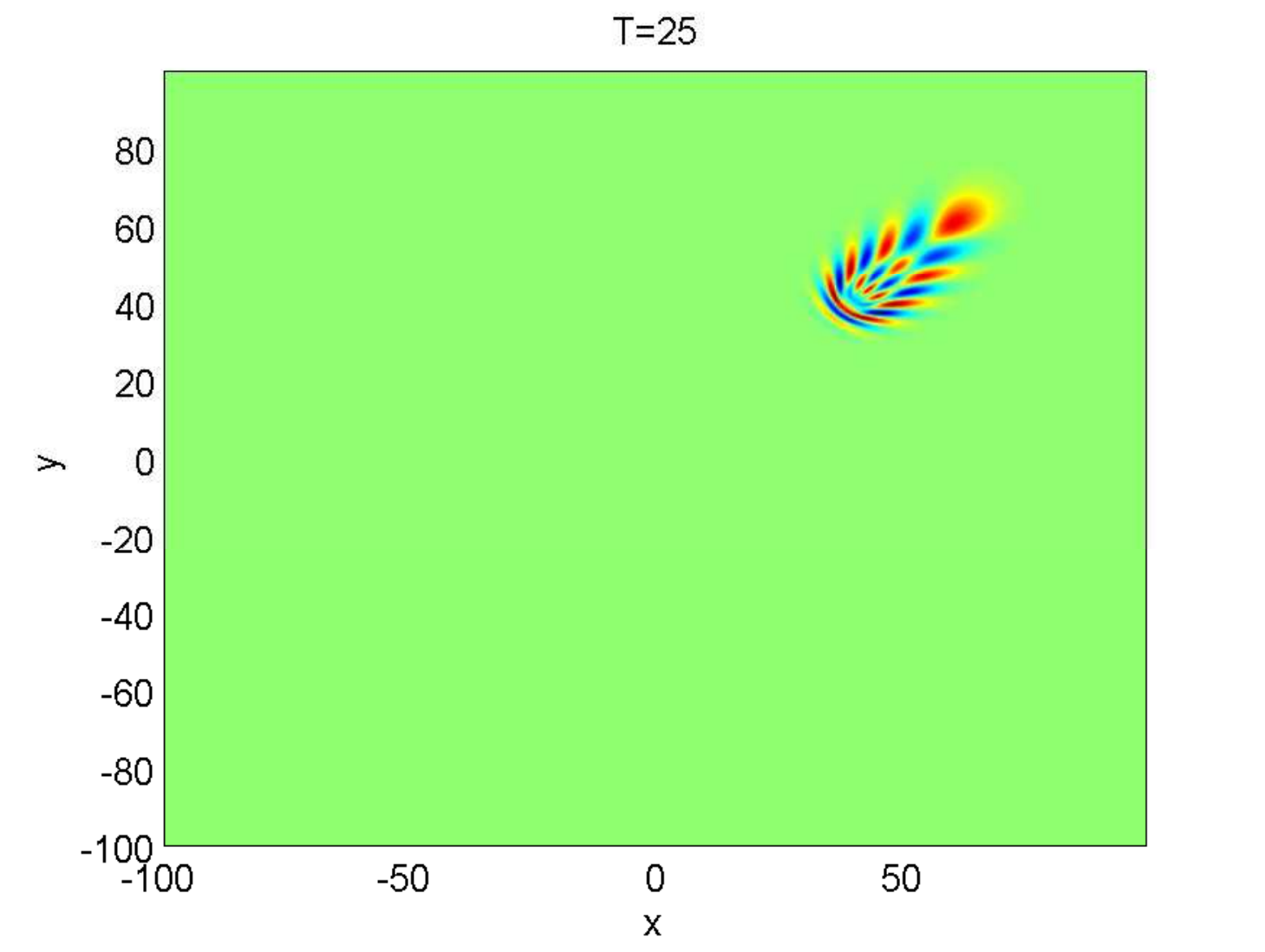}
\end{minipage}
\begin{minipage}[t]{60mm}
\includegraphics[width=60mm]{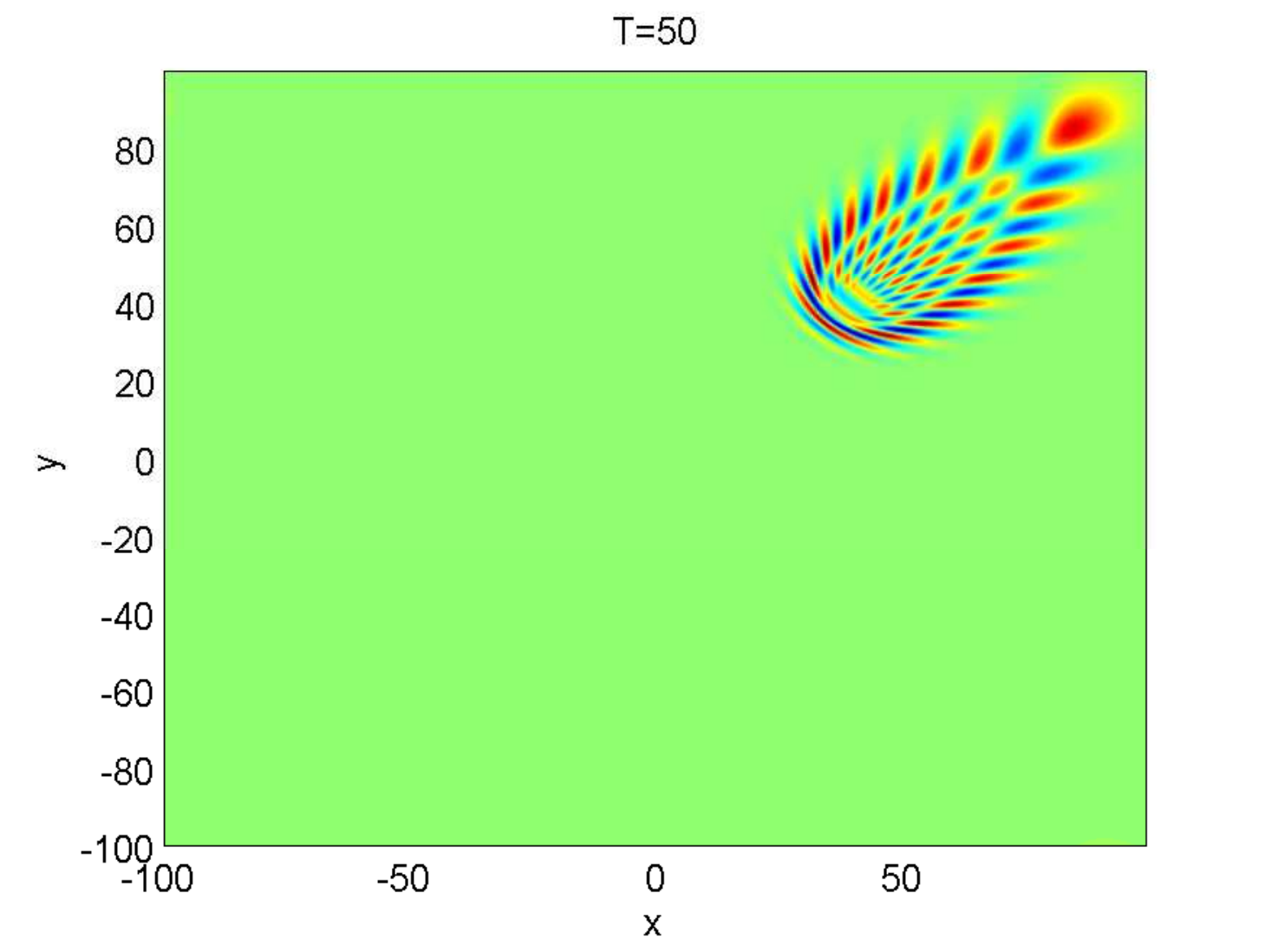}
\end{minipage}
\caption{ The profiles of $u$ provided by 6th-LEP-PCS at times $t=0,5,10,20,25$ and 50 with the time step $\tau=0.1$ and the Fourier node $512\times 512$, respectively, for the RLW equation in 2D \eqref{RLW-equation}.}\label{RLW2d-example3-fig1}
\end{figure}

\begin{figure}[H]
\centering\begin{minipage}[t]{60mm}
\includegraphics[width=60mm]{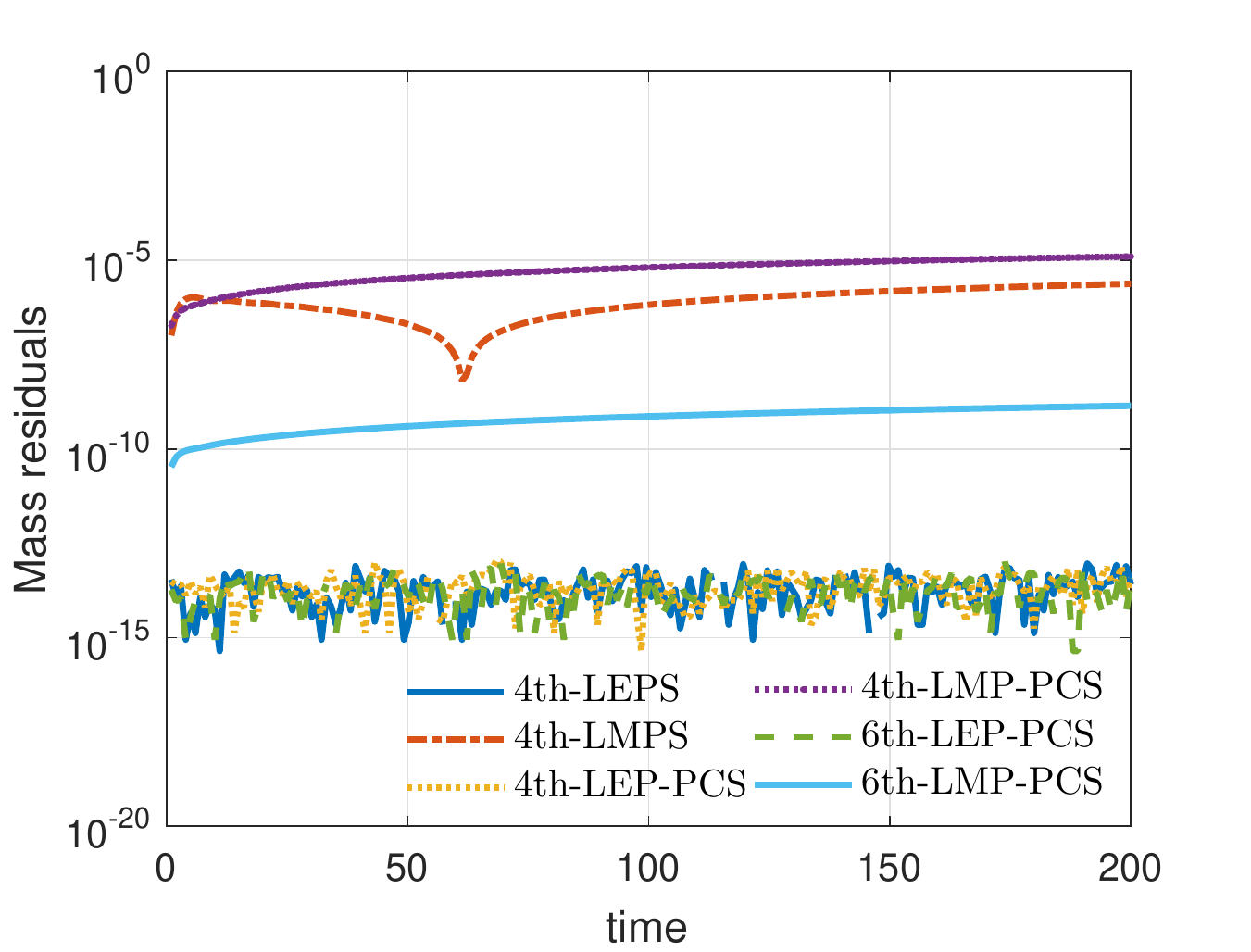}
\caption*{(a) Mass}
\end{minipage}
\begin{minipage}[t]{60mm}
\includegraphics[width=60mm]{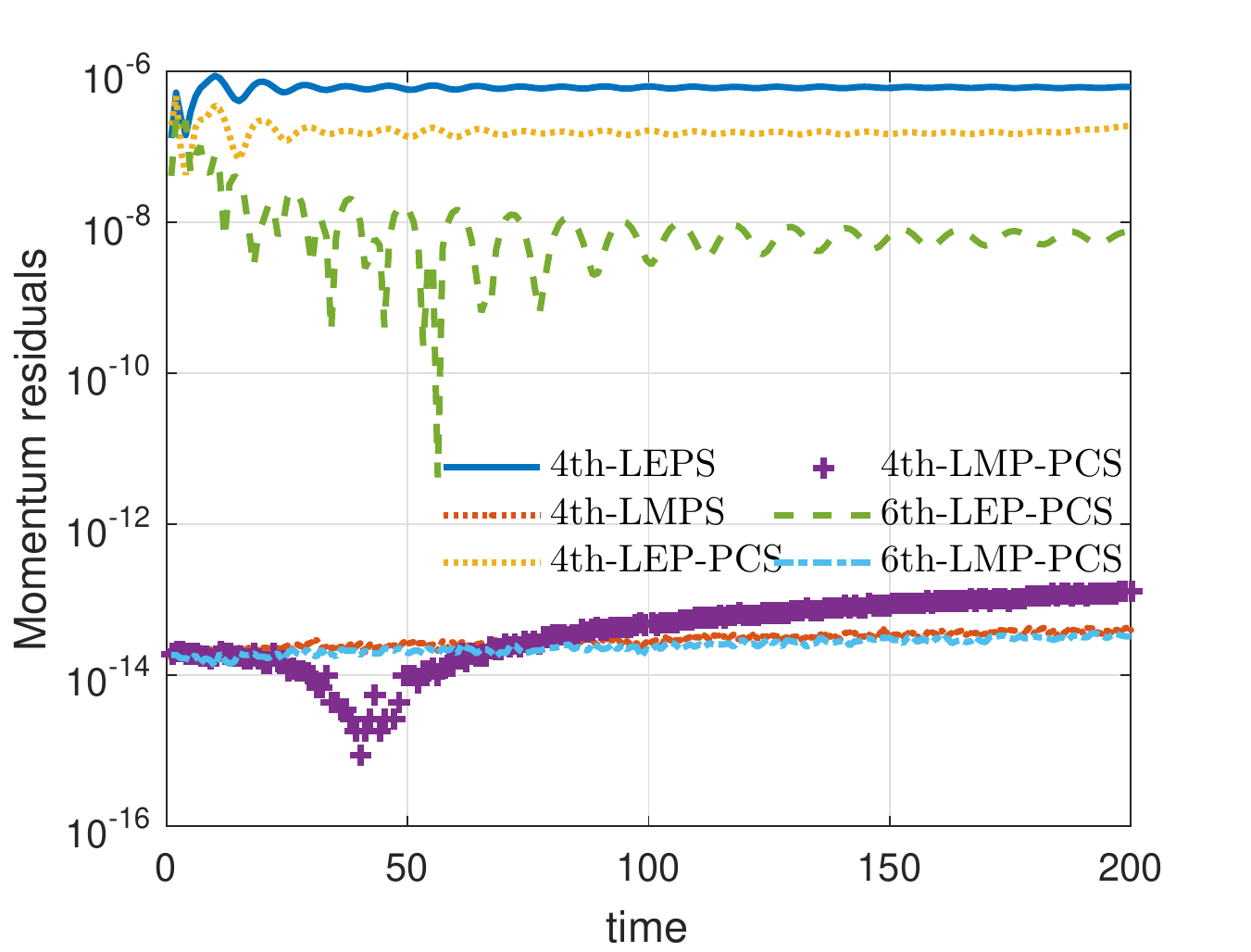}
\caption*{(b) Momentum}
\end{minipage}
\centering\begin{minipage}[t]{60mm}
\includegraphics[width=60mm]{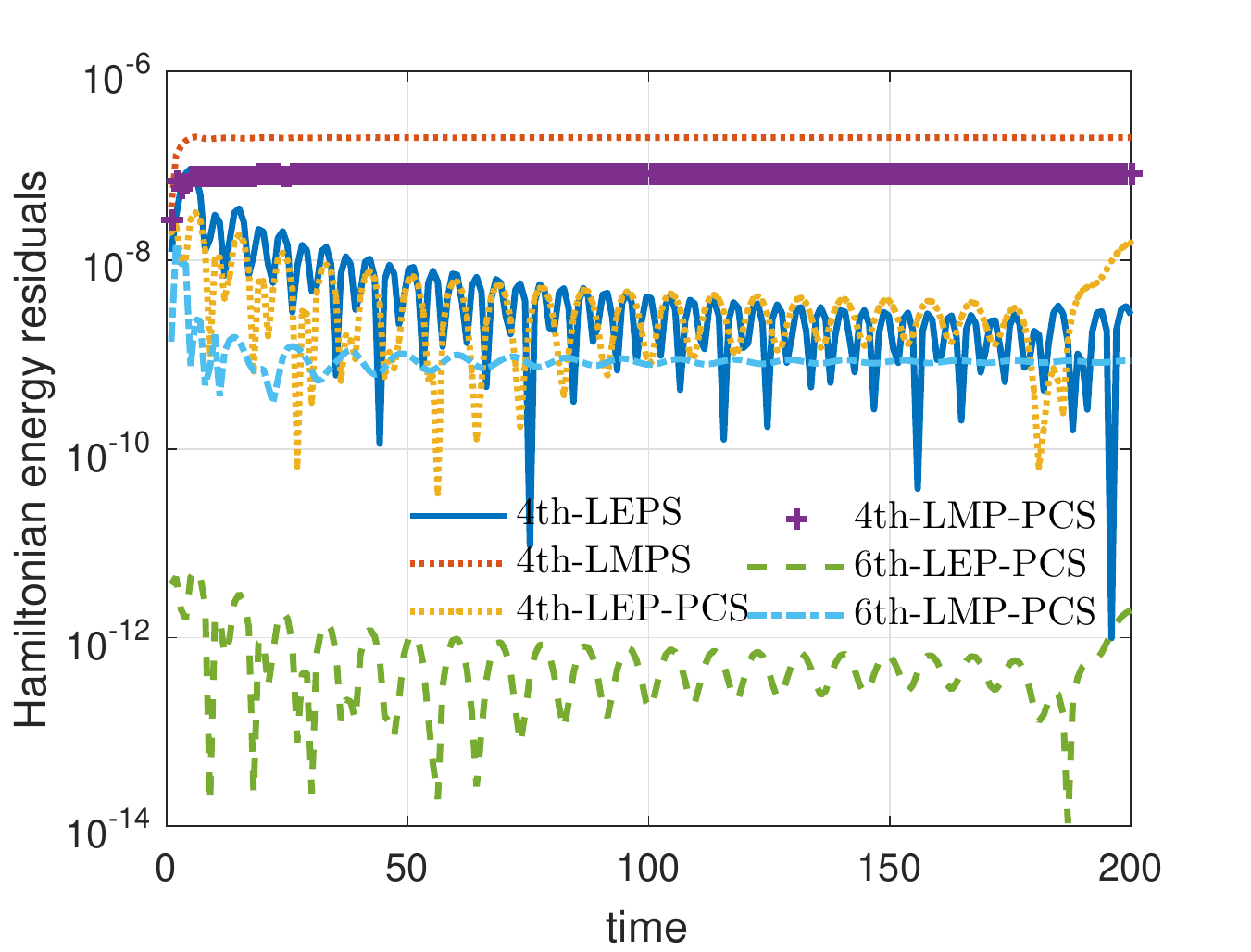}
\caption*{(c) Hamiltonian energy}
\end{minipage}
\begin{minipage}[t]{60mm}
\includegraphics[width=60mm]{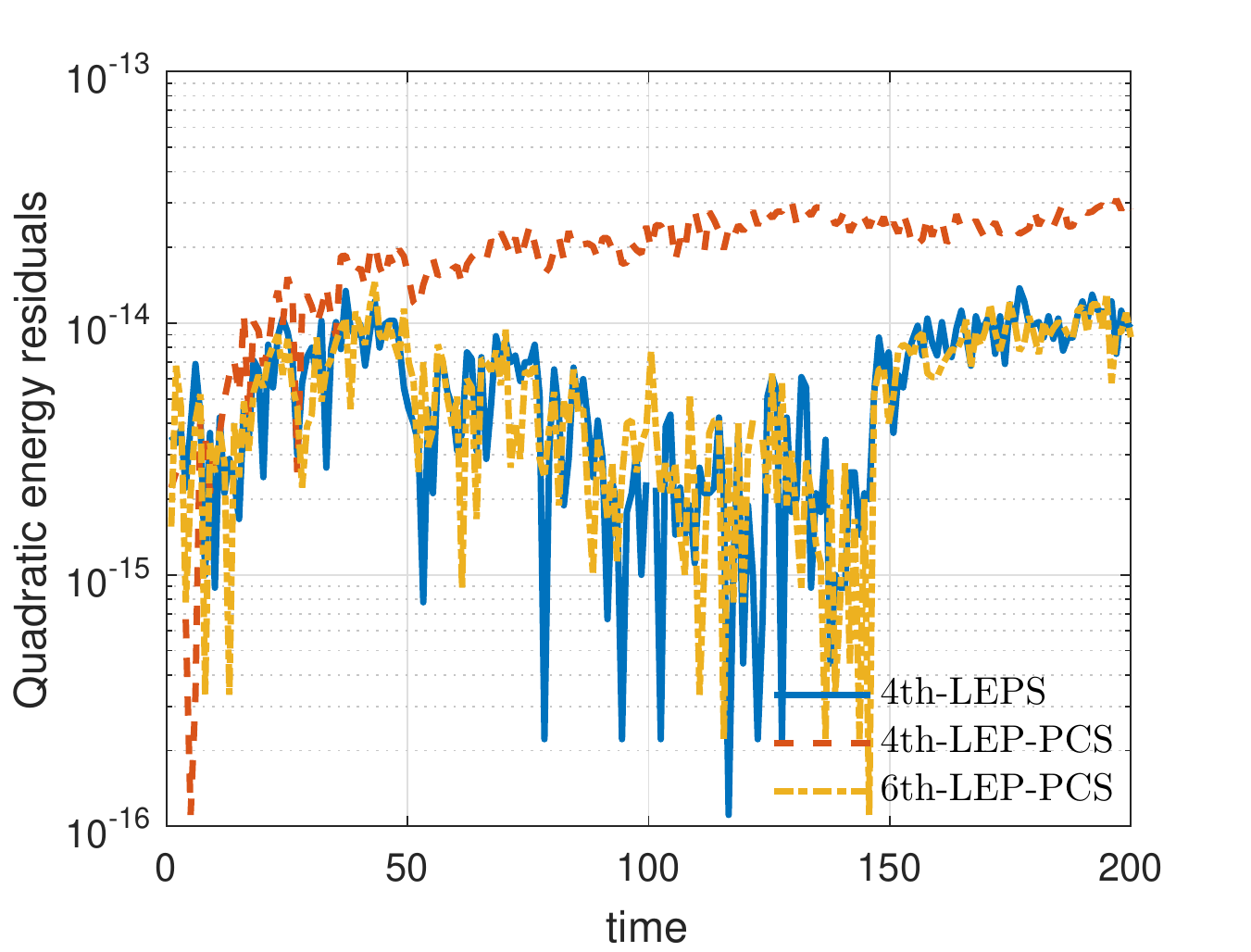}
\caption*{(d) Quadratic energy}
\end{minipage}
\caption{ The residuals on discrete conservation laws using the different numerical schemes with the time step $\tau=0.1$ and the Fourier node $512\times 512$, respectively, for the two-dimensional RLW equation \eqref{RLW-equation}.}\label{RLW2d-example3-fig2}
\end{figure}

\section{Conclusions}\label{RLW-section-4}

In this paper, two new class of structure-preserving schemes are developed for solving the RLW equation \eqref{RLW-equation}. These schemes can reach arbitrarily high-order accurate in time and are linear so that they are easy to be implementd and computationally efficient. The proposed momentum-preserving scheme is based on the use of the extrapolation/prediction-correction technique and the symplectic RK method in time, together with the standard Fourier pseudo-spectral method in space, respectively. The proposed energy-preserving scheme is based on the combination of the energy quadratization approach with the linearized idea used in the construction of the linear  momentum-preserving scheme. We show that such scheme can conserve the discrete quadratic energy. Extensive numerical examples are addressed to illustrate
the efficiencies and accuracies of our new schemes. Numerical results indicate that the momentum-preserving scheme is more robust at
larger time steps than the energy-preserving one, the energy-preserving scheme is much more efficient than the momentum-preserving scheme, and the predictor-corrector scheme performs much better than the extrapolation scheme.

We conclude this paper with second remarks. First, the proposed energy-preserving scheme involves solving linear equations with complicated variable coefficients at every time step and  this limitation can be removed by using the idea of the scalar auxiliary variable (SAV) approach \cite{CWJ-CPC-2021,SXY18,SXY17}. For more details, please refer to \cite{JCQS-arXiv-2021}.
{Second, it is known that any Runge-Kutta method conserves all linear invariants \cite{CIZ-MC-1997}, however, we should note that the proposed linear momentum-preserving schemes cannot preserve the discrete mass since the time discretization is done using the extrapolation/prediction-
correction technique. Although the higher order fully-implicit, mass and momentum conserving schemes can be obtained by applying a high-order symplectic Runge-Kutta method to the system \eqref{RLW-momentum-equation} (see \cite{Frasca-Caccia-FCM-2021}), it is challenging to construct high-order linear ones which can preserve both the mass and momentum. This will be an interesting topic for future studies.}

% Here, we should note that the proposed energy-preserving can only preserve a discrete analogue of the Hamiltonian energy, thus, it is an interesting question is whether it is possible to construct linearly implicit schemes which can preserve the discrete Hamiltonian energy, which will be an interesting topic for future studies.

%Last but not least, the presented strategy can be directly extended to propose high-order linearly implicit
%energy-preserving exponential integrators for the Hamiltonian PDEs given by
% For more details, please refer to appendix.
%However, we should note that the proposed schemes can not preserve the discrete mass conservation, Furthermore,

\section*{Acknowledgments}
The authors would like to express sincere gratitude to the referees for their insightful
comments and suggestions. This work is supported by the National Natural Science Foundation of China (Grant Nos. 11901513, 11971481, 12071481), the National Key R\&D Program of China (Grant No. 2020YFA0709800), the Natural Science Foundation of Hunan (Grant No. 2021JJ40655), the Yunnan Fundamental Research Projects (Nos. 202101AT070208, 202001AT070066,
202101AS070044), the High Level Talents Research Foundation Project of Nanjing Vocational College of Information Technology (Grant No. YB20200906), and the Foundation of Jiangsu Key Laboratory for Numerical Simulation of Large Scale Complex Systems (Grant No. 202102).
 Jiang and Cui are in particular grateful to Prof. Yushun Wang and Dr. Yuezheng Gong for fruitful discussions.

\bibliographystyle{plain}

%\bibliography{ref}

\end{document}